\documentclass[a4paper,oneside,openany,11pt]{article}
\usepackage[nointlimits]{amsmath}
\usepackage{amsfonts}
\usepackage{amssymb}
\usepackage{amsmath}
\usepackage{amsthm}
\usepackage{enumerate}
\usepackage{latexsym}
\usepackage{graphicx}
\usepackage{layout}
\newcommand{\RR}{\mathbb R}
\def \CC {\mathbb{C}}

\def \NN {\mathbb{N}}

\def \beqn {\begin{equation}}
\def \eeqn {\end{equation}}
\def \bi {\begin{itemize}}
\def \ei {\end{itemize}}
\def \ben {\begin{enumerate}}
\def \een {\end{enumerate}}
\def \beq {\begin{eqnarray*}}
\def \eeq {\end{eqnarray*}}
\def \beqn {\begin{eqnarray}}
\def \eeqn {\end{eqnarray}}
\def \ds {\displaystyle}

\def\pats{\left(}
\def\patd{\right)}

\def\coz{z}

\def\mz2a{\pabs \coz\pabd^{2\al}}
\def\coza{z^{\al+1}}

\def\pcz{\pats \coz\patd}

\def\pzn{\pats \coz_{n}\patd}

\def\pabs{|}
\def\pabd{|}

\def\pnsz{\psi_{n}\pcz}
\def\pnnsz{\psi_{n}\pzn}
\def\psz{\psi\pcz}

\def\pnszp{\psi_{n}^{'}\pcz}
\def\pnnszp{\psi_{n}^{'}\pzn}
\def\pszp{\psi^{'}\pcz}

\def\numza_2{\pabs\unoalfa\psz+z\pszp \pabd^{2}}
\def\numnza_2{\pabs\pats 1+\al_{n}\patd\pnsz+z\pnszp \pabd^{2}}

\def\numnnza_2{\pabs\pats 1+\al_{n}\patd\pnnsz+\coz_{n}\pnnszp \pabd^{2}}

\def\numzb{\pabs-\unoalfa\psz+z\pszp \pabd}
\def\numzb_2{\pabs-\unoalfa\psz+z\pszp \pabd^{2}}
\def\numnzb_2{\pabs-\pats 1+\al_{n}\patd\pnsz+z\pnszp \pabd^{2}}

\def\numzc{\pabs e^{\ds\coza\psz}\pats\unoalfa\psz+z\pszp \patd\pabd}
\def\numzc_2{\pabs e^{\ds\coza\psz}\pats\unoalfa\psz+z\pszp \patd\pabd^{2}}
\def\numnzc_2{\pabs e^{\ds\coza\pnsz}\pats\unoalfa\pnsz+z\pnszp \patd\pabd^{2}}

\newcommand{\ba}{\begin{array}}
\newcommand{\ea}{\end{array}}
\newcommand{\beqa}{\begin{eqnarray}}
\newcommand{\eeqa}{\end{eqnarray}}

\newcommand{\dimo}{\smallskip\noindent{\bf Proof:\ }}

\newtheorem{thm}{Theorem}[section]
\newtheorem{lemma}[thm]{Lemma}
\newtheorem{prop}[thm]{Proposition}
\newtheorem{coro}[thm]{Corollary}
\newtheorem{defn}[thm]{Definition}

\newtheorem{rem}[thm]{Remark}

\numberwithin{equation}{section}

\begin{document}

\title{Blow-up analysis for a cosmic strings equation}
\author{Gabriella Tarantello \footnote{Dipartimento di Matematica.
Universit\`a degli Studi di Roma "Tor Vergata", Via della Ricerca Scientifica, 00133 Rome, Italy. Email: tarantel@mat.uniroma2.it}
\thanks{Supported by PRIN09 project: \textit{Nonlinear elliptic problems in the study of vortices and related topics},  PRIN12 project: \textit{Variational and Perturbative Aspects of Nonlinear Differential Problems} and  FIRB project: \textit{Analysis and Beyond.}}}
\maketitle
\section{Abstract}
\setcounter{equation}{0}
In this paper we develop a blow-up analysis for solutions of an elliptic PDE of Liouville type over the plane. Such solutions describe the behavior of cosmic strings (parallel in a given direction) for a W-boson model coupled with Einstein's equation.
We show how the blow-up behavior of the solutions is characterized, according to the physical parameters involved, by new and surprising phenomena. For example in some cases, after a suitable scaling, the blow-up profile of the solution is described in terms of an equations that bares a geometrical meaning in the context of  the "uniformization"  of the Riemann sphere with conical singularities.

\section{Introduction}
\setcounter{equation}{0}
In this paper, we focus on the blow-up analysis for soluions of a Liouville type equation describing the behavior of  selfgravitating cosmic strings for a massive W-boson model coupled with Einstein's equation, introduced in \cite{ao1}. More precisely, for the model in \cite{ao1}, it is possible to use Einstein's equation together with a suitable ansatz, in order to reduce the analysis of the corresponding selfgravitating cosmic string located at the origin and parallel to the $x_3$-direction to the study of the following elliptic problem:
\begin{equation}\label{prob}
\begin{cases}
-\Delta u=e^{au}+\left|x\right|^{2N}e^{u}\quad\mbox{in }\RR^{2},\\
\displaystyle{\int_{\RR^{2}} \left(e^{au}+\left|x\right|^{2N}e^{u}\right)}<\infty,
\end{cases}
\end{equation}
with $a>0$ a physical parameter and $N \in \mathbb{N}$ the string's multiplicity, see \cite {ao1},  \cite{y} and \cite{pot} for details. 

For a solution $u$ of \eqref{prob} the value, 
\begin{equation}\label{prob1bis}
\beta:=\dfrac{1}{2\pi}\displaystyle{\int_{\RR^{2}} \left(e^{au}+\left|x\right|^{2N}e^{u}\right)}
\end{equation}
relates to the (finite) string's energy, and our main concern will be to identify the (sharp) range of $\beta$' s for which problem \eqref{prob}-\eqref{prob1bis} is solvable. We mention that some existence results  concerning \eqref{prob} are contained in \cite{c}, \cite{c1}, \cite{y} and \cite{cgs}.

Here, we shall build our investigation upon the work in \cite{pot,pt2}, where the authors characterize completely the radial solvability of \eqref{prob}-\eqref{prob1bis}.\\ To be more specific we observe that, (as shown in \cite{cgs} and \cite{cl}), every solution $u$ of \eqref{prob} and \eqref{prob1bis} satisfies:
\begin{equation}\label{probd}
-\beta\log\left(\left|x\right|+1\right)-C\leq u(x)\leq-\beta\log\left(\left|x\right|+1\right)+C,
\end{equation}
\begin{equation}\label{probd2}
r \partial_r u \to-\beta\quad\mbox{ , }\quad \partial_{\vartheta}u\to 0\quad\mbox{ as }r\to+\infty,
\end{equation}
with suitable $C>0$ (depending on $u$), and $(r, \vartheta)$ the polar coordinates in $\RR^{2}$. Thus, from \eqref{probd} we find that,
\begin{equation}\label{001bis}
\beta=\dfrac{1}{2\pi}\displaystyle{\int_{\RR^{2}} \left(e^{au}+\left|x\right|^{2N}e^{u}\right)}>\max\left\{\dfrac{2}{a}, 2(N+1)\right\}.
\end{equation} 
The two values on the right hand side of \eqref{001bis} coincide for $a=\dfrac{1}{N+1}$, and in this case , the equation in \eqref{prob} acquires the following scaling invariance property:
\begin{equation}\label{eq1.6}
u(x)\to u_{\lambda}(x):=u(\lambda x)+2(N+1)\log\lambda
=u(\lambda x)+\dfrac{2}{a}\log\lambda,
\end{equation}
for every $\lambda>0$.

In turn, as shown in \cite{cgs}, one can use a Pohozaev's type inequality (in the usual way) to find that, for problem \eqref{prob} the following holds:
\begin{equation}
\text{ if } a=\dfrac{1}{N+1} \text { then } \dfrac{1}{2\pi}\displaystyle{\int_{\RR^{2}} \left(e^{au}+\left|x\right|^{2N}e^{u}\right)}=4(N+1)\left(=\dfrac{4}{a}\right).
\end{equation}

Furthermore, for $a=\dfrac{1}{N+1}$, problem \eqref{prob}  gains an additional invariance with respect to involution, and such an invariance property is inherited by the corresponding solution. This fact is shown in  \cite{cgs} ( after the approach of  \cite{pt}) together with other relevant facts about the solutions of \eqref{prob}. 

To justify the nature of the values $\beta=4(N+1)$ and $\beta=\dfrac{4}{a}$, we recall that for the Liouville type equation: 
\begin{equation}\label{lspeintro}
\begin{cases}
-\Delta u=\left|x\right|^{2N}e^{bu}\quad\mbox{in }\RR^{2},\\
\displaystyle{\int_{\RR^{2}}\left|x\right|^{2N}e^{bu}}<+\infty,
\end{cases}
\end{equation}
with $N>-1$ and $b>0$, all solutions are explicitly identified and they satisfy: 
$$
\displaystyle{\dfrac{1}{2\pi}\int_{\RR^{2}}\left|x\right|^{2N}e^{bu}=\dfrac{4(N+1)}{b}},
$$
see \cite{cl,cl1,pt,tar}. In particular:
\begin{equation}\label{lspeintrovvv}
\mbox{if } b=1, \mbox{ then }\quad\dfrac{1}{2\pi}\int_{\RR^{2}} \left|x\right|^{2N}e^{u}=4(N+1),
\end{equation}

while,
\begin{equation}\label{lspeintrovvvns}
\mbox{if } b=a \mbox{ and } N=0, \mbox{ then }\quad \dfrac{1}{2\pi}\int_{\RR^{2}} e^{au}=\dfrac{4}{a}.
\end{equation}

In fact, when $0<a \neq \dfrac{1}{N+1}$, then, under either one of the scaling indicated in \eqref{eq1.6}, problem \eqref{prob} can be interpreted as a 'perturbation' of \eqref{lspeintro} respectively with $b=a$ and $N=0$, or with $b=1$.

Also notice that for $a=1$ problem \eqref{prob} reduces to a much studied equation in the context of the assigned Gauss curvature problem in $\RR^{2}$. In this respect, many interesting results have been established especially in the radial case, see for example \cite{lin},  \cite{chlin}, \cite{ck1},  \cite{det},  \cite{blt} and references therein.

Therefore, for problem \eqref{prob}, the true interesting and novel case to investigate occurs when 
\begin{equation}\label{001star1}
0<a\neq\dfrac{1}{N+1}\quad\mbox{ and }\quad a\neq 1.
\end{equation}

To this purpose, we observe that if $u$ saisfies \eqref{prob}- \eqref{prob1bis}, then by virtue of \eqref{probd}, it is possible to derive the following Pohozaev's identity:

\begin{equation}\label{pohov}
2N\int_{\RR^2}\left|x\right|^{2N}e^u dx+2\left(\dfrac{1}{a}-1\right)\int_{\RR^2}e^{au}{dx}=\pi\beta(\beta-4),
\end{equation}

see \cite{cgs} for details, and as a consequence we obtain:

\begin{equation}\label{lim}
\dfrac{1}{2\pi}\int_{\RR^2}\left|x\right|^{2N}e^{u}=\dfrac{\beta(4-a\beta)}{4[1-a(N+1)]}\quad\mbox{and}\quad
\dfrac{1}{2\pi}\int_{\RR^2}e^{au}=\dfrac{a\beta(\beta-4(N+1))}{4[1-a(N+1)]}.
\end{equation}

Therefore, from \eqref{lim} we find the following additional necessary condition for the solvability of \eqref{prob}--\eqref{prob1bis}:
\begin{equation}\label{prob2}
\beta\in\left(\min\left\{4/a, 4(N+1)\right\},\max\left\{4/a, 4(N+1)\right\}\right).
\end{equation}

At this point, it is natural to ask whether \eqref{001bis} and \eqref{prob2} are also sufficient for the solvability of \eqref{prob}.

A complete answer to this question has been provided by Poliakovsky and Tarantello \cite{pot} in the context of radially symmetric solutions. See also \cite{pt2} where the same authors deal with a general class of 'cooperative' systems of Liouville type which include \eqref{prob} as a degenerate case. In the context of Liouville systems, we mention also the related results contained in \cite{csw}, \cite{ck1}, \cite{ck}, \cite{lwy}, \cite{lz}, \cite{jost}, \cite{sw1} and \cite{sw2}, which have motivated the work in  \cite{pt2}. \\  From \cite{pot} and \cite{pt2} we know the following:

\begin{thm}[Theorem~1.1, \cite{pot}, Theorem~5.1, \cite{pt2}]\label{jjj7} Let $N>-1$ and $0<a\neq\dfrac{1}{N+1}$, then problem \eqref{prob}--\eqref{prob1bis} admits a radial solution if and only if: 
\begin{equation}\label{servet1}
\beta\in\left(\max\left\{4(N+1),\,\dfrac{4}{a}-4(N+1)\right\},\dfrac{4}{a}\right)\quad\mbox{when}\quad 0<a<\dfrac{1}{N+1},
\end{equation}

and,

\begin{equation}\label{servet2}
\beta\in\left(\max\left\{\dfrac{4}{a},\,4(N+1)-\dfrac{4}{a}\right\}, 4(N+1)\right)\quad\mbox{when }\quad a>\dfrac{1}{N+1}.
\end{equation}

Furthermore for any $\beta$ satisfying \eqref{servet1} or \eqref{servet2} the corresponding radial solution of \eqref{prob}--\eqref{prob1bis} is {\bf{unique}} and {\bf{nondegenerate}}. 
\end{thm}
\qed

By comparing \eqref{servet1} and \eqref{servet2} with \eqref{prob2}, we see that, 
\begin{equation}
\text{ if }\dfrac{1}{2(N+1)}\leq{a}\neq{\dfrac{1}{N+1}}\leq{\dfrac{2}{N+1}}\text{ then }\eqref{prob}-\eqref{prob1bis}\mbox{ is solvable iff } \eqref{prob2}\mbox{ holds.} 
\end{equation}
Furthermore, for $N\in (-1,0]$ all solutions of \eqref{prob} are necessarily radially symmetric (see Theorem~1.4.1 in \cite{f}), and  therefore Theorem~\ref{jjj7} also provides a complete answer to the question of solvability of  \eqref{prob}-\eqref{prob1bis} in this case.\\

On the other hand for $N>0$, we know from \cite{pot} that  \eqref{prob} admits also non-radial solutions. Thus, it becomes an interesting issue to investigate whether or not the existence interval specified in \eqref{servet1} or \eqref{servet2} remains 'sharp', also for $0<a<\dfrac{1}{2(N+1)}$ or $a>\dfrac{2}{N+1}$ and when we take into account also (possible) {\bf{non radial}} solutions. \\To contribute in this direction,  we develop a general blow--up analysis for problem \eqref{prob}, with the goal to detect (as limits) the sharp bounds  on the $\beta$'s which allow for the solvability of problem \eqref{prob}--\eqref{prob1bis}, beyond the radial situation.

In this respect, we recall that the (existence) intervals specified in \eqref{servet1} and \eqref{servet2}, have been identified in \cite{pot} via a blow--up argument applied to solutions of the radial 'Initial Value Problem' corresponding to equation \eqref{prob}, when the initial datum tends to $+\infty$ (blow--up) or to $-\infty$ (blow--down), see \cite{pot} for details. 

Furthermore, we hope that the understanding of the blow--up behaviour of solutions to \eqref{prob} may help (via a degree argument) to find more general strings located at different points other than those treated here, which are just superimposed at the origin. \\

Thus, the main focus of this paper will be to analyse a solution sequence $\left\{u_k\right\}$ satisfying:

\begin{equation}\label{probprimo}
\begin{cases}
-\Delta u_k=e^{au_k}+\left|x\right|^{2N}e^{u_k}\quad\mbox{in }\RR^{2},\\
\beta_k=\dfrac{1}{2\pi}\displaystyle{\int_{\RR^{2}}\left(e^{au_k}+\left|x\right|^{2N}e^{u_k}\right)dx},
\end{cases}
\end{equation}
such that, 
\begin{equation}\label{0010star}
\beta_k\to \beta^{\infty}\quad\mbox{ as }\quad k\to+\infty.
\end{equation}

Furthermore, to express the fact that we cannot pass to the limit along $u_k$,  we  assume the following, 

\begin{equation} \label{002bis}
\sup_{|x| < R}u_k \to+\infty\quad \text {or}\quad \sup_{|x| \geq R} \left(u_{k}+\beta_k\log\left|x\right|\right)\to+\infty,\qquad\mbox{as }k\to+\infty.
\end{equation}
for some $R>0$. \\
Indeed, we shall see that condition \eqref{002bis} does imply that "blow-up" or "vanishing"  occurs for  $u_k$.

We aim to prove that this situation can occur only for specific "limiting" values of $\beta^{\infty}$ in \eqref{0010star}, which  should give an indication about the range of $\beta$ in \eqref{prob1bis} for which \eqref{prob} is solvable, or else about the nature of  possible non-radial solutions. \\
Notice that, from \eqref{001bis} and \eqref{prob2}, we know already that $\beta^{\infty}$ must satisfy:
\begin{equation}\label{betakkkkkkkk}
\beta^{\infty}\geq\max\left\{\dfrac{2}{a},2(N+1)\right\}
\end{equation}
and
\begin{equation}\label{betakkkkkkkk2}
\min\left\{\dfrac{4}{a}, 4(N+1)\right\}\leq\beta^{\infty}\leq\max\left\{\dfrac{4}{a}, 4(N+1)\right\}.
\end{equation}

Clearly, the inequality \eqref{betakkkkkkkk} confirms the fact that we need to distinguish between the cases: 

\begin{equation}\label{casia}
0<a<\dfrac{1}{N+1}\quad\mbox{ and }\quad a>\dfrac{1}{N+1}.
\end{equation}

Since the lower bounds of $\beta^{\infty}$ given in \eqref{betakkkkkkkk} and \eqref{betakkkkkkkk2} coincide exactly when
$$
a=\dfrac{1}{2(N+1)}\quad\mbox{ and }\quad a=\dfrac{2}{N+1}
$$
we can anticipate that also those values will  play a role in the bolw-up analysis discussed below.\\

On the basis of a suitable Harnack type inequality, we obtain the following property which describes a typical behaviour for solutions of Liouville-type equations subject to blow-up. 

\begin{prop}\label{prop1.2}
Let $N>-1$ and assume \eqref{001star1}. If $\left\{u_k\right\}$ is a sequence of solutions satisfying \eqref{probprimo}, \eqref{0010star} and \eqref{002bis}, then\\

either (blow-up)

\begin{equation}\label{01}
\exists \; R_0 > 0  :  \sup_{B_{R_0}} u_k \to +\infty  \text{ as } k \to \infty, 
\end{equation}

or (vanishing)

\begin{equation}\label{02}
\sup_{B_{R}} u_k \to -\infty \text{ as } k \to \infty, \quad \forall \; R > 0.
\end{equation}

\end{prop}
\qed

We start to discuss the case: $0<a<\dfrac{1}{N+1}$, where we provide rather complete results as follows:

\begin{thm}\label{teo1introA}
Let $N>0$ and $0<a<\dfrac{1}{N+1}$. If  $\left\{u_k\right\}$ satisfies \eqref{probprimo}, \eqref{0010star} and if \eqref{01} holds then, 

\begin{equation}\label{eq010}
\beta^\infty = \max \left\{ 4(N+1), \frac{4}{a} - 4(N+1) \right\}.
\end{equation}

\end{thm}
\qed

Clearly, \eqref{eq010} is consistent with \eqref{servet1}.\\

On the other hand, when \eqref{02} holds we obtain:

\begin{thm}\label{teo1introB}
Let $N>0$ and $0<a<\dfrac{1}{N+1}$. Assume that $\left\{u_k\right\}$ satisfies \eqref{probprimo}, \eqref{0010star}, \eqref{002bis} and \eqref{02}, then we have:

\begin{enumerate}[i)]

\item If $0 < N < 1$ then $\beta^\infty = \frac{4}{a}$.

\item If $N \geq 1$ then either $\beta^\infty = \frac{4}{a}$ or  $\beta^\infty$ takes either  one of the following values:

\begin{equation}\label{03}
\beta^\infty = \frac{2}{a} \left( 1 + \sqrt{1-4am(1-a(N+1))} \right), \text{ with } m \in \mathbb{N}: 2 \leq m < N+1 \text{ or } m=N+1,
\end{equation}

\begin{equation}\label{04}
\begin{split}
\beta^\infty = &\frac{2}{a} \sqrt{\left( 1 - \frac{2(m-1)(1-a(N+1))}{N} \right)^2 + \frac{4(m-1)ma}{N} (1-a(N+1))}\\
& + \frac{2}{a} \left( 1 - \frac{2(m-1)(1-a(N+1))}{N} \right) 
\end{split}
\end{equation}

with $m \in \mathbb{N}: 2 \leq m \leq  1+ \frac{1}{2a} \left( \sqrt{(1-a(N+1))^2 + \frac{Na}{1-a}}-(1-a(N+1)\right)$ and in this case $a > a_N$ (defined in \eqref{eq532} below).

\end{enumerate}
\end{thm}
\qed

We refer to Theorem \ref{teo5.6} for a more detailed statement.

\begin{rem}\label{rem0.1}
We suspect that it should be possible to rule out alternative \eqref{04}, in the sense that it can occur only for $m=1$ in order to account for the value: $\beta^\infty = \frac{4}{a}$.\\ 
On the contrary, \eqref{03} seems very likely to occur for the remaining values of $ m=2,...., N+1$. Indeed the expression of $ \beta^\infty $ in \eqref{03} is smaller than $\frac{4}{a}$ and it is decreasing with respect to $m$  with least value at $\beta^\infty = \max \left\{ 4(N+1)\frac{4}{a} - 4(N+1) \right\}$, attained for $m=N+1$. \\ Even more interestingly, we point out that, \eqref{03} relates the blow--up behaviour of $u_k$ to a ''singular'' Liouville equation in the plane (see \eqref{eq566}) that bears a geometrical meaning in the context of the uniformization of Riemann surfaces with conical singularities, see \cite{sg}. In other words, with \eqref{03} we express an interesting connection between a possible blow--up profile of $u_k$ and the existence of a conformal metric on the Riemann sphere with constant curvature equal to one and assigned conical singularities.
For a detailed discussion of this aspect, we refer to the proof of Theorem \ref{teo5.6} in Sction 5. By the well known uniqueness and non degeneracy properties of solutions of   \eqref{eq566}, see \cite{lt}, and the advanced "perturbation" techniques available in literature (see e.g. \cite{bp}, \cite{c}, \cite{ct}, \cite{ct1}, \cite{egp}and  \cite{pkm}) it should be possible to exhibit explicit examples of solution sequences of \eqref{probprimo}, \eqref{0010star} which satisfies \eqref{002bis} and whose blow-up behaviour is characterised by \eqref{03}.
\end{rem}

In case $a>\dfrac{1}{N+1}$ we prove the following:

\begin{thm}\label{teo2introA}
Let $N>0$ and $\frac{1}{N+1} < a < 1$. If  $\left\{u_k\right\}$ satisfies \eqref{probprimo}, \eqref{0010star} and \eqref{002bis} and if \eqref{02} holds, then we have:  

\begin{equation}\label{05}
\beta^\infty = 4(N+1).
\end{equation}

\end{thm}
\qed

Alternatively, when blow-up occurs, in the sense that \eqref{01} holds, we prove the following:

\begin{thm}\label{teo2introB}
Let $N>0$ and $\frac{1}{N+1} < a < 1$. Assume that  $\left\{u_k\right\}$  satisfies \eqref{probprimo}, \eqref{0010star} and \eqref{01}, then we have:

\begin{enumerate}[i)]

\item If $0 < N < 1$ then $\beta^\infty = \frac{4}{a}$.

\item If $N \geq 1$  and $\frac{1}{N+1} < a \ne 1 \leq \frac{2}{N+1}$ then $\beta^\infty = \frac{4}{a}$ or it  satisfies one of the following :

\begin{equation}\label{03bis}
\begin{split}
\beta^\infty = &\frac{2}{a} \left( 1 + \sqrt{1-4am(1-a(N+1))} \right), \text{ with }  a \in \left(\frac{1}{N+1},\frac{1}{2} \right)  \\
& \text{ and } m \in \mathbb{N}: 2 \leq m \leq N+1
\end{split}
\end{equation}

\begin{equation}\label{04bis}
\begin{split}
&\beta^\infty = \sum_{j=1}^{m} \beta_j , \quad \sum_{j=1}^{m} \beta_{j}^{2} = \frac{\beta^\infty (4N - (1-a)\beta^\infty)}{a(N+1)-1} \text{ with } \beta_{j} \in \left[4, \frac{4}{a} \right], \\
&  \; \max_{j=1,...,m} \beta_j \geq \frac{2}{a} \text{ and } m \in \mathbb{N}: 2 \leq m \leq N + 1- \max \left\{ 0, \frac{1-2a}{2a} \right\} 
\end{split} 
\end{equation}
In particular, the value $\beta^\infty = 4(N+1)$ can be attained only if $N \in \mathbb{N}$ and $m = N+1$ in \eqref{03bis} (when $a \in \left( 0, \frac12 \right)$) or in \eqref{04bis} (when $a \geq \frac12$).

\item if $N>1$ and $\frac{2}{N+1} < a < 1$ then $\beta^\infty \geq 2(N+1)$ and $\beta^\infty$ is given by either one of  \eqref{03bis} and  \eqref{04bis} or it satisfies:

\begin{equation}\label{06}
2(N+1) \leq \beta^\infty \leq 4(N+1) - \frac{4}{a}.
\end{equation}
\end{enumerate}
\end{thm}
\qed

We observe that, if $N \in (0,1)$ and $0 < a \ne \frac{1}{N+1} < 1$ then the values of $\beta^\infty$ as given above, identify exactly the end points of the (sharp) interval range of $\beta$ for the \underline{radial} solvability of \eqref{prob}, as specified in Theorem \ref{jjj7}. Thus in this case, it is reasonable to expect that actually $u_k$ is itself radially symmetric (for large $k$). The issue of radial symmetry for solutions of \eqref{prob} has been addressed in \cite{tar?}, where in particular (by following \cite{blt}), it is shown that also for $N=1$ we have: 
$\beta^\infty = \frac{4}{a}$.

\begin{rem}
As already mentioned, the value $\beta^\infty$ in \eqref{03bis} arises in connection with the existence of a conformal metric on the Riemann sphere with constant curvature equal to one and assigned conical singularities. But in this case, one of the conical angle is larger than $2\pi$, and 
in view of the work in \cite{troy1}, \cite{troy2}, \cite{lt}, \cite{er1}, \cite{egt1},\cite{egt2},\cite{egt3}, \cite{uy}, \cite{cwx}, this feature may induce geometrical obstructions such to prevent the accurance of  \eqref{03bis}.
On the contrary in Section 6 we give an explicit example  of a solution sequence that admits the blow-up behaviour described in  \eqref{04bis}. Furthermore, in the given example, all the $\beta_j$'s coincide and are explicitly identified, and we suspect that this is the only possible instance which gives rise to \eqref{04bis}
\end{rem}

Finally, concerning the case $a>1$, we face a much more delicate  situation whose detailed analysis requires further investigation to be pursued in the future. For the moment we point out the following:

\begin{thm}\label{teo3introA}
Let $N>0$ and $a>1$. Assume that  $\left\{u_k\right\}$  satisfies \eqref{probprimo}, \eqref{0010star} and \eqref{01}, then we have:

\begin{enumerate}[i)]

\item If $0 < N < 1$ and $1 < a \leq \frac{2}{N+1}$ then $\beta^\infty = \frac{4}{a}$.

\item If  $a > \max \left\{ 1, \frac{2}{N+1}\right\}$ then \eqref{06} holds.

\end{enumerate}
\end{thm}
\qed

Even more intricate is the case when $a>1$ and \eqref{02} holds, where we have:

\begin{thm}\label{teo3introB}
Let $N>0$ and $a>1$. Assume that  $\left\{u_k\right\}$  satisfies \eqref{probprimo}, \eqref{0010star},  \eqref{002bis} and \eqref{02}, then we have:

\begin{enumerate}[i)]

\item If $0 < N < 1$ and $1 < a \leq \frac{2}{N+1}$ then $\beta^\infty = \frac{4}{a}$ or $\beta^\infty = 4(N+1)$.

\item If  $a > \max \left\{ 1, \frac{2}{N+1}\right\}$ then $\beta^\infty = 4 (N+1)$  or $\beta^\infty$  satisfies \eqref{06}  unless it takes one of the following values:

\begin{equation}\label{07}
\beta^\infty = 2 \left( N + 1 + \sqrt{(N+1)^2+ \frac{4m}{a^2}(1-a(N+1))} \right), \text{ with }  m \in \mathbb{N} \text{ and } a > 2;
\end{equation}

\begin{equation*}
\beta^\infty \text{ satisfies \eqref{04bis} with } \beta_{j} \in \left[ \frac{4}{a}, 4 \right],   \max_{j=1,...,m} \beta_j \geq 2 \text{ and } 2 \leq m \leq a(N+1) - \max \left\{0, \frac{a-2}{2} \right\}.
\end{equation*}

\end{enumerate}
\end{thm}
\qed

\begin{rem}\label{rem0.2}
Again we observe that, \eqref{07} occurs in connection with the existence of a conformal Riemann sphere with conical singularities and constant curvature equal to one. In view of the comments in Remark~\ref{rem0.1}  we suspect that actually \eqref{07} can occur only for $m=1$, in account of the value: $\beta^\infty = 4(N+1) - \frac{4}{a}$.
\end{rem}

In connection with Remark~\ref{rem0.2}, or more generally with the condition \eqref{06},  we recall that a {\bf{radial}} solution $u$ of \eqref{prob} satisfies the necessary conditions:

\begin{equation}\label{minnn}
\dfrac{1}{2\pi}\int_{\RR^{2}}\left|x\right|^{2N}e^{u}<4(N+1)\quad\mbox{ and }\quad\dfrac{1}{2\pi}\int_{\RR^{2}}e^{au}<\dfrac{4}{a},
\end{equation}
see \cite{pt2}.

Therefore, for {\bf{radial}} $u_k$ we can use \eqref{lim} together with \eqref{minnn} and obtain,

\begin{eqnarray}
\label{eq0.8} &\beta^{\infty}\geq\dfrac{4}{a}-4(N+1)\quad\mbox{ for }\quad 0<a<\frac{1}{2(N+1)}\\
\notag &\beta^{\infty}\geq4(N+1)-\dfrac{4}{a}\quad\mbox{ for }\quad a>\frac{2}{(N+1)}.
\end{eqnarray}

 Consequently, for {\bf{radial}} $u_k$, the condition \eqref{06} just turns into the following identity: $\beta^{\infty}=4(N+1)-\dfrac{4}{a}$, and similarly we check that indeed \eqref{07} can hold only with $m=1$.\\  Consistently, we point out that, the other alternatives for $\beta^\infty$ indicated above do occur in account of the (possible) non--radial features of $u_k$.  

\begin{rem} It is an interesting open problem to see whether \eqref{minnn}  remains valid also for \underline{non--radial} solutions. If so, one obtains (as for radial solutions) that the statements of the results above should be  improved accordingly.
\end{rem}

Next, we proceed to interpret the above results in terms of blow--up sets. To this purpose, for a solution $u_k$ of \eqref{probprimo}, we define (via Kelvin transformation) the function: 
$$
\hat{u}_{k}(x):=u_{k}\left(\dfrac{x}{\left|x\right|^{2}}\right)+\beta_{k}\log\dfrac{1}{\left|x\right|}.
$$
By virtue of \eqref{probd}, $\hat{u}_{k}$ is well defined in $\RR^{2}$ and satisfies:
\begin{equation}\label{kininttr}
\begin{cases}
-\Delta \hat{u}_{k}=\dfrac{1}{\left|x\right|^{2(N+2)-\beta_k}}e^{\hat{u}_{k}}+\dfrac{1}{\left|x\right|^{4-a\beta_{k}}}e^{a\hat{u}_{k}}=:\hat{f}_{k}\quad\mbox{in }\RR^{2}\\
\displaystyle{\dfrac{1}{2\pi}\int_{\RR^{2}}\left(\dfrac{1}{\left|x\right|^{2(N+2)-\beta_{k}}}e^{\hat{u}_{k}}+\dfrac{1}{\left|x\right|^{4-a\beta_{k}}}e^{a\hat{u}_{k}}\right)}
=\displaystyle{\dfrac{1}{2\pi}\int_{\RR^{2}}\left(\left|x\right|^{2N}e^{u_k}+e^{au_k}\right)}=\beta_k.
\end{cases}
\end{equation}
By using the standard definition of  blow-up point, as given by Brezis--Merle in \cite{bm} (see Definition~\ref{dblup} below), we may consider the (possibly empty) blow--up set of $\left\{u_k\right\}$, and  denote it by $S$. Analogously, we denote by $\hat{S}$ the (possibly empty) blow--up set of $\left\{\hat{u}_{k}\right\}$.
\begin{rem}\label{rem0.2bis}
The condition \eqref{002bis} simply states that,
\begin{equation}\label{00star00}
S\cup\hat{S}\neq\emptyset,
\end{equation}
and in particular it implies that, if $S=\emptyset$ then  $\hat{S}=\left\{0\right\}$.
Notice also that if $z_0\in S$ and $z_0\neq 0$ then $\dfrac{z_0}{\left|z_0\right|^{2}}\in\hat{S}$ 
\end{rem} 
We point out that in order to develop a blow--up analysis for problem \eqref{kininttr} around the \underline{origin} (in the spirit of \cite{bm}, \cite{bt}, \cite{ls} and \cite{barm}) one needs to require that:
$$\beta^{\infty} > \max \left\{ \frac{2}{a}, 2(N+1) \right\}.$$
In account of \eqref{betakkkkkkkk}, we shall indicate in Lemma~\ref{lem42} how to deal with the situation where  $\beta = \max \left\{ \frac{2}{a}, 2(N+1) \right\}$.\\

Furthermore, if $0<a<1$ then we shall show that in most cases we have: $S=\left\{0\right\}$ and/or $\hat{S}=\left\{0\right\}$. But for $a>1$, this should not be expected any longer, since now the first term in the right hand side of \eqref{probprimo} overpowers the second one, and so the role of the origin becomes irrelevant.\\

More precisely, we establish the following:
\begin{thm}\label{zerozeroc}
Let $N>0$ and assume \eqref{001star1}. Suppose that $u_k$ satisfies \eqref{probprimo}, \eqref{0010star} and its blow--up set $S\neq\emptyset$. We have: 
 \begin{enumerate}[i)]
 \item for $0<a<1$ then for $S$ one of the following alternatives holds :
\begin{itemize}
\item[$\bullet$] $S=\left\{0\right\}$. 
\item[$\bullet$]  $N\in\NN$ and $S$ is formed by the vertices of a $(N+1)$-regular polygon, namely in complex notation: $S=\left\{z_1,...,z_{N+1}\right\}\subset\RR^{2}\setminus\left\{0\right\}$,  with \begin{equation}\label{scalintronuova}
z^{N+1}_{k}=\xi_{0}, \quad k=1,...,N+1,
\end{equation}
and $\xi_{0}=(-1)^Nz_1\cdot...\cdot z_{N+1}$.
\item[$\bullet$]  In case $0\in S$ and $ S\setminus \left\{ 0 \right\}\neq\emptyset$ then $N>2$, $ \dfrac{2}{N}<a<1$ and $\beta^\infty=2(N+1+m)$, where $m$ is the number of points in $S\setminus \left\{ 0 \right\}$ and it satisfies:  $1\leq m < N+1 - \frac{2}{a}$. In particular, $2(N+1)<\beta^\infty \leq 4(N+1) - \frac{4}{a}$, in this case.
\end{itemize}
\item For $a>1$ then $S=\left\{z_0\right\}$, for some $z_0\in\RR^{2}$.

\end{enumerate}

\end{thm} 
\qed 

To understand the nature of alternative \eqref{scalintronuova}, one should compare it with the blow-up behaviour exhibited by solutions of \eqref{lspeintro} and described in Remark~\ref{rembbbbbis} below. Thus, alternative \eqref{scalintronuova} indicates the possibility that  an analogous  blow-up behaviour could be attained also by solutions of \eqref{prob}. \\

In addition we have: 

\begin{prop}\label{prop1.3}
Let $N>0$ and assume \eqref{001star1}. If $u_k$ satisfies \eqref{probprimo}, \eqref{0010star} and \eqref{002bis} we have:
\begin{enumerate}[i)]
\item for $0<a<\dfrac{1}{2(N+1)}$ or $a>\max\left\{1,\dfrac{2}{N+1}\right\}$, then necessarily $0\in\hat{S}$.
\item for $\dfrac{1}{2(N+1)}<a\neq\dfrac{1}{N+1}<\dfrac{2}{N+1}\ \text {and S}  \neq\emptyset$ then $0\notin\hat{S}$.

\end{enumerate}

\end{prop} 
\qed  

We conclude with the observation that the blow--up analysis for solutions of \eqref{probprimo}, becomes most delicate when we deal with the situation described in part $i)$ of Proposition~\ref{prop1.3}.\\

This paper is organised as follows, in Section 3 we introduce some general tools useful for the blow-up analysis developed in Section 4 and 5. In particular, in Section 5 we investigate the nature of the blow-up sets $S$ and $\hat{S}$ . 
The remaining Sections 6, 7, and 8 are devoted to establish the results stated above  respectively in the cases: $0<a<\frac{1}{N+1}$, $\frac{1}{N+1}<a<1$ and $a>1$.

\underline{Acknowledgement}: We wish to express our gratitude to Roberto Tauraso and Carlo Pagano for their insight in the proof of Theorem \ref{new}.

\section{Useful facts}
\setcounter{equation}{0}

As well known (see \cite{cl}, \cite{cl1}, \cite{pt} and also \cite{tar}), for $b>0$ and $N>-1$ every solution of the following problem:

\begin{equation}\label{plsc2}
\begin{cases}
-\Delta w=\left|x\right|^{2N}e^{bw}\quad\mbox{in }\RR^{2},\\
\displaystyle{\int_{\RR^{2}}\left|x\right|^{2N}e^{bw}}<+\infty,
\end{cases}
\end{equation}

takes (in complex notations) the form:

\begin{equation}\label{profg}
w(z)={\rm{log}}\left[\frac{8(N+1)^{2}\mu}{b\left(1+\mu\left|z^{N+1}-c\right|^{2}\right)^{2}}\right]^{1/b}
\end{equation} 

where $\mu>0$, $c\in\CC$ and $c=0$ when $N\notin\NN\cup\left\{0\right\}.$

In particular we have,

\begin{equation}\label{inpa}
\frac{1}{2\pi}\displaystyle{\int_{\RR^{2}}\left|x\right|^{2N}e^{bw}}=\frac{4(N+1)}{b}.
\end{equation}

\begin{rem}\label{rembbbbbis}
When $N\in\NN$, we see that the solution $w=w_{\mu}$ given in \eqref{profg} with $c\neq 0$, is not radially symmetric and actually it 'concentrates' exactly at the $(N+1)-$roots of $c\in\CC$ as $\mu\to+\infty$. We should keep this fact in mind in order to justify the statement $(i)$ in Theorem~\ref{new} established in the sequel.
\end{rem}

For the more general Liouville type equation:

\begin{equation}\label{preq1}
\begin{cases}
-\Delta u=e^{au}+\left|x\right|^{2N}e^u\quad\mbox{  in }\RR^2\\
\displaystyle{\int_{\RR^2} e^{au}+\left|x\right|^{2N}e^u \,dx<\infty}
\end{cases}
\end{equation}
a classification results of the type \eqref{profg} is not  available, however analogous qualitative informations about the solution are available as follows.

\begin{prop}\label{verde1}
Let u be a solution of \eqref{preq1} with 
\begin{equation}\label{1star}
\beta=\frac{1}{2\pi}\displaystyle{\int_{\RR^2}\left(e^{au}+\left|x\right|^{2N}e^u\right)dx }.
\end{equation}
Then the following estimates hold:
\begin{itemize}
\item[(i)]\begin{equation}\label{estib}
-\beta\,\mbox{\emph {log}}(\left|x\right|+1)-C\leq u(x)\leq -\beta\,\mbox{\emph {log}}(\left|x\right|+1)+C\quad \mbox {in }\RR^2
\end{equation}
with a suitable constant $C=C(a,\beta,u(0),N)>0$, and
\begin{equation}\label{dd2estib}
ru_r\to-\beta\quad\mbox{ and }\quad u_{\vartheta}\to 0\quad\mbox{ as }r\to+\infty,
\end{equation}
with $(r, \vartheta)$ the polar coordinates in $\RR^{2}$.

In particular, 
\begin{equation}\label{lowbeta}
\beta>\mbox{max }\left\{2(N+1),\frac{2}{a}\right\}.
\end{equation}
\item[(ii)]  \textbf{Green's representation formula}: for any $x,\overline{x}\in \RR^{2}$
\begin{itemize}
\item[(a)]
\begin{equation}\label{greenrep}
u(x)-u(\overline{x})=\frac{1}{2\pi}\int_{\RR^{2}}\log\left(\frac{\left|\overline{x}-y\right|}{\left|x-y\right|}\right)\left(\left|y\right|^{2N}e^{u(y)}+e^{au(y)}\right)dy,
\end{equation}
\item[(b)] 
\begin{equation}
\nabla u(x)=-\frac{1}{2\pi}\int_{\RR^{2}}\frac{x-y}{\left|x-y\right|^{2}}\left(\left|y\right|^{2N}e^{u(y)}+e^{au(y)}\right)dy.
\end{equation}
\end{itemize}
\end{itemize}
\end{prop}
\qed

The above results follow as in \cite{cl}, and details may be found in \cite{cgs} or \cite{f}.

We also mention the following useful local identity of Pohozaev's type:
\begin{thm}[Pohozaev's identity]\label{thmpooho} Let $N>-1$ and $u$ be a solution of \eqref{preq1}. Then for any $r>0$ the
following identity holds:
\begin{equation}\label{poho2}
\begin{aligned}
r\int_{\partial B_r}&\left(\frac{\left|\nabla u\right|^2}{2}-(\nu\cdot\nabla u)^2\right)\mbox{d}\sigma\\
&=\frac{r}{a}\int_{\partial B_r}
e^{au}\mbox{d}\sigma +r^{2N+1}\int_{\partial B_r} e^u\mbox{d}\sigma-\frac{2}{a}\int_{B_r}e^{au}-2(N+1)\int_{B_r}\left|x\right|^{2N}e^{u}
\end{aligned}
\end{equation}
with $\nu$ is the outward normal vector to $\partial B_r$.
\end{thm}
\qed  

Identity \eqref{poho2} follows in the usual way, and it has been derived in \cite{cgs} in case $N\geq 0$, while in \cite{f} it has been shown how to extend it for the more general case $N>-1$.

Actually the asymptotic estimates \eqref{estib} and \eqref{dd2estib} allow one to pass to the limit in \eqref{poho2} as $r\to \infty$ and obtain:

\begin{coro}
Let $u$ be a solution of \eqref{preq1} with $\beta=\frac{1}{2\pi}\displaystyle{\int_{\RR^2}\left(e^{au}+\left|x\right|^{2N}e^u\right)}$. Then
$$
2N\int_{\RR^2}\left|x\right|^{2N}e^u dx+2\left(\frac{1}{a}-1\right)\int_{\RR^2}e^{au}{dx}=\pi\beta(\beta-4).
$$
In particular,

\begin{equation}
\begin{split}
\label{lim2} \mbox{if } a \ne \frac{1}{N+1}, \mbox{ then } \quad &\frac{1}{2\pi}\int_{\RR^2}\left|x\right|^{2N} e^{u} =\frac{\beta(4-a\beta)}{4[1-a(N+1)]},\\
\mbox{ and } \quad &\frac{1}{2\pi}\int_{\RR^2} e^{au} =\frac{a\beta(\beta-4(N+1))}{4[1-a(N+1)]}
\end{split}
\end{equation}

\begin{equation}\label{lim22}
\mbox{if } a = \frac{1}{N+1}, \mbox{ then } \quad \beta=4(N+1)(=\frac{4}{a}).
\end{equation}

\end{coro}
\qed

\begin{rem}\label{remn1}
From \eqref{lim2} it follows that, for $a\neq \frac{1}{N+1}$ a necessary condition for the solvability of \eqref{preq1}-\eqref{1star} is given by: 

\begin{equation}\label{intt}
\min\left\{4(N+1), \frac{4}{a}\right\}<\beta<\max\left\{4(N+1), \frac{4}{a}\right\},
\end{equation}

which must hold together with \eqref{lowbeta}. Since \eqref{intt} implies \eqref{lowbeta} exactly when

\begin{equation}\label{2.14bis}
\frac{1}{2(N+1)} \leq a\neq\frac{1}{N+1}\leq\frac{2}{N+1}.
\end{equation}

from Theorem 1.1 (see \cite{pot}, \cite{pt2}), we know that when \eqref{2.14bis} holds, then \eqref{intt} gives a \underline{necessary} and  \underline{sufficient} condition for the solvability of \eqref{preq1}-\eqref{1star}.
\end{rem}

\section{Local blow--up analysis}\label{bblu}
\setcounter{equation}{0}

As discussed above, we shall assume that,

\begin{equation}\label{eq3.0}
N > -1,\quad 0 <  a \ne \frac{1}{N+1} \; \text{and} \; a \ne 1.
\end{equation}
unless is otherwise specified.

Let $u_k$ satisfy: 

\begin{equation}\label{kappa}
\begin{cases}
-\Delta u_k=e^{au_k}+\left|x\right|^{2N}e^{u_k}=:f_k\\
\beta_k:=\frac{1}{2\pi}\displaystyle{\int_{\RR^{2}}\left(e^{au_k}+\left|x\right|^{2N}e^{u_k}\right)dx}.
\end{cases}
\end{equation} 

so that from \eqref{lowbeta} and \eqref{intt} we have,

$$
\beta_k>\max\left\{\frac{2}{a}, 2(N+1)\right\}\quad\mbox{and}\quad\min\left\{\frac{4}{a}, 4(N+1)\right\}<\beta_k<\max\left\{\frac{4}{a}, 4(N+1)\right\}.
$$

Therefore, by taking a subsequence if necessary, we can always suppose that,

\begin{equation}\label{131bis}
\beta_k\to\beta^\infty := \beta \quad\mbox{ as }\quad k\to+\infty,
\end{equation} 
with
\begin{equation}\label{toclim2}
\beta\geq\max\left\{\frac{2}{a},2(N+1)\right\} \quad \mbox{ and } \quad \min\left\{\frac{4}{a}, 4(N+1)\right\}\leq\beta\leq\max\left\{\frac{4}{a}, 4(N+1)\right\}.
\end{equation}

Furthermore, from \eqref{lim2} we also know that,

\begin{equation}\label{33bis}
\lim_{k\to\infty} \frac{1}{2\pi}\int_{\RR^2}\left|x\right|^{2N}e^{u_k}=\frac{\beta(4-a\beta)}{4[1-a(N+1)]} \,
 \mbox{ and }\, \lim_{k\to\infty}\frac{1}{2\pi}\int_{\RR^2}e^{au_k}=\frac{a\beta(\beta-4(N+1))}{4[1-a(N+1)]}.
\end{equation}

As in \cite{bm}, we give the following notion of blow--up point:

\begin{defn}\label{dblup}
A point $x_0$ is called a blow--up point for $u_k$ if there exists a sequence $\{x_k \} \subset \mathbb{R}^2$: $x_k \to x_0$ and $u_k (x_k) \to + \infty$.
\end{defn}

\begin{prop}\label{prop17}
If $u_k$ satisfies \eqref{kappa}, \eqref{131bis}, then its blow--up set $S$ may contain only a finite number of points. Furthermore, for every $x_0 \in S$,
\begin{equation}\label{propar}
\beta (x_0) := \; \displaystyle{\lim_{r\to 0} \liminf_{k \to + \infty} \left( \frac{1}{2\pi} \int_{B_r(x_0)} f_k(x) dx\right)} \geq \min\left\{ 2(N^- +1), \frac{2}{a}\right\},
\end{equation}
with $N^-=\min \left\{ 0, N \right\}$. In particular, for $N \geq 0$ we have: $\beta (x_0) \geq \frac{2}{\max\{a,1\}}$
\end{prop}

In order to establish \eqref{propar} we recall the following well known fact established first in \cite{bm} (see e.g. Lemma~5.2.1 of \cite{tar2}), which here we state in a form suitable for our purposes:

\begin{lemma}\label{lem21}
Let $\Omega \subset \mathbb{R}^2$ be a bounded open set and let $u_k$ satisfy:
\begin{equation}\label{prob345}
\begin{cases}
-\Delta u_k=f_k \in L^1(\Omega)\\
\displaystyle{\limsup_{k\to \infty} \left( \Vert u_k^+\Vert_{L^1(\Omega)}+ \Vert f_k\Vert_{L^1(\Omega)} \right) < +\infty }
\end{cases}
\end{equation}
For every, 
$$0<p<\left(\frac{1}{4\pi} \displaystyle{\limsup_{k\to \infty} \Vert f_k\Vert_{L^1(\Omega)}} \right)^{-1},$$
there exists a constant $C_p>0$ such that:
$$\displaystyle{\int_{\Omega}e^{pu_k}\leq C_p}$$
\end{lemma}
\qed

\begin{rem}\label{rem21}
In order to use Lemma \ref{lem21} for $u_k$ satisfying \eqref{kappa} (or for $\hat{u}_k$ defined in \eqref{kin0} below), we point out that, if the following condition holds:

$$\displaystyle{\int_{\Omega}W_ke^{au_k} + \int_{\Omega}\frac{1}{W^q_k}\leq C_{0},}$$

with $q>0$, $a>0$, $W_k \geq 0$ and $C_{0}>0$; then we may derive that, $\Vert u_k^+ \Vert_{L^1(\Omega)} \leq C$, with a suitable constant $C>0$,  (see the estimate $(5.3.6.)$ in \cite{tar2}).

Actually the above condition will be useful also for ''scaled'' versions of $u_k$.
\end{rem}

\dimo (of Proposition~\ref{prop17})
To establish \eqref{propar}, we observe that, if $\Omega \subset \mathbb{R}^2$ and $f_k:=e^{au_k}+\vert x \vert^{2N} e^{u_k}$ satisfy:
$$ \displaystyle{ \limsup_{k \to + \infty} \frac{1}{2\pi} \int_{\Omega} f_k(x) dx < \min\left\{ 2(N^- +1), \frac{2}{a}\right\}},
$$
then, on the basis of Remark~\ref{rem21}, we can apply Lemma~\ref{lem21} to $u_k$ in $\Omega$. Thus we obtain that, $f_k$ is uniformly bounded in $L^p(\Omega)$, for suitable $p>1$, and by 
standard elliptic regularity theory, we conclude that $u_k^+$ is uniformly bounded in $L^{\infty}_{loc}(\Omega)$. Therefore $\Omega$ cannot contain any blow--up point of $u_k$ in this case. 
Consequently, for $x_0 \in S$ then necessarily \eqref{propar} must hold, and in view of \eqref{131bis}, the set $S$ must be finite.
\qed

We show next how to improve \eqref{propar} on the basis of a suitable Alexandrov Bol's inequality (see \cite{ban}, Theorem 6.4 in \cite{s}), which in particular implies the following: 

\begin{thm}\label{ab}
Let $\Omega\subset\RR^{2}$ be a bounded domain with smooth boundary $\partial\Omega$. If $p\in C^{2}(\Omega)\cap C^{0}(\overline{\Omega})$ satisfies
\begin{equation}\label{logp}
-\Delta {\rm{log}}p\leq p\quad \mbox{ in }\quad\Omega,
\end{equation}
and $\Sigma:=\displaystyle{\int_{\Omega}p(x)dx}< 8\pi$, then we have:
\begin{equation}\label{alb}
\max_{\overline{\Omega}}p\leq\left(1-\frac{\Sigma}{8\pi}\right)^{-2}\max_{\partial\Omega}p.
\end{equation}
\end{thm}

See \cite{s} for a proof.

\begin{prop}\label{teoalmenoimp}
Suppose that $\left\{u_k\right\}$ sastifies \eqref{kappa}, with $N \geq 1$. If $x_0\in S$ then for $\beta (x_0)$ in \eqref{propar} there holds:

\begin{equation}\label{almenoimp}
\beta ( x_0 ) \geq \frac{4}{\max\left\{1,a\right\}}
\end{equation}

\end{prop}

The estimate \eqref{almenoimp} will be a consequence of the following lemma which is of independent interest.

\begin{lemma}\label{135bbis}

Under the assumption of Proposition \ref{teoalmenoimp}, let $\eta_k=\max\left\{1,a\right\}\left(e^{au_k}+\left|x\right|^{2N}e^{u_k}\right)$ then

\begin{equation}\label{1310bis}
-\Delta\log\eta_k\leq\eta_k.
\end{equation}

\end{lemma}

\dimo
Set,
\begin{equation}\label{prov2}
V_{k}(x)=\left|x\right|^{2N}e^{(1-a)u_k}
\end{equation}
and 
\begin{equation}\label{prov02}
\xi_{k}(x):=u_{k}(x)+\frac{1}{a}\log\left(1+V_{k}(x)\right),
\end{equation}
so that
\begin{equation}\label{prov01}
-\Delta u_k=\left(1+V_{k}(x)\right)e^{au_k}=e^{a\xi_k}.
\end{equation}
In order to compute $\Delta\xi_k$,
we observe that, 
\begin{equation}\label{prov3}
\Delta\log\left(1+V_{k}(x)\right)=\frac{\Delta V_{k}(x)}{1+V_{k}(x)}-\frac{\left|\nabla V_{k}(x)\right|^{2}}{\left(1+V_{k}(x)\right)^{2}},
\end{equation}
and 
\begin{equation}\label{prov5}
\Delta V_{k}=(1-a)V_{k}(x)\Delta u_k+\frac{\left|\nabla V_{k}(x)\right|^{2}}{V_{k}(x)}
\end{equation}
Using \eqref{prov5} in \eqref{prov3} we get:
\begin{equation}\label{prov6}
\Delta\log\left(1+V_{k}(x)\right)=\left(\frac{(1-a)V_{k}(x)}{1+V_{k}(x)}\right)\Delta u_k+\frac{\left|\nabla V_{k}(x)\right|^{2}}{V_{k}(x)\left(1+V_{k}(x)\right)^{2}}.
\end{equation}
Therefore from \eqref{prov02} and \eqref{prov6} we obtain,
\begin{equation}\label{prov7}
\begin{aligned}
-\Delta\xi_k&=-\Delta u_k\left(\frac{a+V_{k}(x)}{a\left(1+V_{k}(x)\right)}\right)-\frac{\left|\nabla V_{k}(x)\right|^{2}}{aV_{k}(x)\left(1+V_{k}(x)\right)^{2}}\\
&\leq e^{a\xi_k}\left(\frac{a+V_{k}(x)}{a\left(1+V_{k}(x)\right)}\right),\quad\mbox{in }\RR^{2}.
\end{aligned}
\end{equation}
Observe that, in case $a\in(0,1]$, then from \eqref{prov7} we find,
\begin{equation}\label{prov8}
-\Delta\xi_k\leq \frac{1}{a}\left(\frac{a+V_{k}(x)}{1+V_{k}(x)}\right)e^{a\xi_k}\leq\frac{1}{a}e^{a\xi_k}.
\end{equation}
Consequently, in this case we have: $\eta_{k}:=e^{a\xi_k}$ and it satisfies:
\begin{equation}\label{prov9}
-\Delta\log\eta_k\leq \eta_k,
\end{equation}
and \eqref{1310bis} is established.

On the other hand, for $a>1$, the previous computations yields to,
$$
-\Delta\xi_k\leq \left(\frac{1+\frac{1}{a}V_{k}(x)}{1+V_{k}(x)}\right)e^{a\xi_k}\leq e^{a\xi_k},
$$
and the desired conclusion follows for $\eta_k=e^{a\xi_k+\log a}$.
\qed 

{\bf{Proof of Proposition~\ref{teoalmenoimp}}}
For $0<a \leq 1$ we can apply Theorem~\ref{ab} to $\eta_k$ in $B_{\delta}(x_0)$ with $x_0\in S$ and $\delta>0$ sufficiently small, so that $S \cap \overline{B}_{\delta}(x_0)=\left\{x_0\right\}$. Hence, $\max_{\partial B_{\delta}(x_0)} u_k\leq C$ and as a consequence, $\max_{\partial B_{\delta}(x_0)} \eta_k \leq C$. Thus if by contradiction we suppose that,
$$\limsup_{k\to+\infty}\frac{1}{2\pi}\displaystyle{\int_{B_{\delta}(x_0)}\eta_{k}<4},$$ 

then from \eqref{alb} we would find that $\max_{\overline{B}_{\delta/2}(x_0)}\eta_k\leq C$, in contradiction with the fact that $x_0\in S$.

Hence $\beta(x_0) \geq 4$, and \eqref{almenoimp} is established for $a\in(0,1]$.

For $a>1$, the argument above applies to $\eta_k=e^{a\xi_{k}+\log a}$ and it yields to the following:

$$
a \beta (x_0) := \lim_{\delta \to 0^+} \liminf_{k\to+\infty}\frac{a}{2\pi}\int_{B_{\delta}(x_0)}e^{a\xi_k} \geq 4.
$$

and so \eqref{almenoimp} is established also for $a > 1$.$\hfill\Box$

\begin{rem}
Actually, by using a sharper version of Theorem \ref{ab} where \eqref{logp} is assumed only in the sense of distributions, it is possible to obtain an analogous improvement of \eqref{propar} for any $N > -1$. We refer to \cite{BaCa} for details. 
\end{rem}

Next, we establish uniform estimates for $u_k$, away from the blow--up set.

\begin{prop}\label{1310}
Let $\left\{u_k\right\}$ be a sequence of solution for \eqref{kappa}, \eqref{131bis}, and let $S$ (possibly empty) be its blow--up set. For any compact set $K\subset\RR^{2}\setminus S$ there exists $C=C(K)>0$ such that,
\begin{equation}\label{cnblon}
\max_{K}u_k-\min_{K}u_k\leq C.
\end{equation}
\end{prop}

\dimo
We shall discuss the more delicate case where $S\neq \emptyset$, since for $S=\emptyset$ then \eqref{cnblon} follows by similar yet simpler arguments. For a compact set $K\subset\RR^{2}\setminus S$, 
let $\delta=\mbox{dist}\left(K, S\right)>0$, and set,
 
$$\mathcal{N}_{\delta}=\mathcal{N}_{\delta}(S)=\left\{x\in\RR^{2}: \,\,\mbox{dist}\left(x, S\right)<\frac{\delta}{2} \right\}.$$

Furthermore, we take $R>\delta$ large enough, so that: $K\cup\mathcal{N}_{\delta}\subset B_{R}$ and for $x, \overline{x} \in K$ we use \eqref{greenrep} to estimate:
\begin{equation}\label{grepo}
\begin{aligned}
&\left|u_{k}(x)-u_{k}(\overline{x})\right|\leq\frac{1}{2\pi}\int_{\RR^2}\left|\log\frac{\left|\overline{x}-y\right|}{\left|x-y\right|}\right|f_k(y)dy\leq\\
&\frac{1}{2\pi}\int_{\mathcal{N}_{\delta}}...+\frac{1}{2\pi}\int_{B_{2R}\setminus\mathcal{N}_{\delta}}...+\frac{1}{2\pi}\int_{\left\{\left|x\right|\geq 2R\right\}}=:I_{1,k}+I_{2,k}+I_{3,k}.
\end{aligned}
\end{equation}
Since for $y\in\mathcal{N}_{\delta}$ we have:
$$
\left|\log\frac{\left|\overline{x}-y\right|}{\left|x-y\right|}\right|\leq\log\frac{4R}{\delta},
$$ 
we find,
$$
I_{1,k}:=\frac{1}{2\pi}\int_{\mathcal{N}_{\delta}}\left|\log\frac{\left|\overline{x}-y\right|}{\left|x-y\right|}\right|f_k(y)dy\leq \log\left(\frac{4R}{\delta}\right)\int_{\RR^{2}}f_{k}(y)dy\leq C.
$$
Here and in the following $C>0$ denotes a constant whose value may change from line to line, but is always independent of the specific points $x$ and $\overline{x}$.

For the second integral we note that $\sup_{B_{2R}\setminus\mathcal{N}_{\delta}}f_k\leq C$ and so:
\begin{equation}
\begin{aligned}
I_{2,k}&:=\frac{1}{2\pi}\int_{B_{2R}}\left|\log\frac{\left|\overline{x}-y\right|}{\left|x-y\right|}\right|f_k(y)dy\leq
C\int_{B_{2R}\setminus\mathcal{N}_{\delta}}\left|\log\frac{\left|\overline{x}-y\right|}{\left| x-y\right|}\right|dy\\
&\leq C \left( \int_{B_1(x)}...+\int_{B_1(\overline{x})}...+\int_{B_{2R}\setminus \left( B_1(x) \cup B_1(\overline{x}) \right)} \right) \\
&\leq C\left(\int_{\left|x-y\right|<1}\log\left(\frac{1}{\left|x-y\right|}\right) dy +\int_{\left|\overline{x}-y\right|<1}\log\left(\frac{1}{\left|\overline{x}-y\right|} \right) dy +  R^2(\log R +1)  \right)\leq C.
\end{aligned}
\end{equation}

Finally, since for any $\left|y\right|\geq 2R$ and $x, \overline{x} \in K \subset B_R$ we have:

$$
\left| \log \frac{\left|x-y\right|}{\left|\overline{x}-y\right|} \right|\leq \log 4,
$$

and we see that,

\begin{equation}
I_{3,k}:=\frac{1}{2\pi}
\int_{\left|y\right|\geq 2R}\left|\log\left(\frac{\left|\overline{x}-y\right|}{\left|x-y\right|}\right)\right|f_k(y)dy\leq\frac{\log 4}{2\pi}\int_{\left|x\right|\geq 2R}f_k\leq C.
\end{equation}
and \eqref{cnblon} follows.\qed\\

To account for the blow--up behavior of $u_k$ near a given blow--up point, we start to recall the blow--up analysis which is available in \cite{bm,bt,ls,barm} for a sequence of solutions $\left\{u_k\right\}$ satisfying: 
\begin{equation}\label{lirich}
\begin{cases}
-\Delta u_k=\left|x\right|^{2\alpha_k}U_{k}e^{bu_k}+\sigma_k\quad\mbox{in }\Omega,\\
\displaystyle{\int_{\Omega}\left|x\right|^{2\alpha_k}U_{k}e^{bu_k}}\leq C,\\
\displaystyle{\max_{\partial \Omega}} u_k - \displaystyle{\min_{\partial \Omega}} u_k \leq C,
\end{cases}
\end{equation}
where $\Omega\subset\RR^{2}$ is a regular open and bounded domain, and
\begin{equation}\label{limpof}
\alpha_k \rightarrow \alpha > -1.
\end{equation}

In this context, we have:
\begin{prop}\label{knbl}
Let $\left\{u_k\right\}$ be a sequence of solutions for \eqref{lirich}, so that \eqref{limpof} holds. Suppose that,
$$
\begin{aligned}
& \sigma_k \in L^p(\Omega) \mbox{ and } \Vert \sigma_k \Vert_{L^p(\Omega)} \leq C \mbox{ for some } p>1; \\
& U_k \in C^{0,1}(\Omega) \mbox{ and } 0<a_1\leq U_k(x) \leq a_2, \vert \nabla U_k(x) \vert \leq A, \forall x \in \Omega.
\end{aligned}
$$

Then along a subsequence (denoted the same way) only one of the following alternatives hold:
\begin{itemize}
\item[(i)] $\left\{u_k\right\}$ is uniformly bounded in $L^{\infty}_{{{\rm}loc}}(\Omega)$;
\item[(ii)] $\displaystyle{\sup_{\Omega'}u_k\to-\infty}$, \mbox{ for every } $\Omega'\subset\subset\Omega$;
\item[(iii)] there exists a finite set $S=\left\{z_1,...,z_m\right\}\subset\Omega$, of blow--up points such that,
\begin{itemize}
\item[(a)] $\displaystyle{\sup_{\Omega'}u_k\to-\infty}$, \mbox{ for every } $\Omega'\subset\subset\Omega\setminus S$,
\item[(b)] $\frac{1}{2\pi}\left|x\right|^{2\alpha_k} U_k e^{b u_k}\rightharpoonup\displaystyle{\sum^{m}_{j=1}\beta(z_j)\delta_{z_j}}$ weakly in the sense of measures in $\Omega$,\\ with $\beta(z_j)=\frac{4}{b}$ if $z_j \neq 0$ and $\beta(z_j)=\frac{4}{b}(1+\alpha)$ if $z_j = 0$, for some $j \in \{ 1,...,m \}.$.

\end{itemize} 
\end{itemize} 
\end{prop}

\dimo
See \cite{bm}, \cite{barm}, \cite{tar}, \cite{ls}.
\qed

\begin{rem}\label{rem20bis}
By a direct inspection of the proof in \cite{bt} (see also \eqref{1337bis} below), it is possible to weaken the given assumption about $\nabla U_k$ around the origin as follows: $U_k \in C^{0,1}\left(\Omega \setminus \left\{ 0 \right\} \right)$, $\vert \nabla U_k(x) \vert$ is uniformly bounded in $L^{\infty}_{\rm{loc}}\left(\Omega \setminus \left\{ 0 \right\} \right)$ and $\displaystyle{\lim_{r\to 0} \lim_{k \to + \infty} \sup_{B_r} \vert x\cdot\nabla U_{k}\vert =0}$. 

In particular the result above applies to a solutions sequence $u_k$ satisfying \eqref{eq3.0} and \eqref{kappa} with $a = 1$ and $N > -1$. For this reason, in the following we may focus only to the case where $0 < a \ne 1$.
\end{rem}

\begin{coro}\label{coroc1}
Let $\left\{u_k\right\}$ satisfy \eqref{kappa}, \eqref{131bis} and assume \eqref{eq3.0}. Then Proposition~\ref{knbl} applies to $u_k$ considered on every open bounded set  $\Omega\subset\RR^{2}\setminus\left\{0\right\}$, or when $\Omega=B_r$, $r>0$ and, $a \in \left(0,\min\left\{1,\frac{1}{N+1}\right\}\right) \cup \left(\max\left\{1,\frac{1}{N+1}\right\},+\infty\right)$.

Furthermore,  property $(b)$ holds with $\beta(z_j)=4(N+1)$, when $0<a<\min\left\{1, \frac{1}{N+1}\right\}$ and $z_j = 0$, and with $\beta(z_j)=\frac{4}{\max\{1,a\}}$, for $j\in \left\{1,...,m \right\}$, in the other cases.
\end{coro}

\dimo
We only check the case $N\geq 0$, as the case $-1<N<0$ follows in analogous way.

If $0<a<1$ and $\Omega\subset\RR^{2}\setminus\left\{0\right\}$ then we easily check that we can apply Proposition~\ref{knbl} with $b=1$, $U_k\equiv 1$, $\alpha_k=N$ and $\sigma_k=e^{au_k}$.
We show next that the same holds for $\Omega=B_r$ and $0<a<\frac{1}{N+1}$. Indeed, we only need to check that $e^{au_k}$ is uniformly bounded in $L^{p}(\Omega)$ for a suitable exponent $p>1$. To this purpose we take $1<p<\frac{1}{a(N+1)}$, and use Holder inequality to estimate:

\begin{equation}\label{corohom}
\begin{aligned}
\int_{\Omega}e^{apu_k}&=\int_{\Omega}\left(\left|x\right|^{2N}e^{u_k}\right)^{ap}\frac{1}{\left|x\right|^{2Nap}} \, dx\\
&\leq\left(\int_{\Omega}\left|x\right|^{2N}e^{u_k}\right)^{ap}\left(\int_{\Omega}\frac{1}{\left|x\right|^\frac{2Nap}{1-ap}}\right)^{1-ap}\leq C \left(\int_{\Omega}\frac{1}{\left|x\right|^\frac{2Nap}{1-ap}}\right)^{1-ap}\leq C.
\end{aligned}
\end{equation}

Similarly, if $a>1$ then we check that Proposition~\ref{knbl} applies for any bounded set $\Omega \subset \mathbb{R}^2$, with $b=a$, $U_k\equiv 1$, $\alpha_k=0$ and $\sigma_k=\left|x\right|^{2N}e^{u_k}$.\\
\qed

By virtue of Corollary~\ref{coroc1}, (see also Remark \ref{rem20bis}) it remains to investigate the blow--up behavior of $u_k$ only when blow--up occurs at the \underline{origin} , and $\min\left\{1,\frac{1}{N+1} \right\} < a < \max\left\{1,\frac{1}{N+1} \right\},$  $N \neq 0$. 

As we shall see, a blow--up at the origin is the must likely situation to occur when $S\neq \emptyset$. and $0<a<1$\\

In order to account also for the behavior of $u_k$ at $\infty$, we use the Kelvin transform and define:
\begin{equation}\label{kin0}
\hat{u}_{k}(x)=u_{k}\left(\frac{x}{\left|x\right|^{2}}\right)+\beta_k{\rm {log}}\frac{1}{\left|x\right|}.
\end{equation}
Observe that $\hat{u}_{k}$ extends smoothly at the origin, (see estimate \eqref{estib}), and satisfies:
\begin{equation}\label{kin}
\begin{cases}
-\Delta \hat{u}_{k}=\frac{1}{\left|x\right|^{2(N+2)-\beta_k}}e^{\hat{u}_{k}}+\frac{1}{\left|x\right|^{4-a\beta_k}}e^{a\hat{u}_{k}}=\hat{f}_{k}\quad\mbox{in }\RR^{2}\\
\frac{1}{2\pi}\displaystyle{\int_{\RR^{2}}\left(\frac{1}{\left|x\right|^{2(N+2)-\beta_k}}e^{\hat{u}_{k}}+\frac{1}{\left|x\right|^{4-a\beta_k}}e^{a\hat{u}_{k}}\right)}
=\frac{1}{2\pi}\displaystyle{\int_{\RR^{2}}\left(\left|x\right|^{2N}e^{u_k}+e^{au_k}\right)}=\beta_k.
\end{cases}
\end{equation}

As above, we define the blow--up set $\hat{S}$ for $\hat{u}_k$ as given by
\begin{equation}\label{blowupseth}
\hat{S}=\left\{ \hat{x}_0: \exists \, \hat{x}_k \to \hat{x}_0 \mbox{ and } \hat{u}(\hat{x}_k) \to +\infty \right\}.
\end{equation}

Clearly $\hat{S}$ is also finite, and we have:

\begin{equation}\label{obseblow}
x_0\in S \mbox{ and } x_0\neq 0 \Leftrightarrow \frac{x_0}{\left|x_0\right|^{2}}\in\hat{S}.
\end{equation}

Obviously, Proposition \ref{1310} applies to $\hat{u}_k$ in any domain $\Omega \subset \mathbb{R}^2 \setminus \{ 0 \}$. Thus, in practice we are left to investigate the blow--up behavior for $u_k$ and $\hat{u}_k$ only when it occurs at the origin.

To this purpose, by arguing for $\hat{u}_k$ as in Corollary~\ref{coroc1}, we already know the following:

\begin{coro}\label{corostar}
Let $0\neq N>-1$ and $\min\left\{ 1,\frac{1}{N+1}\right\} < a \leq \max\left\{ 1,\frac{1}{N+1}\right\}$. If $\hat{u}_k$ in \eqref{kin0} admits a blow--up point at the origin and $\beta>\max\left\{2(N+1),\frac{2}{a}\right\}$ then for $\delta>0$ sufficiently small,
$$ \frac{1}{2\pi}\hat{f}_k\rightharpoonup\beta_{\infty}\delta_0\mbox{ weakly in the sense of measures in } B_{\delta},$$
with $\beta_{\infty}=2\beta-4\max\left\{ N+1,\frac{1}{a}\right\}$.
\end{coro}

\dimo
We shall prove that, if $N>0$ and $\frac{1}{N+1}<a<1$, then Proposition~\ref{knbl} applies to $\hat{u}_k$ with $U_k\equiv1, b=1, \alpha_k=\frac{\beta_k}{2}-(N+2) \to \frac{\beta}{2}-(N+2)>-1$ and $\sigma_k = \frac{e^{a\hat{u}_k}}{\vert x \vert^{4-a\beta_k}}$. While for $-1<N<0$ and $1<a<\frac{1}{N+1}$ then the Proposition~\ref{knbl} applies with $U_k\equiv 1$, $b=a$, $\alpha_k=-(2-\frac{a}{2}\beta_k)\to -(2-\frac{a}{2}\beta)>-1$ and $\sigma_k=\frac{e^{\hat{u}_{k}}}{\left|x\right|^{2(N+2)-\beta_k}}$. We discuss the first case only as the second one follows similarly. Indeed, for $r>0$ and 

\[
q=\left\{\begin{array}{lll}
\frac{1}{2-a(N+1)}  & \mbox{for }\frac{1}{N+1}<a\leq\frac{2}{N+2},\\
\frac{1}{a} & \mbox{for } a>\frac{2}{N+2},\\
\end{array}\right.
\]

we take $1<p<q$ and use Holder inequality to estimate

\begin{footnotesize}

\begin{equation}\label{richiamo}
\| \sigma_{k} \|_{L^p (B_r)}^{p} = \int_{B_r}\left(\frac{e^{\hat{u}_k}}{\left|x\right|^{\frac{4}{a}-\beta_k}}\right)^{ap}\leq\\
\left(\int_{B_r}\left|x\right|^{\beta_k-2(N+2)}e^{\hat{u}_k}\right)^{ap}\left(\int_{B_r}\left(\frac{1}{\left|x\right|^{\frac{4}{a}-2(N+2)}}\right)^\frac{ap}{1-ap}\right)^{1-ap} \leq C,
\end{equation}

\end{footnotesize}

for suitable $C>0$.
\qed

\begin{rem}\label{doublestar}
Notice that the assumption: $\beta> \max\left\{\frac{2}{a}, 2(N+1) \right\}$ in Corollary~\ref{corostar} is rather ''natural'' within the framework of the blow--up analysis developped in \cite{bm,bt,ls,barm}. On the other hand, from \eqref{toclim2}, we can guarantee such assumption only when $\frac{1}{2(N+1)}<a \neq \frac{1}{N+1} < \frac{2}{(N+1)}$. While in the other cases we can only ensure that: $\beta \geq \max\left\{\frac{2}{a}, 2(N+1) \right\}$. However, we shall show in Lemma \ref{lem42} how to handle the situation where, $\beta= \max\left\{\frac{2}{a}, 2(N+1) \right\}$. Actually, for $0 < a < \frac{1}{N+1}$ we show that $\beta = \frac{2}{a}$ can occur only for $a=\frac{1}{2(N+1)}$, see Lemma \ref{lem56}.
\end{rem}

In the remaining cases not covered by Corollary~\ref{coroc1} and Corollary~\ref{corostar} we shall describe the blow--up behavior for $u_k$ or $\hat{u}_k$ around the origin, by considering the following general problem:

\begin{equation}\label{starn}
\begin{cases}
-\Delta u_k=\left|x\right|^{2\alpha_{k,1}}V_{k}(x)e^{au_k}+\left|x\right|^{2\alpha_{k,2}}U_{k}(x)e^{u_k}=:g_k\quad\mbox{in }\Omega\\
\displaystyle{\int_{\Omega}\left(\left|x\right|^{2\alpha_{k,1}}V_{k}(x)e^{au_k}+\left|x\right|^{2\alpha_{k,2}}U_{k}(x)e^{u_k}\right)}\leq C,
\end{cases}
\end{equation}

supplemented by the condition:

\begin{equation}\label{starn2}
\displaystyle{\max_{\partial \Omega}u_k-\min_{\partial \Omega}u_k\leq C},
\end{equation}

with $\Omega \subset \mathbb{R}^2$ a regular open bounded domain. As for problem \eqref{lirich}, we assume:

\begin{equation}\label{limpofn}
V_k, U_k\in C^{0,1}(\Omega): 0<a_1\leq V_{k}(x)\leq b_1,\quad 0<a_2\leq U_{k}(x)\leq b_2\quad \forall \, x\in \Omega,
\end{equation}

\begin{equation}\label{limpof2n}
\left|\nabla V_{k}\right|+\left|\nabla U_{k}\right|\leq A\quad\mbox{in }\Omega,
\end{equation}

and

\begin{equation}\label{limpof3n}
\alpha_{k,1}\to\alpha_1>-1 \quad\mbox{ and }\quad \alpha_{k,2}\to\alpha_2>-1.
\end{equation}

By a direct extension of the arguments in \cite{bt}, we establish the following:

\begin{prop}\label{concefirst}
Let $\left\{u_k\right\}$ satisfy \eqref{starn} and \eqref{starn2} with $\Omega=B_{\delta}$ and assume \eqref{limpofn}, \eqref{limpof2n} and \eqref{limpof3n}. If $0<a\neq \frac{\alpha_1+1}{\alpha_2+1}$, and the origin is the only blow--up point for $u_k$ in $B_{\delta}$ then, along a subsequence, the following holds:  
\begin{equation}\label{concefirst1}
\frac{1}{2\pi}g_k\rightharpoonup\beta_0\delta_{0} \mbox{ weakly in the sense of measures in }B_{\delta},
\end{equation}
with $\beta_0=\beta_{0,1}+\beta_{0,2}$ such that:
\begin{equation}\label{concefirst2}
\beta_{0,1}:=\lim_{r\to0}\lim_{k\to+\infty}\frac{1}{2\pi}\int_{B_{r}}\left|x\right|^{2\alpha_{k,1}}V_{k}e^{au_k(x)}=\frac{a\beta_0(\beta_0-4(\alpha_2+1))} {4((\alpha_1+1)-a(\alpha_2+1))},
\end{equation}

\begin{equation}\label{concefirst3}
\beta_{0,2}:=\lim_{r\to0}\lim_{k\to+\infty}\frac{1}{2\pi}\int_{B_{r}}\left|x\right|^{2\alpha_{k,2}}U_{k}e^{u_k(x)}=\frac{\beta_0(4(\alpha_1+1)-a\beta_0)}{4((\alpha_1+1)-a(\alpha_2+1))}
\end{equation}
In particular,
\begin{equation}\label{concefirst4}
\displaystyle{\max_{K}u_k\to-\infty}, \quad \forall K \subset\subset B_{\delta} \setminus \left\{0 \right\},
\end{equation}
and
\begin{equation}\label{concefirst5}
\min\left\{\frac{4(\alpha_1+1)}{a}, 4(\alpha_2+1)\right\}\leq \beta_0 \leq \max\left\{\frac{4(\alpha_1+1)}{a}, 4(\alpha_2+1)\right\}.
\end{equation}
\end{prop}
\dimo

According to the given assumptions, there exists a finite measure $\mu$ in $B_{\delta}$ such that, along a subsequence, the following holds:
\begin{equation}\label{133star}
\frac{1}{2\pi} g_k\rightharpoonup\mu\mbox{ weakly in the sense of measures in }B_{\delta}. 
\end{equation}

Since, away from the origin, the sequence $u_k$ is uniformly bounded from above, by means of elliptic estimates we see that,

\begin{equation}\label{133star2}
\mu=2 \pi \beta_0\delta_0+\Phi, \hbox{ with } \Phi\in L^{1}(B_{\delta}),
\end{equation}

and

\begin{equation}\label{328bis}
\beta_0=\lim_{r\to0}\lim_{k\to+\infty}\frac{1}{2\pi}\int_{B_{r}}g_k(x) dx.
\end{equation}

Furthermore, as in Proposition~\ref{prop17} we easily check that necessarily,

\begin{equation}\label{133star3}
\beta_0\geq\min\left\{\frac{2(1+\alpha_1^-)}{a},2(1+\alpha_2^-)\right\},
\end{equation}

with $\alpha^-_j=\min\left\{0,\alpha_j\right\}$, $j=1,2$.

To check that actually we have $\Phi\equiv 0$ in \eqref{133star2}, we define:

$$\varphi_k=u_k-\inf_{\partial {B_{\delta}}}u_k,$$

satisfying:

\begin{equation}
\begin{cases}
-\Delta \varphi_k=g_k \quad\mbox{in }{B_{\delta}},\\
0\leq\varphi_k\leq C  \quad\mbox{on }{\partial {B_{\delta}}},
\end{cases}
\end{equation}

with suitable $C>0$. Therefore, we can use Green's representation formula for $\varphi_k$, and standard elliptic estimates, to check that along a subsequence, the following holds: 

\begin{equation}\label{GR}
\varphi_k\rightarrow \varphi=\beta_0{\rm log}\frac{1}{\left|x\right|}+\psi \quad \mbox{uniformly in }C^2_{\rm{loc}}({B_{\delta}}\setminus\left\{0\right\}), 
\end{equation}

with $\psi$ a smooth function in ${B_{\delta}}$.

In particular,

\begin{equation}\label{grrr2}
\nabla u_k=\nabla\varphi_k\rightarrow\nabla \varphi=-\beta_0\frac{x}{\left|x\right|^2}+\nabla\psi,
\end{equation} 

uniformly in $C^1_{\rm{loc}}({B_{\delta}}\setminus\left\{0\right\})$.

Analogously to \eqref{poho2}, for the solution $u_k$ of \eqref{starn} the following (local) Pohozaev's identity can be established in the usual way :

\begin{equation}\label{nuopo}
\begin{aligned}
&r\int_{\partial B_{r}}\left(\frac{\left|\nabla u_k\right|^{2}}{2}-\left(x\cdot\nabla u_k\right)^{2}\right)d\sigma=
\frac{r}{a}\int_{\partial B_{r}}\left|x\right|^{2\alpha_{k,1}}V_{k}(x)e^{au_k}d\sigma+r\int_{\partial B_{r}}\left|x\right|^{2\alpha_{k,2}}U_{k}(x)e^{u_k}d\sigma
\\
&-\frac{2(\alpha_{k,1}+1)}{a}\int_{B_{r}}\left|x\right|^{2\alpha_{k,1}}V_{k}(x)e^{au_k}-2(\alpha_{k,2}+1)\int_{B_{r}}\left|x\right|^{2\alpha_{k,2}}U_{k}(x)e^{u_k}\\
&-\frac{1}{a}\int_{B_{r}}\left(x\cdot\nabla V_{k}\right)\left|x\right|^{2\alpha_{k,1}}e^{au_k}-\int_{B_r}\left(x\cdot\nabla U_k\right)\left|x\right|^{2\alpha_{k,2}}e^{u_k}, \quad \mbox{for } 0<r\leq \delta.
\end{aligned}
\end{equation}

From \eqref{limpofn}, \eqref{limpof2n}, we see that,
\begin{equation}\label{1337bis}
\left|\int_{B_{r}}\left(x\cdot\nabla V_{k}\right)\left|x\right|^{2\alpha_{k,1}}e^{au_k}\right|\leq Cr\quad\mbox{ and }\left|\int_{B_r}\left(x\cdot\nabla U_k\right)\left|x\right|^{2\alpha_{k,2}}e^{u_k}\right|\leq Cr
\end{equation}
with a suitable constant $C>0$. Therefore, by using \eqref{grrr2} and \eqref{1337bis}, along a subsequence, we can pass to the limit in \eqref{nuopo} and obtain:
\begin{equation}\label{stermu}
\begin{aligned}
&\lim_{r\to 0}\lim_{k\to+\infty}\left(
\frac{r}{a}\int_{\partial B_{r}}\left|x\right|^{2\alpha_{k,1}}V_{k}(x)e^{au_k}d\sigma+r\int_{\partial B_{r}}\left|x\right|^{2\alpha_{k,2}}U_{k}(x)e^{u_k}d\sigma\right)\\
&=-\pi\beta^{2}_{0}+\frac{4\pi(\alpha_{1}+1)}{a}\beta_{0,1}+4\pi(\alpha_2+1)\beta_{0,2},
\end{aligned}
\end{equation}
with $\beta_{0,1}$ and $\beta_{0,2}$ defined in \eqref{concefirst2} and \eqref{concefirst3} respectively.

{\bf{CLAIM:}}
\begin{equation}\label{claimbis}
\displaystyle{\inf_{B_{\delta}} u_k = \inf_{\partial B_{\delta}} u_k \to - \infty  }
\end{equation}

To establish \eqref{claimbis}, we argue by contradiction and assume that $\displaystyle{\inf_{\partial B_{r}}u_k>-M}$. Then, by using \eqref{GR} and Fatou's Lemma, we obtain:
\begin{equation}\label{336bis}
C>\limsup_{k\to+\infty}\left(\int_{B_{r}}\left|x\right|^{2\alpha_{k,1}}V_{k}(x)e^{au_k}\right)\geq C_M
\limsup_{k\to+\infty}\int_{B_{r}}\left|x\right|^{2\alpha_{k,1}}e^{\varphi_k}\geq C\int_{B_r}\left|x\right|^{2\alpha_1-a\beta_0}e^{\psi},
\end{equation}
and consequently,
\begin{equation}\label{stimetaa}
0<\beta_{0}<\frac{2(1+\alpha_1)}{a}.
\end{equation}

Similarly, we check that,
\begin{equation}\label{stimetaa2}
0<\beta_{0}<2(1+\alpha_2).
\end{equation}

If $\alpha_j \leq 0$ for every $j=1,2$, then \eqref{stimetaa} and \eqref{stimetaa2} would already contradict \eqref{133star3}. On the other hand, if $\alpha_j > 0$ for some $j=1,2$, then from \eqref{stimetaa}, \eqref{stimetaa2} we see that,
\begin{equation}\label{propabar}
\lim_{r\to 0}\lim_{k\to+\infty}\left(\frac{r}{a}\int_{\partial B_{r}}\left|x\right|^{2\alpha_{k,1}}V_{k}(x)e^{au_k}d\sigma\right) = 0,
\end{equation}
and
\begin{equation}\label{propabar2}
\lim_{r\to 0} \lim_{k\to+\infty} \left( r\int_{\partial B_{r}}\left|x\right|^{2\alpha_{k,2}}U_{k}(x)e^{u_k}d\sigma\right)=0.
\end{equation}

Thus, from \eqref{stermu}, we find:
\begin{equation}\label{ellipse}
-\beta^{2}_{0}+\frac{4(\alpha_{1}+1)}{a}\beta_{0,1}+4(\alpha_2+1)\beta_{0,2}=0,
\end{equation}
with $\beta_0=\beta_{0,1}+\beta_{0,2}$.

But, in case $\alpha_{2}+1\geq\frac{\alpha_{1}+1}{a}$, then from \eqref{ellipse} we get
\begin{equation}\label{ellipse2}
\beta_0\left(\beta_0-\frac{4(\alpha_{1}+1)}{a}\right)=4\left[(\alpha_{2}+1)-\frac{\alpha_{1}+1}{a}\right]\beta_{0,2}\geq 0,
\end{equation}
that implies $\beta_0\geq\frac{4(\alpha_{1}+1)}{a}$, in contradiction to \eqref{stimetaa}.

Similarly, in case $\alpha_{2}+1<\frac{\alpha_{1}+1}{a}$, then \eqref{ellipse} implies:
\begin{equation}\label{ellipse3}
\beta_0\left(\beta_0-4(\alpha_{2}+1)\right)=4\left[\frac{\alpha_{1}+1}{a}-(\alpha_{2}+1)\right]\beta_{0,1}\geq 0,
\end{equation}
in contradiction to \eqref{stimetaa2}. So \eqref{claimbis} is established in any case. \\
Thus, from \eqref{GR} we get in particular that \eqref{133star2} holds with $\Phi \equiv 0 $, and also \eqref{concefirst1} and \eqref{concefirst4} are satisfied.
In addition, \eqref{claimbis} allows us to check that \eqref{propabar} and \eqref{propabar2} continue to hold, and consequently we can still guarantee the validity of \eqref{ellipse}. \\
At this point, \eqref{concefirst2} and \eqref{concefirst3} can be easily deduced.
\qed
\begin{rem}
It follows from the arguments above, that the assumption \eqref{limpof2n} can be relaxed around the origin, as pointed out in Remark~\ref{rem20bis}.
\end{rem}

By combining the results above, for problem \eqref{kappa} we have:

\begin{coro}\label{clocp0}
Assume \eqref{eq3.0},  \eqref{131bis} and let $\left\{u_k\right\}$ satisfy \eqref{kappa}. If {\bf {zero}} is a blow--up point for $\left\{u_k\right\}$, then (up to a subsequence) and for $\delta>0$ sufficiently small the following holds:
\begin{equation}
\frac{1}{2\pi}f_k\rightharpoonup\beta_0\delta_0\quad\mbox{ weakly in the sense of measure in } B_{\delta},
\end{equation}
with $\beta_0=\beta_{0,1}+\beta_{0,2}$, where
\begin{equation}\label{locp01b}
\beta_{0,1}:=\lim_{r\to0}\lim_{k\to+\infty}\frac{1}{2\pi}\int_{B_r}e^{au_k}=\frac{a\beta_0\left[\beta_0-4(N+1)\right]}{4\left[1-a(N+1)\right]},
\end{equation}
and
\begin{equation}\label{locp02b}
\beta_{0,2}:=\lim_{r\to0}\lim_{k\to+\infty}\frac{1}{2\pi}\int_{B_r}\left|x\right|^{2N}e^{u_k}=\frac{\beta_0\left[4-a\beta_0\right]}{4\left[1-a(N+1)\right]}.
\end{equation}

In particular,
\begin{equation}\label{346bis}
\min\left\{\frac{4}{a},4(N+1)\right\} \leq \beta_0 \leq \max\left\{\frac{4}{a},4(N+1)\right\},
\end{equation}
and
$$ \max_{K} u_k \to -\infty, \quad \forall K \subset \subset B_{\delta}\setminus\left\{ 0 \right\}, \mbox{ as } k\to+\infty $$

Furthermore,
\begin{equation}\label{locp03b}
\hbox{if } 0<a<\min\left\{1,\frac{1}{N+1}\right\}, \hbox{ then } \beta_0=4(N+1) \hbox{ and } \int_{B_{\delta}} e^{au_k} \to 0;
\end{equation}
and
\begin{equation}\label{locp04b}
\hbox{if } a>\max\left\{1,\frac{1}{N+1}\right\}, \hbox{ then } \beta_0=\frac{4}{a} \hbox{ and } \int_{B_{\delta}} \vert x \vert^{2N}e^{u_k} \to 0
\end{equation}
as $k\to\infty$.
\end{coro} 
\qed
\begin{rem}
Observe that for $N>0$, the conclusion \eqref{locp04b} holds also for $a=1$.
\end{rem}

Concerning the sequence $\hat{u}_k$ defined in \eqref{kin0} we find:

\begin{coro}\label{clocpinf}
Assume \eqref{eq3.0} and \eqref{131bis} and let $\hat{u}_k$ satisfy \eqref{kin}. If $\beta > \max\{ \frac{2}{a}, 2(N+1)\}$ and {\bf {zero}} is a blow--up point for $\left\{\hat{u}_k\right\}$ then, for $\delta>0$ sufficiently small, the following holds (along a subsequence):
\begin{equation}
\frac{1}{2\pi}\hat{f}_k\rightharpoonup\beta_{\infty}\delta_0\quad\mbox{ weakly in the sense of measures in } B_{\delta},
\end{equation}
\begin{equation}\label{locp01bi}
\beta_{1,\infty}:=\lim_{r\to0}\lim_{k\to+\infty}\frac{1}{2\pi}\int_{B_r}\frac{1}{\left|x\right|^{4-a\beta_k}}e^{a\hat{u}_{k}}=\frac{a\beta_{\infty}\left[2\beta-4(N+1)-\beta_{\infty}\right]}{4\left[1-a(N+1)\right]}
\end{equation}
and
\begin{equation}\label{locp02bi}
\beta_{2,\infty}:=\lim_{r\to0}\lim_{k\to+\infty}\frac{1}{2\pi}\int_{B_r}\frac{1}{\left|x\right|^{2(N+2)-\beta_k}}e^{\hat{u}_{k}}=\frac{\beta_{\infty}\left[4+a\beta_{\infty}-2a\beta\right]}{4\left[1-a(N+1)\right]}.
\end{equation}
In particular,
\begin{equation}\label{eq359bis}
\beta_{\infty}\geq 2\left(\beta-\max\left\{\frac{2}{a}, 2(N+1)\right\}\right)>0 \hbox{ and } \max_{K}\hat{u}_k \to -\infty, \forall K \subset \subset B_{\delta}\setminus\left\{0 \right\}.
\end{equation}
Furthermore, 
\begin{equation}\label{eq360bis}
\mbox{if } \min\left\{1,\frac{1}{N+1}\right\} \leq a \leq \max\left\{1,\frac{1}{N+1}\right\}, \mbox{ then }
\beta_{\infty}=2\left(\beta-\max\left\{ \frac{2}{a}, 2(N+1)\right\} \right);
\end{equation}
and so,
\begin{equation}\label{locp03bi}
\hbox{for } N>0, \mbox{ we have: } \, \int_{\vert x \vert \leq \frac{1}{\delta}} e^{au_k} = \int_{\vert x \vert \leq \delta} \frac{1}{\vert x \vert^{4-a\beta_k}} e^{a\hat{u}_k} \to 0,
\end{equation}
while,
\begin{equation}\label{locp04bi}
\hbox{for } -1<N\leq0 \mbox{ we have: } \, \int_{\vert x \vert \leq \frac{1}{\delta}} \vert x \vert^{2N} e^{u_k} = \int_{\vert x \vert \leq \delta} \frac{1}{\vert x \vert^{2(N+2)-\beta_k}} e^{\hat{u}_k} \to 0,
\end{equation}
as $k\to\infty$.
\end{coro}

As a useful consequence of Proposition~\ref{1310}, Corollary~\ref{coroc1} and Corollary~\ref{clocpinf}, we have:

\begin{lemma}\label{lemma307}
Assume \eqref{eq3.0}, \eqref{kappa} and \eqref{131bis}. If $\beta > \max\{ \frac{2}{a}, 2(N+1)\}$ and $S \cup \hat{S}\neq \emptyset$, then, as $k\to +\infty$, the following holds:
\begin{enumerate}[(i)]
\item If $S=\emptyset$ then $\displaystyle{\max_{\vert x \vert \leq R}u_k\to -\infty}$, $\forall R>0$.
\item If $0 \notin \hat{S}$ then $\displaystyle{\max_{\vert x \vert \leq r}\hat{u}_k\to -\infty}$, and \\ $\displaystyle{\lim_{k \to \infty} \int_{\vert x \vert \geq \frac{1}{r}} f_k(x) dx= \lim_{k \to \infty} \int_{\vert x \vert \leq r} \hat{f}_k(x) dx}=0$, $\forall r>0$ sufficiently small.
\item For every $K \subset \subset \mathbb{R}^2\setminus S$ we have $\displaystyle{\max_{K}u_k\to -\infty}$.
\end{enumerate}
\end{lemma}
\qed

%
%

We point out the following version of Lemma~\ref{lemma307}, which covers the case where $\beta=\max \left\{\frac{2}{a},2(N+1)\right\}$.

\begin{lemma}\label{lem42}
Assume \eqref{eq3.0} and let $u_k$ satisfy \eqref{kappa} and \eqref{131bis}. If $\beta=\max \left\{\frac{2}{a},2(N+1)\right\}$ and we set,

\begin{equation}\label{4.2.0}
\beta_{\infty}:=\lim_{R\to +\infty} \lim_{k \to + \infty} \frac{1}{2\pi} \int_{\vert x \vert \geq R} f_{k}(x) \, dx 
\end{equation}

then the following holds:

\begin{enumerate}[(i)]

\item $0 < a \leq \frac{1}{2(N+1)}$ or $a \geq \frac{2}{N+1}$.

\item For every $0<\varepsilon<r$ sufficiently small,

\begin{equation}\label{eq400}
\displaystyle{\max_{\varepsilon\leq \vert x \vert \leq r} \hat{u}_k \to -\infty} \mbox{ and in particular } \displaystyle{\inf_{B_r} \hat{u}_k = \inf_{\partial B_r} \hat{u}_k \to -\infty} \mbox{ as } k\to+\infty.
\end{equation}

Moreover, part $(i)$ and $(iii)$ of Lemma \ref{lemma307} continue to hold.

\item 

\begin{equation}\label{eq401}
\begin{split}
&\mbox{if } 0<a\leq\frac{1}{2(N+1)} \mbox{ and } 0 \notin \hat{S} \text{ or more generally } \beta_{\infty} < 2(1+N^-),\text{ with } N^-=\min \left\{ 0, N \right\}, \\
&\displaystyle{\text{ then } \lim_{r\to 0}\lim_{k\to +\infty} \int_{\vert x \vert\geq \frac{1}{r}} \vert x \vert^{2N} e^{u_k} =\lim_{r\to 0}\lim_{k\to +\infty} \int_{\vert x \vert\leq r} \vert x \vert^{\beta_k-2(N+2)} e^{\hat{u}_k}=0};
\end{split}
\end{equation}

\begin{equation}\label{eq402}
\begin{split}
&\mbox{ if either } N>1 \text{ and } \frac{2}{N+1} \leq a < 1,  \text{ or } 1 \ne a \geq \max \left\{ 1, \frac{2}{N+1} \right\} \text{ and } 0 \notin \hat{S} \text{ or  } \beta_\infty < \frac{2}{a},\\
&\mbox{ then } \displaystyle{\lim_{r\to 0}\lim_{k\to +\infty} \int_{\vert x \vert\geq \frac{1}{r}} e^{au_k} =\lim_{r\to 0}\lim_{k\to +\infty} \int_{\vert x \vert\leq r} \vert x \vert^{a\beta_k-4} e^{a\hat{u}_k}=0}
\end{split}
\end{equation}

\end{enumerate}

\end{lemma}

\dimo
Clearly, part $(i)$ is a direct consequence of \eqref{toclim2}. Concerning $(ii)$ and $(iii)$,  we only consider the case where:

\begin{equation}\label{preq403}
0<a\leq\frac{1}{2(N+1)} \mbox{ and } \beta_k \to \beta=\frac{2}{a} \mbox{ as } k\to+\infty,
\end{equation}

as the other cases follow in a similar way.

By virtue of the estimate:

\begin{equation}\label{preq404}
\displaystyle{e^{a\inf_{B_r}\hat{u}_k}} \int_{B_r} \vert x \vert^{a\beta_k-4} \leq \int_{B_r} \vert x \vert^{a\beta_k-4} e^{a\hat{u}_k} \leq C,
\end{equation}

from \eqref{preq403}, we deduce that necessarily, $\displaystyle{\inf_{B_r} \hat{u}_k \to -\infty}$, as $k\to+\infty$. On the other hand, for $r>0$ sufficiently small and $\varepsilon \in (0,r)$ we find a constant $C_\varepsilon > 0$: $\| \hat{u}_{k} \|_{L^\infty (B_r \setminus B_\varepsilon)} \leq C_\varepsilon$. 

Therefore, by means of Harnack inequality (e.g. see Proposition~5.2.8. of \cite{tar2}) we deduce that: $\displaystyle{\max_{\varepsilon\leq \vert x \vert \leq r} \hat{u}_k \to -\infty}$, as $k\to+\infty$, and $(ii)$ is established.

To establish $(iii)$, we simply observe that, if $0 \notin \hat{S}$ then \eqref{eq401} follows quite easily. While under the assumption: $\beta_{\infty} < 2(1+N^-)$, we can use Lemma \ref{lem21} together with Remark \ref{rem21} for $\hat{u}_k$ in $B_r$, and conclude that, 

\begin{equation}\label{4.6bis}
\| e^{\hat{u}_k} \|_{L^p (B_r)} \leq C,
\end{equation}

for suitable $p>1$ and $C>0$.

Hence by letting $q=\frac{p}{p-1}$, we can estimate:

\begin{equation}\label{4.6.1}
\int_{\vert x \vert\leq r} \vert x \vert^{\beta_k-2(N+2)} e^{\hat{u}_k}\leq C \left( \int_0^{r} t^{(\beta_k-2(N+2))q+1} \, dt \right)^{\frac1q} = C \left( \frac{r^{l_k}}{l_k} \right)^{\frac1q} , 
\end{equation}

with $l_k = \beta_k - 2 (N+1) + (\beta_k - 2(N+2))\frac{1}{p-1} \to (\frac2a - 2 (N+1)) + (\frac2a - 2(N+2))\frac{1}{p-1}>0$, as $k\to+\infty$.  So, from \eqref{4.6.1} we deduce  \eqref{eq401}.\\
\qed

In concluding this section, we point out some useful generalizations of the results stated above, which allow us to account for the scaling properties of problem \eqref{kappa} under either one of the following transformations:

$$
u_k(x) \mapsto u_k(Rx)+\frac{2}{a}\log R \quad \mbox{ or } \quad u_k(x) \mapsto u_k(Rx)+2(N+1)\log R,
$$
$R>0$. 

To this purpose we replace \eqref{limpofn} respectively with the following assumptions:

\begin{equation}\label{conclu1}
V_k(x)=\varepsilon_{2,k}V_{1,k}(x) \mbox{ with } U_k \mbox{ and } V_{1,k} \mbox{ satisfying } \eqref{limpofn}, \eqref{limpof2n} \mbox{ and } \lim_{k \to \infty} \varepsilon_{2,k} = 0,
\end{equation}
or

\begin{equation}\label{conclu3}
U_k(x)=\varepsilon_{1,k}U_{1,k}(x) \mbox{ with } U_{1,k} \mbox{ and } V_k \mbox{ satisfying } \eqref{limpofn}, \eqref{limpof2n} \mbox{ and } \varepsilon_{2,k}\to 0,
\end{equation}

For example, we can anticipate that the conclusion of Proposition~\ref{concefirst} remains valid if \eqref{conclu1} holds with:

\begin{equation}\label{conclu2}
 \lim_{k\to+\infty} \frac{2(\alpha_{1,k}+1)}{a}:=\frac{2(\alpha_1+1)}{a} \geq \alpha_2+1:=\lim_{k\to+\infty}\alpha_{2,k}+1>0
\end{equation}

or \eqref{conclu1} holds with:

\begin{equation}\label{conclu4}
\lim_{k\to+\infty} \alpha_{2,k}+1 :=\alpha_2+1 \geq \frac{\alpha_1+1}{2a}:=\frac{1}{2a}\lim_{k\to+\infty}\left(\alpha_{1,k}+1\right)>0
\end{equation}

In this direction we have:

\begin{prop}\label{prop29}
Let $u_k$ satisfy \eqref{starn} with $a>0$. Assume \eqref{limpof3n} and that either \eqref{conclu1} or \eqref{conclu3} hold. If $x_0$ is a blow--up point for $u_k$ in $\Omega$, then,

\begin{equation}\label{eq357}
\beta(x_0):=\lim_{r\to 0} \lim_{k\to + \infty} \frac{1}{2\pi}\int_{B_r(x_0)}g_k(x) dx \geq \min \left\{ 2(N^- +1), \frac{2}{a} \right\}.
\end{equation}

So the blow--up set of $u_k$ in $\Omega$ contains at most finitely many points.

Furtermore, if we assume \eqref{starn2} then the following holds:

a) for $x_0 = 0$ and $a \ne \frac{\alpha_1+1}{\alpha_2+1}$, the identities \eqref{concefirst2}, \eqref{concefirst3} and \eqref{concefirst5} remain valid with $\beta_0 = \beta (x_0)$;

b) for $x_0 \ne 0$, $a \ne 1$ and $\beta_0 = \beta (x_0)$ we have:
\begin{equation}\label{locp11b2}
\lim_{r\to0}\lim_{k\to+\infty}\frac{1}{2\pi}\int_{B_r(x_0)}\vert x\vert^{2\alpha_{k,1}} V_k(x)e^{au_k}dx=\frac{a\beta_0\left(\beta_0-4\right)}{4\left(1-a\right)},
\end{equation}
\begin{equation}\label{locp12b2}
\lim_{r\to0}\lim_{k\to+\infty}\frac{1}{2\pi}\int_{B_r(x_0)}\vert x\vert^{2\alpha_{k,1}} U_k(x)e^{u_k}dx=\frac{\beta_0\left(4-a\beta_0\right)}{4\left(1-a\right)}.
\end{equation}
and in particular,
\begin{equation}\label{eq326bis2}
\frac{4}{\max \left\{1,a \right\}} \leq \beta_0 \leq \frac{4}{\min\left\{1,a \right\}}.
\end{equation}
\end{prop}

\dimo
First of all we observe that, by virtue of Lemma~\ref{lem21} and Remark~\ref{rem21} the weaker assumption \eqref{conclu1} or \eqref{conclu3} still suffice to ensure that \eqref{eq357} holds and the blow--up set of $u_k$ in $\Omega$ contains at most finitely many points. Hence, for $x_0=0$, we can follow step by step the proof of Proposition~\ref{concefirst} to arrive at the identity \eqref{nuopo}. At this point, we can still deduce \eqref{stermu}. Indeed, if we suppose for example that \eqref{conclu3} holds, then \eqref{propabar2} is now a direct consequence of the fact that: $\Vert U_k \Vert_{L^{\infty}(\Omega)} \to 0$, as $k\to +\infty$. While,  as in the proof of Proposition~\ref{concefirst}, the identity \eqref{propabar} can be derived regardless of the validity of  \eqref{claimbis}. In other words, under the given assumptions we can always deduce \eqref{ellipse} and the desired identities \eqref{concefirst2} and \eqref{concefirst3} follow in this case. In an analogous way, we may handle the case where $x_0\neq 0$. Indeed, after a translation, we can use the same arguments of Proposition \ref{concefirst}, which now we may apply with $\alpha_{kj} \equiv 0$, $j=1,2$, in order to derive \eqref{locp11b2}-\eqref{locp12b2}. Similarly we treat the case where \eqref{conclu1} holds.\\
\qed

\begin{rem}
The case $a=1$ can be handled directly by Proposition \ref{knbl} and yields: $\beta_0 = 4$.
\end{rem}

Notice that Proposition~\ref{prop29} does not provide any information about the "concentration" properties of $u_k$ around its blow--up points, in the sense that the limiting measure $\mu$ in \eqref{133star} may admit the decomposition \eqref{133star2} with $\Phi \not\equiv 0$.

On the other hand, by keeping on following the arguments of the proof of Proposition~\ref{concefirst}, we can easily check that, under the weaker assumption \eqref{conclu1} or \eqref{conclu3} the following holds:

\begin{prop}\label{prop33}
Under the assumptions of Proposition~\ref{prop29}, if $x_0=0$, $a\neq \frac{\alpha_1+1}{\alpha_2+1}$ and we suppose that:
\begin{enumerate}[(i)]
\item \begin{equation}\label{eq368}
\eqref{conclu3} \mbox{ hold with } \beta(x_0)\geq\frac{2(\alpha_1+1)}{a},
\end{equation} or
\item \begin{equation}\label{eq369} 
\eqref{conclu1} \mbox{ hold with } \beta(x_0)\geq 2(\alpha_2+1),
\end{equation}
then the full conclusion of Proposition~\ref{concefirst} holds. In particular, ''concentration'' occurs in the sense that, \eqref{concefirst1} and \eqref{concefirst4} hold.
\end{enumerate}
\end{prop}

Similarly, for a blow--up point other than the origin we have:

\begin{prop}\label{prop33bis}
Under the assumption of Proposition~\ref{prop29}, let $x_0\neq 0$ and $a\neq 1$. We have:

\begin{enumerate}[(i)]
\item \begin{equation}\label{eq370}
if \eqref{conclu3} \mbox{ holds and } \beta(x_0)\geq\frac{2}{a},
\end{equation} 

or 

\item \begin{equation}\label{eq371} 
if \eqref{conclu1} \mbox{ holds and } \beta(x_0)\geq 2,
\end{equation}

 then, for $\delta>0$ sufficiently small, there holds:

\begin{equation*}
\frac{1}{2\pi} g_k \rightharpoonup \beta_0 \delta_{x_0} \text{ weakly in the sense of measures in } B_\delta (x_0)
\end{equation*}

with $\beta_0 = \beta (x_0)$ satisfying \eqref{locp11b2}, \eqref{locp12b2}.

In particular:
\begin{equation}\label{locp13b2}
\max_{K} u_k \to -\infty, \quad \forall K \subset \subset B_{\delta}(x_0)\setminus\left\{ x_0 \right\}
\end{equation}

\end{enumerate}
\end{prop}
As a useful consequence of Proposition~\ref{prop33} we find:

\begin{coro}\label{coro31}
Under the assumptions of Proposition~\ref{prop29}, suppose that either \eqref{conclu1} and \eqref{conclu2} hold or that \eqref{conclu3} and \eqref{conclu4} hold, then the full conclusion of Proposition~\ref{concefirst} holds for $x_0=0$ and $a \ne\frac{\alpha_1+1}{\alpha_2+1}$. 
\end{coro}

\dimo
We establish the desired conclusion when \eqref{conclu3}-\eqref{conclu4} hold; as the other case where \eqref{conclu1}-\eqref{conclu2} holds,  follows similarly. Indeed we check that, by the given assumptions, we have: $\beta(x_0) \geq\frac{2(\alpha_1+1)}{a}$, and then the desired conclusion follows by Proposition~\ref{prop33}. Indeed, on the basis of \eqref{conclu4} we distinguish the following two cases:

\begin{equation}\label{eq358}
\frac{\alpha_1+1}{2a}\leq \alpha_2+1<\frac{\alpha_1 +1 }{a} \quad \mbox{ or } \quad \alpha_2+1\geq \frac{\alpha_1+1}{a}
\end{equation}

By Proposition~\ref{prop33}, if the first inequality holds in \eqref{eq358} then we can use \eqref{concefirst2} to get $\beta_0\geq 4(\alpha_2+1)\geq \frac{2(\alpha_1+1)}{a}$. When the second inequality holds in \eqref{eq358} then we can use \eqref{concefirst3} to get the stronger inequality: $\beta_0\geq\frac{4(\alpha_1+1)}{a}$, and the desired lower bound for $\beta(x_0)$ follows in any case.
\qed

Similarly, we obtain:

\begin{coro}\label{coro34bis}
Under the assumption of Proposition~\ref{prop29}, suppose that $x_0\neq 0$ and $a\neq 1$. Assume that \eqref{conclu1} holds with $a\leq 2$ or that \eqref{conclu3} holds with $a \geq \frac{1}{2}$, then the conclusion of Proposition~\ref{prop33bis} holds.
\end{coro}
\qed

Finally, we describe what happens in situations somewhat complementary to those covered by Proposition~\ref{prop33} and Corollary~\ref{coro31}, or Proposition~\ref{prop33bis} and Corollary~\ref{coro34bis}.

\begin{prop}\label{prop35}
Under the assumptions of Proposition~\ref{prop29} the following holds:
\begin{enumerate}[(i)]
\item If $x_0=0$ and
\begin{enumerate}[(a)]
\item \eqref{conclu3} holds with $\alpha_1 > -1$ and $\beta(x_0)< \frac{2}{a} (1 + \alpha_1^-)$, with $\alpha_1^- = \min \{ 0 , \alpha_1 \}$, then $\beta(x_0)=4(\alpha_2+1)$ and necessarily $\alpha_2+1<\frac{1}{2a} (1 + \alpha_1^-)$, 
\item \eqref{conclu1} holds with $\alpha_2 > -1$ and $\beta(x_0)<2 (1 + \alpha_2^-)$, with $\alpha_2^- = \min \{ 0 , \alpha_1 \}$, then $\beta(x_0)=\frac{4}{a}(\alpha_1+1)$ and necessarily $\alpha_1+1<\frac{a}{2} (1 + \alpha_2^-)$.
\end{enumerate}

\item If $x_0\neq0$ and 
\begin{enumerate}[(a)]
\item \eqref{conclu3} holds and $\beta(x_0)<\frac{2}{a}$, then $\beta(x_0)=4$ and necessarily $0<a<\frac12$;
\item \eqref{conclu1} holds and $\beta(x_0)<2$, then $\beta(x_0)=\frac{4}{a}$ and necessarily $a>2$.
\end{enumerate}
\end{enumerate}
\end{prop}

\dimo
As usual we check $(i)$-$(a)$, as the other cases follow similarly. Indeed, by the given assumption, for $\varepsilon>0$ sufficiently small there exists $k_{\varepsilon}\in\mathbb{N}$ and $r_{\varepsilon}>0$:

$$
\frac{1}{2\pi}\int_{B_r}g_k(x)<\frac{2}{a} (1 + \alpha_1^-) -\varepsilon \quad \forall k\geq k_{\varepsilon}, \, \forall r \in (0,r_{\varepsilon}).
$$

Thus in view of Remark~\ref{rem21}, we can apply Lemma~\ref{lem21} for $u_k$ to check that, for a suitable $p>\frac{1}{(1 + \alpha_1^-)}$ and $C>0$, there holds:

$$
\Vert e^{au_k}\Vert_{L^p(B_r)} \leq C.
$$

As a consequence, we must have: $\beta_{0,1}=0$, and from Proposition~\ref{prop29} we can use \eqref{concefirst2} to conclude that, $\beta_0=\beta(x_0)=4(\alpha_2+1)$ as claimed.
\qed

\section{Preliminary results}
\setcounter{equation}{0}

In this section we show how to use the (local) blow--up analysis established above in order to obtain some useful information about the blow--up sets $S$ and $\hat{S}$ of $u_k$ and $\hat{u}_k$ respectively.

\begin{rem}\label{rem4.0}
Although not always specified, it is understood that most of the limits taken as $k \to \infty$, generally hold along a subsequence.
\end{rem}

From now on we shall suppose that,

\begin{equation}\label{cond305}
\begin{cases}
u_k \mbox{ satisfies } \eqref{kappa} \mbox{ and } \eqref{131bis}, \\ 
N>0, 0 < a \neq \frac{1}{N+1}, a \ne 1 \mbox{ and } S\cup \hat{S}\neq \emptyset.
\end{cases}
\end{equation}

We point out the following simple properties: 

\begin{lemma}\label{2app}
Assume \eqref{cond305}. The following holds:
\begin{itemize}
\item[(i)] \begin{equation}\label{12app}
\mbox{if } 0<a<\frac{1}{N+1}, \mbox{ then } \int_{\left\{\left|x\right|\leq R\right\}}e^{au_k}\to0, \quad\mbox{as }k\to+\infty,\quad \forall R>0;
\end{equation}
\item[(ii)] \begin{equation}\label{13app}
\mbox{if }\frac{1}{N+1}<a<1, \mbox{ then }\int_{\left\{\left|x\right|\geq \varepsilon\right\}}e^{au_k}\to0, \quad\mbox{as }k\to+\infty,\quad \forall \varepsilon>0;
\end{equation} 
\item[(iii)]
\begin{equation}\label{14app}
\mbox{if } a>1, \mbox{ then }\int_{\left\{\left|x\right|\leq R\right\}}\left|x\right|^{2N}e^{u_k}\to0, \quad\mbox{as }k\to+\infty,\quad \forall R>0.
\end{equation} 
\end{itemize}
\end{lemma}

\dimo
If $S=\emptyset$ then both \eqref{12app} and \eqref{14app} follows directly from part $(i)$ of Lemma~\ref{lemma307} (or part (ii) of Lemma~\ref{lem42}). On the other hand, if $S\neq \emptyset$ then for $0<a<\frac{1}{N+1}$, we can use the estimate in \eqref{corohom} in order to deduce that, for $x_0 \in S$ we have:

\begin{equation}\label{eq46}
\lim_{r\to 0}\lim_{k\to +\infty} \int_{B_r(x_0)}e^{au_k}=0.
\end{equation}

Similarly for $a>1$ we can check that if $x_0 \in S \ne \emptyset$ then 

\begin{equation}\label{eq47}
\lim_{r\to 0}\lim_{k\to +\infty} \int_{B_r(x_0)}\vert x \vert^{2N}e^{u_k}=0.
\end{equation}

Therefore \eqref{12app} and \eqref{14app} follow as a consequence of \eqref{eq46}, \eqref{eq47} and part $(iii)$ of Lemma~\ref{lemma307} (or part $(ii)$ of Lemma~\ref{lem42}). Concerning \eqref{13app} we argue similarly, and observe that, in case $S\setminus \left\{0\right\} \ne \emptyset$, then \eqref{eq46} still holds for every $x_0\in S\setminus \left\{0\right\}$. Therefore, as above we can use part $(iii)$ of Lemma~\ref{lemma307} (or part $(ii)$ of Lemma~\ref{lem42}) to obtain: 

\begin{equation}\label{eq3011star}
\displaystyle{\int_{\varepsilon\leq\vert x \vert\leq R}e^{au_k}\to0}, \quad \mbox{ as } k\to+\infty,
\end{equation} 
$\forall \, 0<\varepsilon<R$.

On the other hand, for $\frac{1}{N+1}<a<1$ and $\vert x \vert \geq R$ we can estimate:

\begin{equation}\label{richiamo}
\begin{split}
&\int_{\vert x \vert \geq R}e^{au_k}=\int_{\vert x \vert \leq 1/R}\left(\frac{e^{\hat{u}_k}}{\left|x\right|^{4/a-\beta_k}}\right)^{a}\\
&\leq\left(\int_{\vert x \vert \leq 1/R}\left|x\right|^{\beta_k-2(N+2)}e^{\hat{u}_k}\right)^{a}\left(\int_{\vert x \vert \leq 1/R}\left(\frac{1}{\left|x\right|^{4/a-2(N+2)}}\right)^\frac{a}{1-a}\right)^{1-a}\\
&\leq C\left(\int_{\vert x \vert \leq 1/R}\left(\frac{1}{\left|x\right|^{4-2a(N+2)}}\right)^{\frac{1}{1-a}}\right)^{1-a} \to 0,\quad\mbox{ as }R\to+\infty.
\end{split}
\end{equation}

Therefore, \eqref{13app} follows by taking into account \eqref{eq3011star} and \eqref{richiamo}.\\
\qed

We start to analyze the behavior of $u_k$ when $S\neq\emptyset$. To this purpose, we let:

$$
S=\left\{z_1,...,z_m \right\} \subset \mathbb{R}^2, \mbox{ with } z_i\neq z_j \mbox{ for } i\neq j,
$$
and as before we set,
\begin{equation}\label{410*}
\beta (z_j) := \; \displaystyle{\lim_{r\to 0} \liminf_{k \to + \infty} \left( \frac{1}{2\pi} \int_{B_r(z_j)} f_k(x) dx\right)}.
\end{equation}

According to Corollary~\ref{coroc1} and \ref{clocp0}, along a subsequence, as $k \to\infty$ we have:

\begin{equation}\label{eq411}
f_k \rightharpoonup \sum_{j=1}^{m} \beta(z_j)\delta_{z_j}, \quad\mbox{ weakly in the sense of measures};
\end{equation}
locally in $\mathbb{R}^2$,  and with $\beta(z_j)=\frac{4}{\max\left\{1,a\right\}}$ for $z_j\neq 0$. While, if $z_j=0$ for some $j=1,...m$, then we  set: $\beta_0:=\beta(0)$, and the value of $\beta_0$ is characterized by the properties: \eqref{locp01b}, \eqref{locp02b}, \eqref{346bis}, \eqref{locp03b} and \eqref{locp04b} as given in Corollary~\ref{clocp0}.

As a consequence of \eqref{eq411} we have:

\begin{lemma}\label{gradcol}
Let $\left\{u_k\right\}$ satisfy \eqref{kappa}, and assume \eqref{cond305}. If $S=\left\{z_1,...,z_m\right\}$, then 
\begin{equation}\label{gradcve}
\nabla u_{k}(x)\to\sum^{m}_{j=1}\beta_{j}\frac{x-z_j}{\left|x-z_j\right|^{2}},\quad\mbox{as }k\to+\infty,
\end{equation}
uniformly in $C_{{\rm{loc}}}(\RR^{2}\setminus S)$, with $\beta_{j}=\beta(z_j)$ as specified in Corollary~\ref{coroc1} and Corollary~\ref{clocp0}.
\end{lemma}

\dimo

We need to show that, for a compact set $K\subset\subset \RR^{2}\setminus S$, we have:

\begin{equation}\label{1423bis}
\max_{x\in K}\left|\nabla u_k(x)-\beta_j\sum^{m}_{j=1}\frac{x-z_j}{\left|z-z_j\right|^{2}}\right|\to0\quad\mbox{as }k\to+\infty.
\end{equation}

To establish \eqref{1423bis}, we take $\delta>0$ sufficiently small, so that $S \cap B_{2\delta}(z_j)=\left\{ z_j\right\}$ and $B_{2\delta}(z_j)\cap B_{2\delta}(z_i)=\emptyset$ for $i\neq j$. Let $N_{\delta}:=\cup^{m}_{j=1}B_{\delta}(z_j)$, and we choose $R>\delta>0$ sufficiently large so that, $N_{\delta}\subset \subset B_{R/2}$ and 
$K\subset\subset B_{R/2}\setminus N_{\delta}.$
For $x\in K$, we write:
\begin{equation}\label{1412bis}
\begin{aligned}
\nabla u_k(x)&=\frac{1}{2\pi}\int_{B_{R}\setminus N_{\delta}}\frac{x-y}{\left|x-y\right|^{2}}f_{k}(y)dy+\frac{1}{2\pi}\sum^{m}_{j=1}\int_{B_{\delta}(z_j)}\frac{x-y}{\left|x-y\right|^{2}}f_{k}(y)dy\\&
+\frac{1}{2\pi}\int_{\left\{\left|x\right|\geq R\right\}}\frac{x-y}{\left|x-y\right|^{2}}f_{k}(y)dy=:J_{1}(x)+\sum^{m}_{j=1}I_{j}(x)+J_{2}(x),
\end{aligned}
\end{equation}
and we proceed to estimate each of the terms above.

For the first term we have:

\begin{equation}\label{1442bis}
\left|J_{1}(x)\right|\leq\left(\sup_{B_{R}\setminus N_{\delta}}f_{k}\right)\int_{B_{2R}(x)}\frac{1}{\left|x-y\right|}dy=CR^{2}\left(\sup_{B_{R}\setminus N_{\delta}}f_{k}\right)\to0,\quad\mbox{as }k\to+\infty.
\end{equation}

For $j=1,...,m$, we write:

\begin{equation}\label{grapuf1}
\begin{aligned}
I_{j}(x)-\beta_j\frac{x-z_j}{\left|x-z_j\right|^{2}}&=\frac{1}{2\pi}\int_{B_{\delta}(z_j)}\left(\frac{x-y}{\left|x-y\right|^{2}}-\frac{x-z_j}{\left|x-z_j\right|^{2}}\right)f_{k}(y)dy\\
&+\frac{1}{2\pi}\int_{B_{\delta}(z_j)}\frac{x-z_j}{\left|x-z_j\right|^{2}}\left|\beta_{k,j}-\beta_{j}\right|
\end{aligned}
\end{equation}

with $\beta_{k,j}:=\frac{1}{2\pi}\displaystyle{\int_{B_{\delta}(z_j)}f_k\to\beta_j}$, as $k\to+\infty.$
Observing that for $x\in K$, we have: 

$$ V_j(y):=\left|\frac{x-y}{\left|x-y\right|^{2}}-\frac{x-z_j}{\left|x-z_j\right|^{2}}\right|\in C^{0}(B_{\delta}(z_j)),\quad\mbox{ and } V_j(z_j)=0,$$

then from \eqref{eq411} and \eqref{grapuf1}, we find:

\begin{equation}\label{grapuf2}
\left|I_{j}(x)-\beta_j\frac{x-z_j}{\left|x-z_j\right|^{2}}\right|\leq\frac{1}{2\pi}\int_{B_{2\delta}(z_j)}\left|\frac{x-y}{\left|x-y\right|^{2}}-\frac{x-z_j}{\left|x-z_j\right|^{2}}\right|f_{k}(y)dy+\frac{1}{2\delta}\left|\beta_{k,j}-\beta_j\right|\to0,
\end{equation}

as $k\to+\infty$.

Finally, for the last term in \eqref{1412bis}, we find:

\begin{equation}\label{1417bis}
\left|J_{2}(x)\right|=\int_{\left\{\left|y\right|\geq R\right\}}\frac{1}{\left|x-y\right|}f_{k}(y)dy
\leq \frac{C}{R}\int_{\left\{\left|y\right|\geq R\right\}} f_{k}(y)dy\leq\frac{C}{R},
\end{equation}

for every $R$ sufficiently large, and suitable $C>0$. Thus, by combining \eqref{1442bis}, \eqref{grapuf2} and \eqref{1417bis}, we deduce \eqref{1423bis}.
\qed

\begin{lemma}\label{lem3024bis}
Suppose that $S \setminus \{ 0 \} \ne \emptyset$ and let $\left\{z_1,...,z_n \right\}=S\setminus\left\{ 0\right\}$ with $z_i\neq z_j$ for $i\neq j\in \left\{1,...,n\right\}$. The following holds:
\begin{enumerate}[(i)]
\item If $0<a\neq\frac{1}{N+1}<1$ then

\begin{equation}\label{eq36}
(2N-\beta_0)\frac{z_i}{\vert z_i\vert^2}=4\displaystyle{\left(\sum_{i\neq j =1}^{n}\frac{z_i-z_j}{\left|z_i-z_j\right|^{2}}\right)},\quad\forall i=1,...,n;
\end{equation}

\item If $a> 1$ then,

\begin{equation}\label{eq37}
- \beta_0\frac{z_i}{\vert z_i\vert^2}=\frac{4}{a}\displaystyle{\left(\sum_{i\neq j =1}^{m}\frac{z_i-z_j}{\left|z_i-z_j\right|^{2}}\right)},\quad\forall i=1,...,n;
\end{equation}

where $\beta_0=0$ if $0 \notin S$, while for $0\in S$ then $\beta_0$ is specified by Corollary~\ref{clocp0}.  It is understood that for $n=1$ the right hand side of \eqref{eq36} and \eqref{eq37} must be taken as equal to zero. In particular, if in $(i)$ we have, $0 \notin S$ then necessarily $n \geq 2$.
\end{enumerate}
\end{lemma}

\dimo
Under the given assumptions we know that (along a subsequence) the following holds as $ k\to\infty$:

\begin{equation}\label{eq3040}
\frac{1}{2\pi}f_k\rightharpoonup\beta_0\delta_0+\frac{4}{\max\left\{1,a\right\}}\displaystyle{\sum_{j=1}^{n}\delta_{z_j}}\quad\mbox{ weakly in the sense of measures, locally in } \mathbb{R}^2, 
\end{equation}

\begin{equation}\label{eq3041}
\nabla u_k\to\beta_0\frac{x}{\vert x \vert^2}+\frac{4}{\max\left\{1,a\right\}}\displaystyle{\sum_{j=1}^{n}\frac{x-z_j}{\vert x-z_j\vert^2}}\quad\mbox{ uniformly in } C^0_{\rm{loc}}(\mathbb{R}^2\setminus S), 
\end{equation}

with $\beta_0$ as specified above.

We establish \eqref{eq36} in the more intricate case where $n \geq 2$, as to the case $n=1$ follows by similar yet simpler arguments. For $h=1,2$ we multiply \eqref{kappa} by $\partial_{h}u_k$, and obtain:

\begin{equation}\label{gen11}
\begin{aligned}
-\partial_{h}u_k\Delta u_{k}&=e^{au_k}\partial_{h}u_k+\left|x\right|^{2N}e^{u_k}\partial_{h}u_k\\
&=\partial_{h}\left[\frac{e^{au_k}}{a}+\left|x\right|^{2N}e^{u_k}\right]-\partial_{h}\left(\log\left|x\right|^{2N}\right)\left|x\right|^{2N}e^{u_k}.
\end{aligned}
\end{equation}

For $h=1$, by straightforward computations we find:

$$
\partial_{1}u_k\Delta u_{k}=\partial_{1}\left[\frac{\left(\partial_{1}u_k\right)^{2}}{2}-\frac{\left(\partial_{2}u_k\right)^{2}}{2}\right]+\partial_{2}\left(\partial_{1}u_k\partial_{2}u_k\right).
$$

Therefore, for $\delta>0$ sufficiently small, and for every $i=1,...,n$, we may integrate \eqref{gen11} over $B_{\delta}(z_i)\subset\RR^{2}\setminus\left\{0\right\}$, to obtain:

\begin{equation}\label{gen12}
\begin{aligned}
-&\left(\int_{\partial B_{\delta}(z_i)}\frac{\left(\partial_{1}u_k\right)^{2}-\left(\partial_{2}u_k\right)^2}{2}\nu_{1}+\left(\partial_{1}u_k\partial_{2}u_k\right)\nu_2\right)\\
&=\int_{\partial B_{\delta}(z_i)}\left[\frac{e^{au_k}}{a}+\left|x\right|^{2N}e^{u_k}\right]-\int_{B_{\delta}(z_i)}\partial_{1}\left(\log\left|x\right|^{2N}\right)\left|x\right|^{2N}e^{u_k}.
\end{aligned}
\end{equation}

To handle the left hand side of \eqref{gen12}, we use \eqref{1423bis} and take into account that the function $\Psi_{i}(x):=\beta_0\frac{x}{\left|x\right|^{2}}+4\sum^{m}_{i\neq j=1}\frac{x-z_j}{\left|x-z_j\right|^{2}}$ is regular over $B_{\delta}(z_i)$. Therefore, after explicit calculations, from \eqref{gen12}, we derive:

\begin{equation}
\begin{aligned}
-&\left(\int_{\partial B_{\delta}(z_i)}\frac{\left(\partial_{1}u_k\right)^{2}-\left(\partial_{2}u_k\right)^2}{2}\nu_{1}+\left(\partial_{1}u_k\partial_{2}u_k\right)\nu_2\right)\\
&\to-2\pi\left[4\beta_0\frac{z_{i,1}}{\left|z_i\right|^{2}}+16\sum^{m}_{i\neq j=1}\frac{z_{i,1}-z_{j,1}}{\left|z_i-z_j\right|^{2}}\right]+o(1), \text{ as } k\to+\infty
\end{aligned}
\end{equation}

where we have set: $z_i=(z_{i,1},z_{i,2})$, and $o(1)\to0$ as $\delta\to0$.

Concerning the right hand side of \eqref{gen12}, we see that, $$\frac{e^{au_k}}{a}+\left|x\right|^{2N}e^{u_k}\to 0\quad\mbox{ uniformly in }\partial B_{\delta}(z_i),$$ while,

\begin{equation}\label{eq3045bis}
\displaystyle{\int_{B_{\delta}(z_i)} \partial_{1}\left(\log\left|x\right|^{2N}\right)\left|x\right|^{2N}e^{u_k}\to\left[8\pi\partial_{1}\log\left(\left|x\right|^{2N}\right)\right]_{x=z_i},\quad\mbox{as }k\to+\infty}.
\end{equation}

By using an analogous argument for $h=2$ in \eqref{gen11}, we arrive at the following identity:

\begin{equation}\label{starpol0}
-8\pi\nabla\left(\log\left|x\right|^{2N}\right)(z_i)=-2\pi\left[4\beta_0\frac{z_{i}}{\left|z_i\right|^{2}}+16\sum^{n}_{i\neq j=1}\frac{z_{i}-z_{j}}{\left|z_i-z_j\right|^{2}}\right],
\end{equation}

for every $i=1,...,n$. In other words,

\begin{equation}\label{starpol}
\left(2N-\beta_0\right)\frac{z_i}{\left|z_i\right|^{2}}-4\sum^{n}_{i\neq j=1}\frac{z_i-z_j}{\left|z_i-z_j\right|^{2}}=0\quad i=1,...,n,
\end{equation}

and \eqref{eq36} is established.

To obtain \eqref{eq37}, we argue similarly, only that in this case, by \eqref{14app}, we find:

$$ \displaystyle{\int_{B_{\delta}(z_j)} \partial_{h} \left( \log \left|x\right|^{2N} \right) \left|x\right|^{2N}e^{u_k}\to0, \quad\mbox{as }k\to+\infty},$$

from which we can derive \eqref{eq37} as above.
\qed

\begin{rem}\label{3star}
(a) As a consequence of Lemma \ref{lem3024bis} we see that, if $0 < a \ne \frac{1}{N+1} < 1$ and $S \ne \emptyset$ then 
\begin{equation}\label{eq4.34bis}
\text{either } S = \{ 0 \} \text{ or } S \text{ contains at least \underline{two} points, one of which possibly the origin}.
\end{equation}
More precisely, in complex notation, \eqref{eq36} can be expressed equivalently as follows:

\begin{equation}\label{1424bis}
\frac{(2N-\beta_0)}{z_i}-4\sum^{n}_{i\neq j=1}\frac{1}{{z_i-z_j}}=0, \quad \forall i=1,...,n.
\end{equation}

Hence, by summing up the above identities over $i=1,...,n$, we find:

\begin{equation}\label{starpoll3}
\mbox{if } 0 < a \neq \frac{1}{N+1}<1, \mbox{ then } \quad \frac{4n(n-1)}{2}=(2N-\beta_0)n.
\end{equation}

As a consequence,

\begin{itemize}

\item if $0 \notin S \ne \emptyset$ then necessarily: 

\begin{equation}\label{4.33bis}
N \in \mathbb{N} \quad \text{ and } \quad n=N+1 \geq 2
\end{equation}

\item if $0 \in S$ and $S \setminus \{ 0 \} \ne \emptyset$ then necessarily: 

\begin{equation}\label{4.33bis1}
\beta_0 = 2 (N+1-n), \quad n \in \mathbb{N}
\end{equation}

\end{itemize}

While, if $a >1$ then from \eqref{eq37} we find:

\begin{equation}\label{starpoll2}
\frac{\beta_0}{z_i}+\frac{4}{a}\sum^{m}_{i\neq j=1}\frac{1}{z_i-z_j}=0,\quad \forall i=1,...,m;
\end{equation}

which, as above, now yields to the identity:

$$\frac{2}{a}n(n-1)=-\beta_0 n.$$

with $\beta_0 = 0$ for $0 \notin S$ or $\beta_0 = \frac4a$ for $0 \in S$ Thus, if $a > 1$ and $S \ne \emptyset$ then we deduce that necessarily $S=\left\{ z_0 \right\}$ for some $z_0 \in \mathbb{R}^2$.
\end{rem}

\begin{lemma}\label{lemma3029}
Assume \eqref{cond305} and let $\frac{1}{N+1}\neq a \in \left(0,\frac{2}{N}\right) \cup \left(1,+\infty\right)$. Then $0\in S \Leftrightarrow S=\left\{ 0 \right\}$.
\end{lemma}

\dimo
If $0<a<\frac{1}{N+1}$, then by \eqref{locp03b} we know that, $\beta_0=4(N+1)$ and we can use \eqref{starpoll3} in order to derive that, $S\setminus \left\{ 0 \right\}$ must be empty.

If $a>1$ a similar conclusion follows from Remark \ref{3star}. So we are left to consider the case where:

\begin{equation}\label{444bis}
\frac{1}{N+1}<a <\frac{2}{N}
\end{equation}

To this purpose, we observe that, if $\frac{1}{N+1} < a \leq 1$, then we can still use Corollary~\ref{clocp0} together with \eqref{13app} in order to deduce the following information:

\begin{equation}\label{eq4441}
\beta_0\in\left[\frac{4}{a},4(N+1) \right] \, \mbox{ and } \, \beta_0\left(4(N+1)-\beta_0\right)=\beta\left(4(N+1)-\beta\right).
\end{equation}

On the other hand, if $\frac{1}{N+1} < a \leq \frac{2}{N+1}$ then, $\beta_0\geq 2(N+1)$, and from the identity \eqref{eq4441}, we deduce that necessarily $\beta=\beta_0$. Therefore we must have $S= \left\{ 0 \right\}$ in this case as well. Finally, for $N>1$ and $\frac{2}{N+1}<a<1$, by the identity in \eqref{eq4441} it could still happen that $\frac{4}{a}\leq\beta_0<2(N+1)$ and $\beta=4(N+1)-\beta_0$. In other words, if $S\setminus \left\{ 0 \right\}=\left\{z_1,...,z_n\right\}$ then we have $\beta_0=2(N+1-n)$, and consequently $\beta=2(N+1+n)$, for some $1\leq n \leq N+1-\frac2a < N-1$.

Thus this situation can occur only when $N>2$ and $\frac2N \leq a <1$, and the desired conclusion is established.
\qed

\begin{coro}\label{coro46}
Assume \eqref{cond305}.
\begin{enumerate}[(i)]
\item If $0 < N \leq 2$ then $0 \in S \Leftrightarrow S = \{ 0 \}$.
\item If $0 \in S$ and $\beta_0 \geq 2(N+1)$ then $S = \{ 0 \}$.
\item Let $N>2$ and $\frac2N \leq a <1$. If $0 \in S$ and $S \setminus \{ 0 \} \ne \emptyset$ then $\beta=2(N+1+n)$, $\beta_0=2(N+1-n)$, where $n \in \mathbb{N}$ is the number of blow--up points in $S \setminus \{ 0 \}$ and it satisfies: $1 \leq n \leq N+1- \frac2a < N-1$. 
\end{enumerate}
\end{coro}

Next we show that, for $0<a<1$, the situation where $S \ne \emptyset$ but $0 \notin S$ is unlikely to occur. Indeed,  it requires that $N \in \mathbb{N}$, and the blow--up behavior of $u_k$ coincides with that described in Remark \ref{rembbbbbis} for solutions  of problem \eqref{plsc2}.

\begin{thm}\label{new}
Assume \eqref{cond305} and suppose that $S \ne \emptyset$ and $S \subset \RR^{2} \setminus\left\{0\right\}$. We have:

\begin{itemize}

\item[$(i)$]
if $0<a\neq\frac{1}{N+1}<1$ then $N\in\NN$ and $S=\left\{z_1,...,z_{N+1}\right\}$  where $z_j$ are the vertices of a $(N+1)$-regular polygon, that is
\begin{equation}\label{scalingrot}
z^{N+1}_{j}=\xi_{0}, \quad j=1,...,N+1
\end{equation} 
with $\xi_{0}=(-1)^{N}z_1\cdot...\cdot z_{N+1}$.

In particular $\beta \geq 4(N+1)$, and equality holds for $\frac{1}{N+1} < a < 1$.

\item[$(ii)$] If $a>1$ then $S=\left\{z_0\right\}$, for some $z_0\in\RR^{2}\setminus\left\{0\right\}$, and $\beta\geq\frac{4}{a}$.

\end{itemize}

\end{thm}

\dimo
Clearly $(ii)$ follows directly from Remark~\ref{3star}. Hence, we assume that $0<a\neq\frac{1}{N+1}<1$. Then we can use \eqref{4.33bis} in order to conclude that, $N \in \mathbb{N}$ and $S = \{ z_1, ...,z_{N+1} \} \subset \RR^{2}\setminus\left\{0\right\}$. Therefore $\beta \geq 4 (N+1)$, and in particular, if $\frac{1}{N+1}<a<1$, then in view of \eqref{toclim2} we see that, $\beta=4(N+1)$, as claimed. 

Furthermore, since \eqref{eq36} holds with $\beta_0=0$, we may express it in complex notations as follows:

\begin{equation}\label{starpoll2nu}
2\sum^{N+1}_{i\neq j=1}\frac{z_i}{z_i-z_j}=N,\quad \forall i=1,...,N+1.
\end{equation}

We shall use \eqref{starpoll2nu} to establish \eqref{scalingrot}, whose proof we learnt from Roberto Tauraso and Carlo Pagano. We start to observe that, for $N=1$, then \eqref{scalingrot} follows easily by \eqref{starpoll2nu}. Therefore we assume that $N\geq 2$ and let

\begin{equation}\label{polnnn}
P(z)=\displaystyle{\prod^{N+1}_{j=1}\left(z-z_j\right)}=z^{N+1}+\sum^{N}_{k=0}a_kz^{k},
\end{equation} 

so that, the zeros of $P$ are given exactly by the points $z_j$, $\forall \, j=1,...,N+1$. Furthermore

$$
P'(z_i)=\displaystyle{\prod^{N+1}_{i\neq j=1}\left(z_i-z_j\right)}\quad\mbox{ and }\quad P''(z_i)=2\sum_{\left|A\right|=N-1}\displaystyle{\prod_{j\in A}\left(z_i-z_j\right)},
$$

with $A\subset\left\{1,...,N+1\right\}\setminus\left\{i\right\}$.

Hence, according to \eqref{starpoll2nu} we deduce that,

$$
NP'(z_i)=2\displaystyle{\prod^{N+1}_{i\neq j=1}\left(z_i-z_j\right)\left(\sum_{i\neq k=1}^{N+1}\frac{z_i}{z_i-z_k}\right)}=z_iP''(z_i), \quad \forall i=1,...,N+1.
$$

As a consequence, the following polynomial:

$$
Q(z):=NP'(z)-zP''(z)=\sum^{N+1}_{k=1}k\left(N+1-k\right)a_kz^{k-1}
$$

has degree less than $N$, but admits $N+1$ zeroes; as we have: $Q(z_i)=0$ for $i=1,...,N+1$. Therefore $Q\equiv 0$ and we deduce that, $a_k=0$ for $k=1,...,N$. Thus, from \eqref{polnnn} we have that $\{ z_1,...,z_{N+1} \}$ correspond to the $N+1$ distinct zeroes of the polynomial $P(z)=z^{N+1}+a_0,$ with $a_0=(-1)^{N+1}z_1 \cdot...\cdot z_{N+1} \ne 0$, and \eqref{scalingrot} is established. \\
\qed

From Lemma \ref{2app} and by the above results, we can deduce easily the following:

\begin{coro}\label{cor14bbbis}
Assume \eqref{cond305} with $a \in (0,1)$ .
\begin{enumerate}[(i)]
\item If $N \notin \mathbb{N}$ and $S\neq\emptyset$ then $0 \in S$. In particular, if $0 < N \ne 1 < 2$ then $S = \{ 0 \}$.
\item If $\frac{1}{N+1} < a < 1$ and $0 \notin S$ then $\beta = 4(N+1)$.
\item If $0 < a < \frac{1}{N+1}$ and $S \ne \emptyset$ then there exists $R_0>0$ such that, for any $R\geq R_0$ (along a subsequence) there holds:
\begin{equation}\label{eq450bis}
\lim_{k\to+\infty}\int_{\left\{\left|x\right|\leq R\right\}}\left|x\right|^{2N}e^{u_k}=\lim_{k\to+\infty}\int_{\left\{\left|x\right|\leq R\right\}}f_{k}(x)dx=4(N+1).
\end{equation}
\end{enumerate}
\end{coro}

We can further refine and complete the results of Corollary \ref{cor14bbbis} as follows:

\begin{lemma}\label{cor14bbis}
Assume \eqref{cond305}.

\begin{enumerate}[(i)]

\item If $\frac{1}{2(N+1)}\leq a<\frac{1}{N+1}$ and $S \ne \emptyset$ then $\beta=4(N+1)$. Furthermore for $R>0$ sufficiently large, we have:

\begin{equation}\label{eq3065bis}
\displaystyle{\lim_{k\to+\infty}\int_{\left\{\left|x\right|\geq R\right\}}f_k(x)=0},
\end{equation}
and in particular, if $a\neq\frac{1}{2(N+1)}$ then $0\notin\hat{S}$.

\item If $0<a<\frac{1}{2(N+1)}$ then $0\in \hat{S}$.



\end{enumerate}
\end{lemma}

\dimo


We start to establish $(i)$ in case: $\frac{1}{2(N+1)}<a<\frac{1}{N+1}$, where we need to show that $0\notin\hat{S}$. Indeed, from \eqref{toclim2}, we  know that $\beta>\frac{2}{a}$, and so if we suppose by contradiction that $0\in\hat{S}$, then $\beta=4(N+1)+\beta_{\infty}$  and by Corollary~\ref{clocpinf} we also know that: $\beta_{\infty}\geq2\beta-\frac{4}{a}$. As a consequence, we find: $4(N+1)<\beta\leq\frac{4}{a}-4(N+1)$, which is clearly impossible for $a>\frac{1}{2(N+1)}$.

Therefore $0\notin\hat{S}$, and so \eqref{eq3065bis} holds, and by \eqref{eq450bis} we find that $\beta=4(N+1)$, as claimed.

Next we show that, even when $a=\frac{1}{2(N+1)}$  we have: $\beta=4(N+1)=\frac{2}{a}$, and so \eqref{eq3065bis} holds in this case as well. Indeed, if by contradiction we suppose that $\beta>4(N+1)=\frac{2}{a}$, then, on the basis of part $(ii)$ of Lemma \ref{lemma307}, we see that necessarily $0\in\hat{S}$. Therefore we can argue exactly as above to get the following contradiction: $4(N+1)<\beta\leq\frac{4}{a}-4(N+1)=4(N+1)$.

Concerning $(ii)$, we recall that $S \cup \hat{S} \neq \emptyset$, therefore if $S = \emptyset$ then necessarily $0 \in \hat{S}$. Hence, we assume that $S\neq \emptyset$, and proceed to show that, $0 \in \hat{S}$ 
 with an argument by contradiction. Thus, we suppose that $0 \notin \hat{S}$ and recall that, from \eqref{toclim2}, we have: $\beta\geq\frac{2}{a}$. Immediately, we rule out the possibility that $\beta > \frac{2}{a}$, as otherwise we could use \eqref{eq450bis} and conclude that $\beta = 4 (N+1)$,  which is impossible for $0 < a < \frac{1}{2(N+1)}$. 
On the other hand, if $\beta = \frac{2}{a}$ then by recalling \eqref{eq401},  we can use \eqref{eq450bis} and derive:

$$
\displaystyle{\frac{1}{2\pi}\int_{\mathbb{R}^2} \vert x \vert^{2N}e^{u_k} \to 4(N+1)}, \, \mbox{ as } k\to+\infty.
$$

Hence from \eqref{lim2}, we obtain the following identity:
$$
\frac{\beta(4-a\beta)}{4\left[1-a(N+1)\right]}=4(N+1),
$$

which however is  violated  when $\beta = \frac{2}{a}$ and $0 < a < \frac{1}{2(N+1)}$. In conclusion, $0 \in \hat{S}$ and also $(ii)$ is established. \\
\qed

Concerning the case $\frac{1}{N+1}<a <1 $, on the basis of part $(ii)$ of Corollary \ref{cor14bbbis}, it remains to analyze what happens  when $0 \in S$. In this direction we have:

\begin{lemma}\label{lemma2}
Assume \eqref{cond305}. Let $N\in(0,1)$ and suppose that $S\neq \emptyset$.
\begin{enumerate}[(i)]
\item If $\frac{1}{N+1}<a<1$ then $S=\left\{ 0 \right\}$, $\beta=\beta_0$ ($\beta_0$ as characterized in Corollary~\ref{clocp0}) and $0\notin\hat{S}$.
\item If $1<a\leq \frac{2}{N+1}$ then $\beta=\frac{4}{a}$. Furthermore \eqref{eq3065bis} holds and in particular for $a\neq \frac{2}{N+1}$ we have that $0\notin\hat{S}$.
\item If $a > \frac{2}{N+1}$ then $0\in\hat{S}$.
\end{enumerate}
\end{lemma}

\dimo
The fact that in $(i)$ we have:  $S = \{ 0 \}$   follows directly from Corollary \ref{cor14bbbis}. Furthermore \eqref{eq4441} holds, and since $a<1<\frac{2}{N+1}$, it gives that $\beta_0=\beta>2(N+1)$.

Consequently, from Corollary~\ref{clocpinf} we see that $0\notin\hat{S}$, and $(i)$ is established.

To prove $(ii)$, we observe that for $S\neq\emptyset$ and $a>1$, we have $S=\{z_0\}$ and 

\begin{equation}\label{eq4.43bis}
\beta (z_0) = \beta_0 := \lim_{r\to 0}\lim_{k\to +\infty} \frac{1}{2\pi} \int_{B_r(z_0)} f_k (x) \, dx = \lim_{r\to 0}\lim_{k\to +\infty} \frac{1}{2\pi} \int_{B_r(z_0)} e^{a u_k } \, dx = \frac4a.
\end{equation}

We need to establish that actually, $\beta_0=\beta=\frac{4}{a}$. To this purpose,  we start to discuss the case: $1<a<\frac{2}{N+1}$, where we know that, $\beta\geq\beta_0=\frac{4}{a}>2(N+1)$.

Therefore, if by contradiction we assume that $\beta>\frac{4}{a}$, then $0\in \hat{S}$, and so we could use Corollary~\ref{clocpinf} to deduce that $\beta_{\infty}=\beta-\frac{4}{a}\geq2\beta-4(N+1)$. In other words, $2(N+1)<\beta\leq 4(N+1)-\frac{4}{a}$, which is clearly impossible for $a<\frac{2}{N+1}$. Thus, we have shown that,  $\beta = \beta_0 = \frac{4}{a}$ in this case.

On the other hand, for $a=\frac{2}{N+1}$ we can still claim that: $\beta= \beta_0 =\frac{4}{a}=2(N+1)$. Indeed if by contradiction we assume that, $\beta>2(N+1)$ then by part $(ii)$ of Lemma~\ref{lemma307}, we see that necessarily $0\in\hat{S}$ and at this point we could reach a contradiction as above.

Finally, to establish $(iii)$, we observe first that, if $a>\frac{2}{N+1}>1$, then by \eqref{eq4.43bis} we have that: $\beta_0=\frac{4}{a}<2(N+1) \leq \beta$. 
Thus, if by contradiction we assume that $0\notin\hat{S}$, then for $\beta > 2(N+1)$, we may check as above that it is impossible . On the other hand, if $\beta = 2(N+1)$, then we can use \eqref{eq402} and \eqref{eq4.43bis} in order to deduce  \eqref{eq4441}, which   ( for $\beta_0=\frac{4}{a}<2(N+1)$) implies that:

\begin{equation}\label{eq421}
\beta=4(N+1)-\frac{4}{a}.
\end{equation}

However, \eqref{eq421} cannot hold with $\beta=2(N+1)$ and  $a>\frac{2}{N+1}>1$, and so we conclude that $0\in\hat{S}$ as claimed.\\
\qed

\begin{lemma}
Assume \eqref{cond305}. Let $N\geq 1$ and $S\neq\emptyset$.
\begin{enumerate}[(i)]
\item If $\frac{1}{N+1}<a\leq\frac{2}{N+1}$ and $0\in S$, then $S=\left\{ 0 \right\}$ and $\beta=\beta_0$. In particular, \eqref{eq3065bis} holds and $0\notin\hat{S}$ provided that $a\neq\frac{2}{N+1}$.
\item If $N>1$ and $\frac{2}{N+1} < a <1$ then \underline{either} $\beta = \beta_0 \geq 2(N+1)$ and \eqref{eq3065bis} holds, \underline{or} $\frac{4}{a} \leq \beta_0 < 2(N+1)$, $0 \in \hat{S}$ and $\beta = 4(N+1) - \beta_0$.
In particular, 

\begin{equation}\label{eq444}
\beta  \in \left( 2(N+1), 4(N+1) - \frac{4}{a} \right]
\end{equation}

\item If $a>1$ then $0\in\hat{S}$.
\end{enumerate}
\end{lemma}

\dimo
The fact that $S=\left\{ 0 \right\}$ in $(i)$ follows from Lemma~\ref{lemma3029}, and moreover the rest of the proof of $(i)$ and $(iii)$ follows exactly as in the proof of $(i)$ and $(iii)$ of Lemma~\ref{lemma2}.

Concerning $(ii)$, we see (as in the proof of Lemma \ref{lemma3029}) that \eqref{eq4441} holds. Thus, if $\beta_0 \geq 2(N+1)$ then from \eqref{eq4441} we find that necessarily: $\beta = \beta_0$ and consequently \eqref{eq3065bis} must hold . While, if $\frac{4}{a} \leq \beta_0 < 2(N+1)$ then from \eqref{eq4441} we find that $\beta = 4(N+1) - \beta_0$. In particular, \eqref{eq444} holds and necessarily $0 \in \hat{S}$ in this case.\\
\qed

In addition, we may also conclude the following:

\begin{coro}\label{coro415} 
If  $a > \max \left\{ 1, \frac{2}{N+1} \right\}$ then $0 \in \hat{S}$. Furthermore if $S \ne \emptyset$ then $$2(N+1) \leq \beta \leq 4(N+1) - \frac{4}{a}$$.
\end{coro}

\dimo
The first part of the statement is a direct consequence of the above results. Concerning the second part, recall that, $\beta \geq 2(N+1)$, therefore we only need to check that if $\beta > 2(N+1)$ then $\beta \leq 4(N+1) - \frac{4}{a}$. To this purpose, we simply observe that, since $S \ne \emptyset$ and $0 \in \hat{S}$ then necessarily: $\beta = \frac{4}{a} + \beta_\infty$ with $\beta_\infty \geq 2 \beta - 4 (N+1)$ (see \eqref{eq359bis}), and the desired estimate readily follows.
\qed

\section{The case $0<a<\frac{1}{N+1}$.}
\setcounter{equation}{0}

By virtue of part $(i)$ in Lemma~\ref{cor14bbis}, it remains to identify the limiting value of $\beta=\beta^\infty$ in \eqref{131bis} only in the following situations:

\begin{equation}\label{eq501}
\frac{1}{2(N+1)}\leq a < \frac{1}{N+1} \mbox{ and } S=\emptyset,
\end{equation}

or

\begin{equation}\label{eq502}
0<a<\frac{1}{2(N+1)}, \text{ where we know that,  } 0 \in \hat{S};
\end{equation}

see $(ii)$ of Lemma~\ref{cor14bbis}.

To this purpose, we fix $0 < \varepsilon < 1$ sufficiently small (to be specified later) and let $R_k=R_k(\varepsilon)>0$ be the \underline{unique} value defined by the condition:

\begin{equation}\label{eq50}
\frac{1}{2\pi}\displaystyle{ \int_{\vert x\vert\leq R_k} f_k(x)dx=\max\left\{\frac{2}{a},4(N+1) \right\}-\varepsilon}.
\end{equation}

\begin{lemma}\label{lem53}
Assume that either \eqref{eq501} or \eqref{eq502} holds. For $\varepsilon>0$ sufficiently small we have:

\begin{equation}\label{eq511}
R_k \to +\infty, \mbox{ as } k\to+\infty.
\end{equation}

\end{lemma}

\dimo
If $S=\emptyset$ then \eqref{eq511} is a direct consequence of Lemma~\ref{lemma307}-$(i)$. Hence assume that,

\begin{equation}\label{eq501bis}
0< a < \frac{1}{2(N+1)} \mbox{ and } S\neq\emptyset,
\end{equation}

and let $0<\varepsilon < \frac{2}{a}-4(N+1)$.  Then from \eqref{eq450bis} we find that, for any $R>>1$ sufficiently large, we have:

\begin{equation}\label{eq51}
\frac{1}{2\pi}\displaystyle{ \int_{\vert x\vert\leq R} f_k(x)dx\to 4(N+1)<\frac{2}{a}-\varepsilon=\max\left\{\frac{2}{a},4(N+1) \right\}-\varepsilon}, \mbox{ as } k\to+\infty.
\end{equation}

Consequently, \eqref{eq50} can hold only if $R_k\to+\infty$, as $k\to+\infty$, as claimed.\\
\qed

Define:

\begin{equation}\label{eq52}
v_k(x)=u_k(R_kx)+\frac{2}{a}\log R_k,
\end{equation}

which satisfies 

\begin{equation}\label{eq53}
\begin{cases}
-\Delta v_k=e^{av_k}+\varepsilon_{1,k} \left|x\right|^{2N}e^{v_k}:=f_{1,k}(x) \mbox{ in } \mathbb{R}^2\\
\beta_k:=\frac{1}{2\pi}\displaystyle{\int_{\RR^{2}}f_{1,k}(x)dx} \to \beta \, \text{ as } k \to \infty.
\end{cases}
\end{equation}

with

\begin{equation}\label{eq53bis}
\varepsilon_{1,k} = R_{k}^{-\frac{2}{a}(1-a(N+1))} \to 0 \, \text{ as } k \to \infty
\end{equation}

Moreover, from \eqref{eq50} we have:

\begin{equation}\label{eq54}
\frac{1}{2\pi}\displaystyle{ \int_{\vert x\vert\leq 1} f_{1,k}(x)dx=\max\left\{\frac{2}{a},4(N+1) \right\}-\varepsilon}.
\end{equation}

Next set,

\begin{equation}\label{eq55}
\hat{v}_k(x)=v_k \left(\frac{x}{\vert x \vert^2}\right)+\beta_k\log \frac{1}{\vert x \vert},
\end{equation}

satisfying:

\begin{equation}\label{eq56}
\begin{cases}
-\Delta \hat{v}_k=\vert x \vert^{a\beta_k-4}e^{a\hat{v}_k}+\varepsilon_{1,k}\left|x\right|^{\beta_k-2(N+2)}e^{\hat{v}_k}:=\hat{f}_{1,k}(x) \mbox{ in } \mathbb{R}^2\\
\beta_k:=\frac{1}{2\pi}\displaystyle{\int_{\RR^{2}}\hat{f}_{1,k}(x)dx} \to \beta, \text{ as } k \to +\infty,
\end{cases}
\end{equation}

and we denote by $S_1$ and $\hat{S}_1$ the (possibly empty) blow--up set of $v_k$ and $\hat{v}_k$ respectively. \\
For $x_0 \in S_1$ we let, 

\begin{equation}\label{eq56bis}
\beta_0 (x_0) = \lim_{r\to0}\lim_{k\to+\infty}\frac{1}{2\pi}\int_{B_{r} (x_0)} f_{1,k} (x) \, dx = \beta_{0,1} (x_0) +\beta_{0,2} (x_0)
\end{equation}

with

\begin{equation}\label{eq56tris}
\begin{split}
\beta_{0,1} (x_0) :=\lim_{r\to0}\lim_{k\to+\infty}\frac{1}{2\pi}\int_{B_{r} (x_0)} e^{a v_k(x)} \, dx \\
\; \beta_{0,2} (x_0):=\lim_{r\to0}\lim_{k\to+\infty}\frac{1}{2\pi}\int_{B_{r} (x_0)} \varepsilon_{1,k} |x|^{2N} e^{a v_k(x)} \, dx.
\end{split}
\end{equation} 

For later use, we collect in the following lemma some general facts concerning solutions of \eqref{eq53} with $\varepsilon_{1,k} \to 0$.

By direct inspection, we check that the arguments in the proof of Proposition~\ref{prop17} and Proposition~\ref{1310} extend easily to cover the new sequence $\left\{ v_k \right\}$ in \eqref{eq53}.

Furthermore, Proposition \ref{prop33} as well as Corollary \ref{coro31}, Corollary \ref{coro34bis} and Proposition \ref{prop35} obviously apply to $v_k$ (and $\hat{v}_k$) and imply the following,

\begin{lemma}\label{lem50}
Let $v_k$ satisfy \eqref{eq53} with $N > -1$, $0<a \ne \frac{1}{N+1}$ and $a \ne 1$ and $\varepsilon_{1,k} \to 0$.
Then the blow--up set $S_1$ of $v_k$ contains at most finitely many points and the following holds:
\begin{enumerate}[(i)]

\item for every compact set $K \subset \mathbb{R}^2 \setminus S_1$ there exists a constant $C=C(K)>0$:
\begin{equation}\label{eq5.14}
|v_k (x) - v_k (y)| \leq C \quad \forall \, x,y \in K
\end{equation}

\item Let $x_0 \in S_1$ then

\begin{itemize}

\item for $x_0 = 0$ and $\beta_0 = \beta_0 (0)$, $\beta_{0,j} = \beta_{0,j} (0)$, $j=1,2$, we have:

\begin{equation}\label{eq5.14.1}
\beta_{0,1} = \frac{a \beta_0 (\beta_0 - 4(N+1))}{4(1-a(N+1))} \quad \text{ and } \quad \beta_{0,2} = \frac{\beta_0 (4 - a \beta_0)}{4(1-a(N+1))}
\end{equation}

\item for $x_0 \ne 0$ and $\beta_0 = \beta_0 (x_0)$, we have:

\begin{equation}\label{eq5.14.2}
\beta_{0,1} (x_0) = \frac{a \beta_0 (\beta_0 - 4)}{4(1-a)} \quad \text{ and } \quad \beta_{0,2} (x_0) = \frac{\beta_0 (4 - a \beta_0)}{4(1-a)}
\end{equation}

\item If one of the following conditions hold: 
\begin{eqnarray*}
&(a)& x_0 = 0 \text{ and } a \geq \frac{1}{2(N+1)} \\
&(b)& x_0 \ne 0 \text{ and } a \geq \frac12 \\
&(c)& \beta (x_0) \geq \frac2a
\end{eqnarray*}
then
\begin{equation}\label{eq5.14.3}
\frac{1}{2\pi}f_{1,k} \rightharpoonup \beta_0  \delta_{x_0} \quad\mbox{ weakly in the sense of measures in } B_{\delta} (x_0).
\end{equation}
Moreover if $a>1$ then $\beta (x_0) = \frac{4}{a}$, $\, \forall \, x_0 \in S_1$.

\item In case $\beta_0 (x_0) < \frac{2}{a}$, then we have:
\begin{equation}\label{eq5.18}
\begin{split}
&\text{ if } x_0 = 0 \text{ then } \beta (x_0) = 4(N+1) \text{ and } 0 < a < \frac{1}{2(N+1)}\\ 
&\text{ if } x_0 \ne 0 \text{ then } \beta (x_0) = 4 \text{ and } 0 < a < \frac{1}{2}
\end{split}
\end{equation}
In particular,
\begin{equation}\label{eq5.14.5}
\beta_{0,1} (x_0) = 0
\end{equation}
\end{itemize}

\item If $0 \in \hat{S}_1$ and $\beta > \max \{ \frac{2}{a}, 2(N+1)\}$ then by setting: 
\begin{equation} \label{eq5.14.5bis}
\beta_{\infty} := \lim_{r \to 0} \lim_{k \to +\infty} \frac{1}{2\pi} \int_{|x|<r} \hat{f}_{1,k} (x) \, dx
\end{equation}
we have:
\begin{equation} \label{eq5.14.6}
\lim_{r\to0}\lim_{k\to+\infty}\frac{1}{2\pi}\int_{|x|<r} |x|^{a \beta_k -4} e^{a \hat{v}_k} = \frac{a \beta_\infty (2 \beta - 4(N+1) - \beta_\infty)}{4(1-a(N+1))}
\end{equation}
\begin{equation} \label{eq5.14.7}
\lim_{r\to0}\lim_{k\to+\infty}\frac{1}{2\pi}\int_{|x|<r} \varepsilon_{1,k} |x|^{\beta_k - 2(N+2)} e^{a \hat{v}_k} = \frac{\beta_\infty (4+ a \beta_\infty - 2 a \beta)}{4(1-a(N+1))}
\end{equation}
\end{enumerate}
\end{lemma}
\qed

\vspace{0.3cm}

We are going to use Lemma \ref{lem50} for the sequence $v_k$ in \eqref{eq53}.\\

\begin{lemma}\label{lem56}
Assume \eqref{eq501} or \eqref{eq502}. If $\beta = \frac{2}{a}$ then

\begin{equation}\label{eq525}
N \in \mathbb{N}, \quad a = \frac{1}{2(N+1)} \; \text{ and } \; \lim_{r \to 0} \lim_{k\to+\infty} \frac{1}{2\pi}\int_{|x| \geq R} f_{1,k} (x) = 0
\end{equation}

Furthermore the blow--up set $S_1$ of $v_k$ is formed (up to a rotation) by the $(N+1)$-roots of the unity.
\end{lemma}

\dimo
Since $\beta=\frac{2}{a}$, then necessarily: $0 < a \leq \frac{1}{2(N+1)}$ and therefore we have: $\frac{1}{2\pi}\int_{|x| <1} f_{1,k} (x) \, dx = \frac2a - \varepsilon$. In other words, $\frac{1}{2\pi}\int_{|x| >1} \hat{f}_{1,k} = \frac{1}{2\pi}\int_{|x|<1} f_{1,k} = \varepsilon$, with $\varepsilon > 0$ sufficiently small. By virtue of these facts, we can argue exactly as in Lemma \ref{lem42}. Indeed we can use inequality \eqref{preq404} for $\hat{v}_k$ together with \eqref{eq5.14} in order to deduce that,  for $R_0>1$ sufficiently large, we have:

\begin{equation}\label{eq513}
\displaystyle{\max_{R_0\leq\vert x \vert \leq R} v_k\to-\infty}, \mbox{ as } k\to+\infty, \forall R>R_0.
\end{equation}

Similarly, we can use the estimate \eqref{4.6.1} for $\hat{v}_k$ to obtain:

\begin{equation} \label{eq5.25bis}
\lim_{R \to +\infty} \lim_{k\to+\infty} \frac{1}{2\pi}\int_{|x| \geq R} \varepsilon_{1,k} |x|^{2 N} e^{v_k} = 0
\end{equation}

From those properties, we see that necessarily $S_1\neq\emptyset$. Indeed, if $S_1=\emptyset$, then we could use \eqref{eq5.14} and \eqref{eq513} to get,

\begin{equation*}
\displaystyle{\max_{\vert x \vert \leq R} v_k\to-\infty}, \mbox{ as } k\to+\infty,;
\end{equation*}

$\forall \, R>0$, in contradiction with \eqref{eq54}.

Therefore we let $S_1=\left\{z_1,...,z_m\right\}\subset\mathbb{R}^2$, $m\in\mathbb{N}$, and set,

\begin{equation}\label{eq515}
\beta_j:=\beta_0 (z_j):=\displaystyle{\lim_{r\to 0} \lim_{k \to + \infty} \frac{1}{2\pi} \int_{B_r (z_j)} f_{1,k}(x) dx}, \; j=1,...,m.
\end{equation}

Notice that by Lemma~\ref{lem50}, we know that: $4 \leq \beta_j \leq \beta = \frac{2}{a}$, for every $j=1,...,m$.

Furthermore, from \eqref{eq5.14} and \eqref{eq513}, we see that necessarily "concentration" must occur, in the sense that, as $k\to\infty$,
\begin{equation}\label{eq611}
\frac{1}{2\pi}f_{1,k}\rightharpoonup\displaystyle{\sum_{j=1}^m}\beta_j\delta_{z_j}  \quad\mbox{ weakly in the sense of measures, locally in } \mathbb{R}^2 ,
\end{equation}
Clearly, if $z_j \in S_1$ and $|z_j|<1$ then $\beta_j < \frac{2}{a}$; actually, we see next that such estimate for $\beta_j$ holds for every $j=1,...,m$.

{\bf{CLAIM:}}
\begin{equation}\label{eq61}
\beta_j<\frac{2}{a}, \quad \forall j=1,...,m.
\end{equation}

To establish \eqref{eq61} we argue by contradiction and suppose that there exists $j_0\in\left\{ 1,...,m \right\}$ such that: $\beta_{j_0} \geq \frac{2}{a} = \beta$. Hence in this case we must have that $m=j_0=1$ and,

\begin{equation}\label{eq62}
S_1=\left\{z_1\right\}, \quad \beta_1=\frac{2}{a}=\beta.
\end{equation}

In particular, for every $r>0$:

\begin{equation}\label{eq621}
\lim_{k \to + \infty} \int_{\mathbb{R}^2\setminus B_r(z_1)} f_{1,k}(x) dx=0
\end{equation}

By keeping in mind \eqref{eq54} and \eqref{eq611}, we see that actually: $\vert z_1 \vert =1$.

On the other hand, if we use \eqref{eq5.14.2} for the blow--up point $x_0 = z_1\neq 0$, together with \eqref{33bis} we obtain the following identity:

\begin{equation}\label{eq65}
\frac{a\beta_1(\beta_1-4)}{4(1-a)}=\displaystyle{\lim_{r \to 0} \lim_{k\to+\infty} \frac{1}{2\pi} \int_{B_r(z_1)}e^{av_k}dx}=\displaystyle{\lim_{k\to+\infty} \frac{1}{2\pi} \int_{\mathbb{R}^2}e^{av_k}dx}=$$
$$\displaystyle{\lim_{k\to+\infty} \frac{1}{2\pi} \int_{\mathbb{R}^2}e^{au_k}dx}=\frac{a\beta(\beta-4(N+1))}{4(1-a(N+1))}
\end{equation}

which must hold with $\beta_1=\beta=\frac{2}{a}$, and this is clearly impossible. So \eqref{eq61} is established. Therefore, we can use Lemma \ref{lem50} and from \eqref{eq5.18} and \eqref{eq5.14.5} we deduce:

\begin{equation}\label{eq535}
\beta_j = \begin{cases} 4(N+1) &\mbox{if } z_j = 0 \\ 
4 & \mbox{if } z_j\neq0. \end{cases} \quad \mbox{ and } \quad  \displaystyle{\lim_{r\to 0} \lim_{k \to + \infty} \frac{1}{2\pi} \int_{B_r(z_j)} e^{av_k} =0, \quad \forall \, j=1,...,m}.
\end{equation}

On the basis of the "concentration" property \eqref{eq611}, we can check easily that the conclusion of Lemma~\ref{lem3024bis}, Lemma~\ref{lemma3029} and part $(i)$ of Theorem \ref{new} apply to $v_k$. As a consequence, we conclude that, either $S_1=\left\{ 0 \right\}$ or $N\in\mathbb{N}$, $m=N+1$ and $S_1$ is formed by the vertices of a regular $N+1$-polygon. In any case, by combining \eqref{eq5.25bis}, \eqref{eq611} and \eqref{eq535} we find:

\begin{equation}\label{eq68}
\begin{split} 
\frac{\beta(4- a \beta)}{4(1-a(N+1))} = \lim_{k\to+\infty} \frac{1}{2\pi} \int_{\mathbb{R}^2} \varepsilon_{1,k} |x|^{2N} e^{v_k} \, dx=4(N+1),\\ 
\lim_{R \to +\infty} \lim_{k \to +\infty} \frac{1}{2\pi} \int_{\vert x \vert \geq R} e^{a v_k} \, dx = \frac{a \beta(\beta-4(N+1))}{4(1-a(N+1))}
\end{split}
\end{equation}

But for $\beta = \frac{2}{a}$ we readily check that \eqref{eq68} can hold if and only if:

\begin{equation}\label{eq68bis}
a=\frac{1}{2(N+1)}, \quad \beta=\frac{2}{a}=4(N+1), \mbox{ and } \lim_{R \to +\infty} \displaystyle{\lim_{k\to+\infty} \int_{\vert x \vert \geq R}f_{1,k}dx=0} 
\end{equation}

as claimed. Furthermore, by means of  \eqref{eq54} and \eqref{eq68bis} we can rule out the possibility that $S_1 = \{ 0 \}$.

In conclusion, we must have that $N \in \mathbb{N}$ and (in complex notations) $S_1 = \{ z_1,...,z_{N+1} \} \subset \mathbb{R}^2 \setminus \{ 0 \}$ is formed by the $N+1$-roots of the value: $\xi_0 = (-1)^N z_1 \cdot ... \cdot z_{N+1} \ne 0$. In particular, $|z_j|=|\xi_0|^{\frac{1}{N+1}} \; \forall j = 1,..,N+1$. But from \eqref{eq54} we see that necessarily $|z_j|=1$, $\forall j = 1,..,N+1$, and the proof is completed.\\
\qed

On the basis of Lemma \ref{lem56}, from now on we shall assume the following for $\varepsilon>0$:

\begin{equation}\label{eq5.36bis1}
\begin{split}
&\text{ if } a \ne \frac{1}{2(N+1)}  \text{ then } 0 < \varepsilon < \min \left\{ \beta - \frac2a, \left| 4(N+1) - \frac{2}{a} \right|\right\}\\
&\text{ if } \beta < \frac{4}{a} \text{ then } 0 < \varepsilon < \frac4a - \beta\\
&\text { if } a = \frac{1}{2(N+1)} \text{ and }  \beta > \frac{2}{a}  \text{ then } 0< \varepsilon< \beta - \frac{2}{a}. 
\end{split}
\end{equation}

\begin{prop}\label{prop55}
Assume \eqref{eq501} or \eqref{eq502}, then 

\begin{equation}\label{eq5.36bis}
\lim_{R \to +\infty} \lim_{k \to +\infty} \int_{\vert x \vert \geq R} f_{1,k} (x) \, dx = 0
\end{equation}

More precisely, for $a \ne \frac{1}{2(N+1)}$ or for $a = \frac{1}{2(N+1)}$ and $\beta > \frac{2}{a}$, we have that $0 \notin \hat{S}_1$.
\end{prop}

\dimo
By virtue of Lemma \ref{lem56}, we are left to show that,

\begin{equation}\label{eq5.37bis}
\text{ if } \beta > \frac{2}{a} \text{ then } 0 \notin \hat{S}_1
\end{equation}

In order to establish \eqref{eq5.37bis}, we argue by contradiction and suppose that $\beta > \frac{2}{a}$ and $0 \in \hat{S}_1$. Hence we are in position of applying part $(iii)$ of Lemma \ref{lem50}, and from  \eqref{eq5.14.7} we conclude that, 

\begin{equation}\label{eq5.37.1}
\beta_\infty \geq 2 \beta - \frac{4}{a}.
\end{equation}

with $\beta_\infty$ defined in \eqref{eq5.14.5bis}.

On the other hand from \eqref{eq54} we know also that, 

\begin{equation}\label{eq5.37.2}
\beta_\infty \leq \beta - \max \{\frac{2}{a}, 4(N+1) \} + \varepsilon,
\end{equation}

and we may combine \eqref{eq5.37.1} and \eqref{eq5.37.2} to deduce that,

\begin{equation}\label{eq5.37.3}
\beta \leq \frac{4}{a} - \max \{\frac{2}{a}, 4(N+1) \} + \varepsilon
\end{equation}

But since $\beta > \frac{2}{a}$, we readily see that the inequality \eqref{eq5.37.3} is in contradiction with our choice of $\varepsilon > 0$ in \eqref{eq5.36bis1}.\\
\qed

\begin{thm}\label{teo5.6}
Let $v_k$ be defined in \eqref{eq62} and $S_1$ be its (possibly empty) blow--up set.

\begin{enumerate}[(i)]
\item If $S_1 = \emptyset$ or $S_1 = \{ z_1 \} \subset \mathbb{R}^2 \setminus \{ 0 \}$ then $\beta = \frac{4}{a}$.
\item If $0 \in S_1$ then $S_1 = \{ 0 \}$, $0 < a < \frac{1}{2(N+1)}$ and $\beta = \frac{4}{a} - 4(N+1)$.
\item If $S_1 = \{ z_1,...,z_m \} \subset \mathbb{R}^2 \setminus \{ 0 \}$ with $\mathbb{N} \ni m \geq 2$, then $N \geq 1$ and for $\beta_j$ in \eqref{eq515} the following holds:

\begin{itemize}
\item \underline{either} $\beta_j < \frac{2}{a} \; \forall j=1,...,m$, and we have:

\begin{equation}\label{eq516}
2 \leq m \leq N+1 \; \text{ and } \; \beta = \frac{2}{a} \left( 1 + \sqrt{1-4am(1-a(N+1))} \right),
\end{equation}

moreover the equality $m = N+1$, can occur only for $\frac{1}{2(N+1)} \leq a < \frac{1}{N+1}$ and $N \in \mathbb{N}$, and in this case we have that,  $\beta = 4(N+1)$ and $S_1$ is formed (up to rotation) by the $(N+1)-$ roots of the unity ;

\item \underline{or} $\, \exists \,$  \underline{unique} $j_0 \in \{ 1,...,m \}$ such that $|z_{j_0}|=1$, $\beta_{j_0} \geq \frac{2}{a}$ and we have:

\begin{equation}\label{eq517}
1 \leq m-1 \leq  \frac{1}{2a} \left( \sqrt{(1-a(N+1))^2 + \frac{Na}{1-a}}-1 \right) + \frac{N+1}{2} < \frac{N}{2(1-a)} < \frac{N+1}{2}
\end{equation}

\begin{equation}\label{eq518}
\begin{split}
\beta = &\frac{2}{a} \sqrt{\left( 1 - \frac{2(m-1)(1-a(N+1))}{N} \right)^2 + \frac{4(m-1)ma}{N} (1-a(N+1))}\\
& + \frac{2}{a} \left( 1 - \frac{2(m-1)(1-a(N+1))}{N} \right) 
\end{split}
\end{equation}

In particular in this case there holds,

\begin{equation}\label{eq546bis}
N>1 \text{ and } \beta \geq 2(N+1) + \frac{2}{a} \left( \sqrt{(1-a(N+1))^2 + \frac{Na}{1-a}}\right)
\end{equation}

\end{itemize}

\end{enumerate}
\end{thm}

\dimo If $S_1 = \emptyset$ then, by passing to a subsequence if necessary, we can assume that,
$$
v_k\to v \mbox{ in } C^{2}_{loc}(\mathbb{R}^2)
$$
with $v$ satisfying:

$$\begin{cases}
-\Delta v=e^{av} \text{ in } \mathbb{R}^2\\
e^{a v} \in L^1 (\mathbb{R}^2)
\end{cases}$$

So, from \eqref{inpa} and \eqref{toclim2} we conclude that necessarily $\beta = \frac{1}{2\pi} \int_{\mathbb{R}^2} e^{a v} = \frac4a$, as claimed.

Next we suppose that $S_1 \ne \emptyset$.

{\it \underline{Claim 1}: if $0 \in S_1$ then necessarily $0 < a < \frac{1}{2(N+1)}$ and $\beta_0 := \lim_{r \to 0} \lim_{k\to+\infty} \int_{B_r} f_{1,k} (x) \, dx$ satisfies,}

\begin{equation}\label{eq519}
4(N+1) = \beta_0 \leq \frac{2}{a} - \varepsilon
\end{equation}

Indeed, for $\frac{1}{2(N+1)} \leq a < \frac{1}{N+1}$ the condition \eqref{eq54} reads as follows: $\frac{1}{2\pi} \int_{|x| \leq 1} f_{1,k} (x) \, dx = 4(N+1) - \varepsilon$, $\varepsilon > 0$, and by virtue of \eqref{eq5.14.1} we see that $0 \notin S_1$ in this case. At this point \eqref{eq519} follows from \eqref{eq54} and \eqref{eq5.18}.

{\it \underline{Claim 2}: Let $S_1 = \{ z_1,...,z_m \}$ and suppose that $\exists \, j_0 \in \{ 1,...,m \}$: $\beta_{j_0} \geq \frac{2}{a}$, then $|z_{j_0}|=1$, $0 \notin S_1$, $\beta_j < \frac{2}{a} \; \forall j \ne j_0$ (if any) and,

\begin{equation}\label{eq520}
\begin{split}
\beta = &\frac{2}{a} \sqrt{\left( 1 - \frac{2(m-1)(1-a(N+1))}{N} \right)^2 + \frac{4(m-1)ma}{N} (1-a(N+1))}\\
& + \frac{2}{a} \left( 1 - \frac{2(m-1)(1-a(N+1))}{N} \right) 
\end{split}
\end{equation}

Furthermore, either $m=1$ and $\beta = \frac{4}{a}$ or \eqref{eq517} and \eqref{eq546bis} hold. }\\

We observe first, that by virtue of Lemma \ref{lem50} and Proposition \ref{prop55}, (along a subsequence) we have:

\begin{equation}\label{eq522}
\frac{1}{2\pi} f_{1,k} \rightharpoonup \sum_{j=1}^{m} \beta_j \delta_{z_j} \mbox{ weakly in the sense of measures, and } \beta = \sum_{j=1}^{m} \beta_j
\end{equation}

Moreover from Claim 1 and \eqref{eq519}  we have that necessarily $z_{j_0} \ne 0$. 

We deal first with the case $m=1$, where $\beta = \beta_{j_0}$ and, in view of \eqref{eq54}, we see that necessarily $|z_{j_0}|=1$. By virtue of \eqref{eq5.14.2}, we deduce that,

\begin{equation}\label{eq523}
\frac{a \beta_{j_0} (\beta_{j_0} - 4)}{4(1-a)}=\frac{a\beta(\beta-4(N+1))}{4(1-a(N+1))}.
\end{equation}

Thus, if $\beta = \beta_{j_0}$ then from \eqref{eq523} we derive that $\beta = \frac{4}{a}$, as claimed.

So let $m \geq 2$, and by contradiction assume that there exists $j_1 \ne j_0 \in \{ 1,...,m \}$: $\beta_{j_1} \geq \frac{2}{a}$ and $\beta_{j_0} \geq \frac{2}{a}$.

Since $\beta \leq \frac{4}{a}$, then necessarily: $\beta_{j_1} = \beta_{j_0} = \frac{2}{a}$, $\beta = \frac{4}{a}$ and $m=2$. On the other hand as above we see that, $z_{j_1} \ne 0$, $z_{j_0} \ne 0$ and 

$$\frac{2}{a} \left( \frac{1-2a}{1-a} \right) = \frac{a (\beta_{j_0} (\beta_{j_0} - 4) + \beta_{j_1} (\beta_{j_1} - 4))}{4(1-a)}=\frac{a\beta(\beta-4(N+1))}{4(1-a(N+1))} = \frac{4}{a},$$

which is clearly impossible.

Therefore, $\beta_0 (z_j) := \beta_j < \frac{2}{a}$, $\; \forall \, j \in \{ 1,...,m \} \setminus \{ j_0 \}$ and we can use Lemma \ref{lem50} in order to deduce that, for every $j \ne j_0$ the following holds:

\begin{equation}\label{eq524}
\begin{split}
&\text{ if } z_j = 0 \text{ then } \; 0 < a < \frac{1}{2(N+1)} \; \text{ and } \; \beta_j = 4(N+1),\\ 
&\text{ if } z_j \ne 0 \text{ then } \; 0 < a < \frac{1}{2} \; \text{ and } \; \beta_j = 4.
\end{split}
\end{equation}

Furthermore  \eqref{eq523} continue to hold, with  $\beta = \sum_{j=1}^{m} \beta_j$ and where the values of $\beta_j$, for $j \ne j_0$, are specified in \eqref{eq524} .

In particular \eqref{eq523} implies that, for $m \geq 2$ necessarily $\beta < \frac{4}{a}$ and so, by our choice of $\varepsilon > 0$ in \eqref{eq5.36bis1} and \eqref{eq54}, we find that necessarily $|z_{j_0}|=1$. Furthermore, we observe that $0 \notin S_1$, as otherwise we would have:

$$\beta \geq \beta_{j_0} + 4 (N+1) \quad \text{with }\beta_{j_0} \geq \frac{2}{a},$$

which is impossible by virtue of \eqref{eq523}.

Consequently, 

\begin{equation}\label{eq526}
\beta = \beta_{j_0} + 4(m-1),
\end{equation}

and we can use again \eqref{eq523} together with \eqref{eq526} in order to solve for $\beta_{j_0}$ and obtain:

\begin{equation}\label{eq527}
\begin{split}
\beta_{j_0} = &\frac{2}{a} \sqrt{\left( 1 - \frac{2(m-1)(1-a)}{N} \right)^2 + \frac{4(1-a)a(m-1)}{N} (N+2-m)}\\
&+ \frac{2}{a} \left( 1 - \frac{2(m-1)(1-a)}{N}  \right).
\end{split}
\end{equation}

\normalsize

Therefore,

\begin{equation}\label{eq527bis}
\begin{split}
\beta = &\frac{2}{a} \sqrt{\left( 1 - \frac{2(m-1)(1-a(N+1))}{N} \right)^2 + \frac{4(m-1)ma}{N} (1-a(N+1))}\\
& + \frac{2}{a} \left( 1 - \frac{2(m-1)(1-a(N+1))}{N} \right) 
\end{split}
\end{equation}

as claimed.

At this point, if in \eqref{eq527} we require that $\beta_{j_0} \geq \frac{2}{a}$, then after straightforward calculations (explicitly carried out in part $(i)$ of the Calculus Lemma 1 in the Appendix), we find:

\begin{equation}\label{eq528}
\begin{split}
m \in \mathbb{N}: 1 \leq m-1 \leq  &\frac{1}{2a} \left( \sqrt{(1-a(N+1))^2 + \frac{Na}{1-a}}- (1 - a(N+1)) \right)=\\
&\frac{N}{2(1-a)} \left( \frac{1}{1 - a(N+1) +  \sqrt{(1-a(N+1))^2 + \frac{Na}{1-a}}} \right).
\end{split}
\end{equation}

Since we can easily check that, for $a \in \left( 0, \frac{1}{N+1} \right)$, we have:

\begin{equation}\label{eq529}
\sqrt{(1-a(N+1))^2 + \frac{Na}{1-a}} < 1,
\end{equation}

then, by means of \eqref{eq528}, we easily deduce that, $1 \leq m-1 < \frac{N+1}{2}$. As a consequence, $N > 1$ and  $0 < a < \frac{1}{N+1} < \frac{1}{2}$, consistently with the second condition in \eqref{eq524}.

On the other hand, for $N>1$  we can also check that,  

$$\sqrt{(1-a(N+1))^2 + \frac{Na}{1-a}} > a(N+1)$$

and from \eqref{eq528}, we readily conclude \eqref{eq517}. Finally, the estimate \eqref{eq546bis} follows by a monotonicity property which is explicitly derived in part $(ii)$ of the Calculus Lemma 1 in the Appendix.

Next we analyze what happens in case $\beta_j < \frac{2}{a}$, $\, \forall j=1,...,m$ and prove:\\

{\it \underline{Claim 3}: If $0 \notin S_1 = \{ z_1,...,z_m \}$ and assume that $\beta_j < \frac{2}{a}$, $\, \forall j=1,...,m$, then $m \geq 2$ and the following holds:

\begin{enumerate}[(a)]
\item \underline{either}, 

\begin{equation}\label{eq516}
\beta = \frac{2}{a} \left( 1 + \sqrt{1 - 4 am (1 - a(N+1))} \right),
\end{equation}

\begin{equation}\label{eq516bis}
2 \leq m < N+1
\end{equation}

and in particular $N > 1$ in this case;

\item \underline{or}, 

$$\frac{1}{2(N+1)} \leq a < \frac{1}{N+1}, \; N \in \mathbb{N} \;, \, m = N+1 \, \text{ and } \beta = 4 (N+1);$$ 

and up to a rotation, $S_1$ is formed by the $N+1$-roots of the unity.
\end{enumerate}}

To establish Claim 3 we notice first that, according to Lemma \ref{lem50}, under the given assumption we have: $\beta_j  = 4$, $\, \forall j=1,...,m$, and $0 < a < \frac12$. To proceed further, we need to distinguish between the case where ''concentration''  occurs or not. In case of "concentration", namely when (along a subsequence) the following holds:

$$\frac{1}{2\pi} f_{1,k} \rightharpoonup 4 \sum_{j=1}^{m} \delta_{z_j} \mbox{ weakly in the sense of measures}, $$

then  $ \beta = 4 m$ and we can argue exactly as in Lemma \ref{lem3024bis} and Remark \ref{3star} in order to deduce that necessarily: $N \in \mathbb{N}$, $m = N+1$ and $\beta = 4 (N+1)$. So necessarily in this case we must have that, $\frac{1}{2(N+1)} \leq a < \frac{1}{N+1}$, and exactly as in Lemma \ref{lem56}, we may also conclude that (up to a rotation) $S_1$ is formed by the $N+1$-roots of the unity, as claimed in $(b)$.

In case ''concentration'' does not occur, then in view of \eqref{eq54}, we find that $v_k$ is uniformly bounded, say in $C^1$-norm, over compact subsets of $\mathbb{R}^2\setminus S_1$. 

Consequently, along a subsequence, we find that: 

$$
v_k\to v \mbox{ in } C^{2}_{loc}(\mathbb{R}^2\setminus S_1)
$$

and in account of \eqref{eq5.36bis} and \eqref{eq5.14.5} we have, as $k\to+\infty$,

$$
\int_{\mathbb{R}^2}e^{av_k} \to \int_{\mathbb{R}^2}e^{av}; \quad \frac{1}{2\pi}\int_{\mathbb{R}^2} \varepsilon_{1,k} \vert x \vert^{2N} e^{v_k}\to 4m
$$

and,

$$
\frac{1}{2\pi} \varepsilon_{1,k} \vert x \vert^{2N} e^{v_k} \rightharpoonup\displaystyle{4\sum_{i=1}^m}\delta_{z_j} \quad\mbox{ weakly in the sense of measures on compact sets of } \mathbb{R}^2\setminus S_1,
$$

($\varepsilon_{1,k}$ defined in \eqref{eq53bis}).

In particular, from \eqref{lim}, we deduce the following identities:
\begin{equation}\label{eq562}
\frac{1}{2\pi} \int_{\mathbb{R}^2}e^{av}=\frac{a\beta(\beta-4(N+1))}{4(1-a(N+1))},
\end{equation}
\begin{equation}\label{eq563}
\frac{\beta(4-a\beta)}{4(1-a(N+1))}=\lim_{k\to+\infty}\frac{1}{2\pi}\int_{\mathbb{R}^2}e^{u_k}= \lim_{k\to+\infty}\frac{1}{2\pi}\int_{\mathbb{R}^2}\varepsilon_{1,k}\vert x \vert^{2N} e^{v_k}=4m; 
\end{equation}

with $m\geq 1$, and $\beta= \frac{1}{2\pi}\int_{\mathbb{R}^2}e^{a v}+4m$. Furthermore the limiting function $v$ satisfying:

\begin{equation}\label{eq564}
\begin{cases}
-\Delta v=e^{av}+8\pi\displaystyle{\sum_{i=1}^m}\delta_{z_i}, \, \mbox{ in } \mathbb{R}^2\\
\quad e^{av} \in L^1 (\mathbb{R}^2)
\end{cases}
\end{equation}

From \eqref{eq562} we see that $\beta>4(N+1)$. In addition, we can use \eqref{eq563} to solve for $\beta$, and obtain \eqref{eq516}.

To establish \eqref{eq516bis}, we note first that, if we take $m=1$ in \eqref{eq564}, then we can use \eqref{inpa}, with $w(x)=v(x+z_1)+4\log\vert x \vert$,  $b=a$, and $N=-2a$, in order to find that, 

\begin{equation}\label{eq566}
\frac{1}{2\pi}\int_{\mathbb{R}^2}e^{av}=\frac{4}{a}-8,
\end{equation}
and so, $\beta = \frac{4}{a}-4$ in this case.
But this is impossible, since for  $\beta = \frac{4}{a}-4$ and $m=1$ the identity  \eqref{eq563} is not satisfied.  Thus, we must have that, $m\geq2$. Furthermore,  by recalling that $\beta > 4(N+1)$,  then for $\frac{1}{2(N+1)} \leq a < \frac{1}{N+1}$ we can can derive \eqref{eq516bis} directly from \eqref{eq516}. To treat the case: $0 < a < \frac{1}{2(N+1)}$ we relate problem \eqref{eq564} to the construction of conformal metrics in the Riemann sphere with constant curvature equal to one and conical singularities.

To this purpose let us introduce the function:

\begin{equation}\label{eq565}
u(x)=\frac{a}{2}v(x)+\frac{1}{2}\log(\frac{a}{2}),
\end{equation}

satisfying:

\begin{equation}\label{eq566}
\begin{cases}
-\Delta u=e^{2u}+4\pi a\displaystyle{\sum^{m}_{j=1}}\delta_{z_j} \mbox{ in } \mathbb{R}^2\\
\frac{1}{2\pi}\displaystyle{\int_{\RR^{2}}e^{2u}dx=\frac{a}{2} (\beta-4m)},
\end{cases}
\end{equation}

with $\beta$ given in \eqref{eq516}. In other words,

\begin{equation}\label{eq570}
\frac{1}{2\pi}\int_{\mathbb{R}^2} e^{2u} dx=2-2am-\frac{4am(1-a(N+1))}{1+\sqrt{1-4am(1-a(N+1))}},
\end{equation}

Furthermore, 

\begin{equation}\label{eq567}
u(x)=-\frac{a}{2}\beta\log \vert x \vert + O(1), \, \mbox{ as } \vert x \vert \to+\infty.
\end{equation}

Therefore, (by using complex notations) via "compactification" of $\mathbb{C}$, we can consider the Riemann sphere: $M=\mathbb{C}\cup\left\{ +\infty \right\}$, with the conformal metric $g = e^{2u} dx$, so that $(M,g)$ admits conical singularities at the point $z_j$ with relative angle: $\alpha_j=2\pi(1-2a)\in(0,2\pi), \, j=1,...,m$ and at $\infty$, with angle:

\begin{equation}\label{eq568}
\alpha_{\infty}=2\pi(1+(\frac{a}{2}\beta-2))=2\pi\left( 1-\frac{4am(1-a(N+1))}{1+\sqrt{1-4am(1-a(N+1))}}\right)\in(0,2\pi),
\end{equation}

and away from the singularities, the surface admits constant curvature $K \equiv 1$. In this way, we see that \eqref{eq570} just expresses the Gauss-Bonnet formula for a compact surface $(M,g)$ with conical singularities, as it can be formulated as follows:

\begin{equation}\label{eq569}
\frac{1}{2\pi}\int_{\mathbb{R}^2} e^{2u} dx = \frac{1}{2\pi}\int_{M} K d v_g=\chi(M) + \sum_{j=1}^m \left(\frac{\alpha_j}{2\pi}-1\right) + \left(\frac{\alpha_\infty}{2\pi}-1\right)
\end{equation}

with $\chi(M)$ the Euler characteristic of $M$, and so $\chi(M)=2$ in our case.

Incidentally observe that, for $\frac{1}{2(N+1)}\leq a<\frac{1}{N+1}$, the right hand side of \eqref{eq570} defines a decreasing function with respect to $m$, and it vanishes exactly at $m=N+1$. So the condition $2 \leq m < N+1$ is again easily established in this case. More generally, since each angle at the singularities is positive but less than $2\pi$, on the basis of well known results of Troyanov \cite{troy2} and Luo-Tian \cite{lt} , we can formulate necessary and sufficient conditions on the values: $\theta_i=\frac{\alpha_i}{2\pi}-1=-2a$, $i=1,..,m$, and $\theta_{\infty}=\frac{\alpha_{\infty}}{2\pi}-1=-\frac{4am(1-a(N+1))}{1+\sqrt{1-4am(1-a(N+1))}}$, in order to ensure the existence of such singular sphere-like surface, or equivalently to guarantee the solvability of \eqref{eq566}. More precisely, from \cite{troy2} and \cite{lt} we see that problem \eqref{eq566} admits a solution if and only if:

\begin{equation}\label{eq571}
\theta_i>\sum_{i\neq j=1}^m \theta_j + \theta_\infty, \quad \forall \, i=1,...,m 
\end{equation}

and

\begin{equation}\label{eq571bis}
\theta_{\infty}>\sum_{j=1}^m \theta_j,
\end{equation}

(cfr. \cite{troy2}, \cite{lt}).

We easily check that for $\theta_i$ and $\theta_\infty$ specified above, $m \in \mathbb{N}$ and $0 < a < \frac{1}{N+1}$ the conditions \eqref{eq571} and \eqref{eq571bis} are exactly equivalent to the condition: $2 \leq m < N+1$, as claimed in $(a)$.\\

Finally, we treat the case where $0 \in S_1$. We know that this situation may occur only for $0 < a < \frac{1}{2(N+1)}$ (see Claim 1).\\

{\it \textbf{Claim 4}: if $0 < a < \frac{1}{2(N+1)}$ and $0 \in S_1$ then $S_1 = \{ 0 \}$ and $\beta = \frac{4}{a} - 4(N+1)$.}

By virtue of Claim 2, if by contradiction we let $S_1 \setminus \{ 0 \} = \{ z_1,...,z_m \}$ then for $\beta_j = \beta_0 (z_j)$ given in \eqref{eq515}, we have $\beta_j < \frac{2}{a} \; \forall \; j=1,...,m$. In other words, \eqref{eq524} holds and we can argue exactly as in Lemma \ref{lem3024bis} and Lemma \ref{lemma3029} in order to rule out the possibility that ''concentration''  occurs.

Therefore, as above, we see that (along a subsequence) we have:

\begin{equation}\label{eq7151}
v_k \to v \mbox{ uniformly in } C^{2}_{loc}(\mathbb{R}^2\setminus S_1)
\end{equation}

with $v$ satisfying:

\begin{equation}\label{eq715}
\begin{cases}
-\Delta v = e^{av} + 8\pi (N+1)\delta_0 + 8\pi \displaystyle{\sum^{m}_{j=1}}\delta_{z_j} \mbox{ in } \mathbb{R}^2\\
\frac{1}{2\pi}\displaystyle{\int_{\RR^{2}}e^{av}dx=\frac{a\beta(\beta-4(N+1))}{4(1-a(N+1))}}.
\end{cases}
\end{equation}

Moreover, as in \eqref{eq563}, we find:

\begin{equation}\label{eq716.1}
\frac{\beta(4- a \beta)}{4(1-a(N+1))}=4(N+1+m) \mbox{ and } \beta= \frac{1}{2\pi} \int_{\RR^{2}} e^{av} dx +4(N+1+m).
\end{equation}

By solving for $\beta$ in the first identity of \eqref{eq716.1}, we obtain:

\begin{equation}\label{eq717.1}
\beta=\frac{2}{a}\left(1+\sqrt{1-4a(N+1+m)(1-a(N+1))} \right).
\end{equation}

To ensure the integrability of $e^{a v}$, we must have that $\beta > \frac{2}{a}$, and so $m \in \mathbb{N}$ must satisfy:

\begin{equation}\label{eq718.1}
1\leq m <\frac{(1-2a(N+1))^2}{4a(1-a(N+1)}
\end{equation}

Notice that \eqref{eq718.1} already requires that $a>0$ is sufficiently close to zero. Furthermore as above, by a solution $v$ of \eqref{eq715} we can define $u$ via \eqref{eq566} and obtain the conformal factor for a metric $g= e^{2u} d\,x$ over the Riemann sphere $M=\mathbb{C}\cup\left\{+\infty\right\}$, which admits conical singularities at the origin with angle $\alpha_0=2\pi(1-2a(N+1))\in(0,2\pi)$, at the point $z_j$ with angle $\alpha_j=2\pi(1-2a)\in(0,2\pi), \, j=1,...,m$, and at $\infty$, with angle: 

$$\alpha_{\infty}=2\pi\left( 1-\frac{4a(N+1+m)(1-a(N+1))}{1+\sqrt{1-4a(N+1+m)(1-a(N+1))}}\right)\in(0,2\pi).$$

while, away from the singularities, the surface $M$ admits Gauss curvature $K \equiv 1$. As above, we obtain the following Gauss-Bonnet formula for the given "sigular" surface $(M, g)$:

\begin{small}
\begin{equation}\label{eq621bis}
\frac{1}{2\pi} \int_{\RR^{2}} e^{av} dx = \frac{1}{2\pi} \int_M K d\,\sigma_g = 2-2a(N+1+m)-\frac{4a(N+1+m)(1-a(N+1))}{1+\sqrt{1-4a(N+1+m)(1-a(N+1))}}.
\end{equation}
\end{small}

As already mentioned, from \cite{troy2} and \cite{lt} we obtain the following necessary and sufficient conditions on the values: $\theta_\infty = \frac{\alpha_\infty}{2 \pi} - 1$, $\theta_0 = \frac{\alpha_0}{2 \pi} - 1$ and $\theta_j = \frac{\alpha_j}{2 \pi} - 1$, $j=1,...,m$, for the solvability of \eqref{eq715}:

\begin{equation} \label{eq622a} 
\theta_{\infty} > \theta_0 + \sum_{j=1}^{m} \theta_j \; , \quad \theta_{0} > \theta_\infty + \sum_{j=1}^{m} \theta_j \; , \quad \theta_i > \sum_{i \ne j=1}^{m} \theta_j + \theta_{\infty} + \theta_{0}, \quad \forall \, i=0,...,m
\end{equation}

In other words, problem \eqref{eq715} is solvable with $\beta$ satisfying \eqref{eq717.1}, if and only if, for $m \in \mathbb{N}$ satisfying \eqref{eq718.1}, there holds:

\begin{eqnarray}
\label{eq6220} &\theta_{\infty}&:= \frac{-4a(N+1+m)(1-a(N+1))}{1+\sqrt{1-4a(N+1+m)(1-a(N+1))}}>-2a(N+1+m),\\
\label{eq6221} &-2a(N+1)& > - 2 a m + \theta_\infty,\\
\label{eq6222} &-2a& > - 2 a (N+m) + \theta_\infty
\end{eqnarray}

But actually we can check that, under the given assumptions, the inequality \eqref{eq6220} is \underline{never} satisfied. Thus we conclude that necessarily $S_1=\left\{0\right\}$ as claimed.

Furthermore, by virtue of \eqref{eq54}, we can rule out again the possibility that ''concentration'' may occur. As a consequence, we find that $\beta$ must satisfy \eqref{eq717.1} with $m=0$, which for $0 < a < \frac{1}{2(N+1)}$ gives that: $\beta = \frac4a - 4(N+1)$.

Equivalently, since $v$ in \eqref{eq7151} now satisfies:

$$
\begin{cases}
-\Delta v=e^{av}+8\pi (N+1)\delta_0 \mbox{ in } \mathbb{R}^2\\
\frac{1}{2\pi} \int_{\RR^{2}}e^{av}dx < \infty.
\end{cases}
$$

we can use \eqref{inpa} (applied with $b=a$ and $w(x)=v(x)+4(N+1)\log \vert x \vert$) to find that: $\frac{1}{2\pi} \int_{\RR^{2}}e^{av}dx = \frac{4}{a}-8(N+1)$ and so $\beta=\beta_0+4(N+1)=\frac{4}{a}-4(N+1)$, as claimed.\\

At this point, the proof may be completed simply by a direct consequence of the four Claims established above. \qed

\begin{coro}
If $S\neq\emptyset$ then $\beta = \max \left\{ 4 (N+1), \frac{4}{a}-4(N+1) \right\}$.
\end{coro}

\dimo
In case $\frac{1}{2(N+1)} \leq a < \frac{1}{N+1}$ then the desired conclusion: $\beta = 4(N+1)$, has been already established in part $(ii)$ of Corollary \ref{cor14bbis}. Hence, let $0 < a < \frac{1}{2(N+1)}$, and observe that, if $S \neq \emptyset$ then $0\in S_1$. Indeed if $x_0\in S$, then $\exists \, x_k\to x_0: \, u(x_k)\to \infty$. Therefore if we let, $z_k=\frac{x_k}{R_k}$ we find: $v_k(z_k)=u_k(x_k)+\frac{2}{a}\log R_k \to +\infty$ and $z_k\to 0$, as $k\to+\infty$. So now the desired conclusion follows by part $(ii)$ of Theorem \ref{teo5.6}.
\qed\\

On the other hand, when the blow--up set $S$ of $u_k$ is \underline{empty}, we may also conclude the following:

\begin{coro}\label{cor64}
If $S$ (blow--up set of $u_k$) is empty but $S_1$ (blow--up set of $v_k$ in \eqref{eq52}) is not empty, then ''concentration'' occurs for $v_k$ if and only if:

\begin{equation}\label{eq530}
S_1 \subset \mathbb{R}^{2} \setminus \{ 0 \} \text{ and } \exists \, x_0 \in S_1: |x_0| = 1 \text{ with }  \beta_0 (x_0) \geq \frac2a
\end{equation}

($\beta_0 (x_0)$ defined in \eqref{eq56bis}) unless $\frac{1}{2(N+1)} \leq a < \frac{1}{N+1}$, $N \in \mathbb{N}$, $\beta = 4(N+1)$, and (up to rotation) $S_1$ is formed by the $N+1$ roots of the unity.
\end{coro}

Actually, when \eqref{eq530} occurs, we can complete the description of the blow--up behavior of $v_k$ as follows

\begin{coro}\label{coro57bis}
Let (in complex notation) $S_1 = \{ z_1,...,z_m \}\subset \mathbb{R}^{2} \setminus \{ 0 \}$, and assume that $\beta_1 = \beta (z_1) \geq \frac2a$, ( i.e. \eqref{eq530} holds ).  Then either $m=1$ and $\beta = \frac4a$, or $N>1$,  $m \geq 2$,  and the following holds:

\begin{eqnarray}
\label{eq533*} \frac{N (1 - a \left( \beta_1 / 4\right))}{2(1-a)} = \sum_{j=2}^{m} \frac{z_1}{z_1 - z_j}\\
\notag \frac{N}{2} = \sum_{i \ne j=2}^{m} \frac{z_i}{z_i - z_j} + \frac{\beta_1}{4} \sum_{j=2}^{m} \frac{z_i}{z_i - z_1},\; j=2,...,m.
\end{eqnarray}

Furthermore,

{\footnotesize $$\beta_1 = \frac{2}{a} \left( 1 - \frac{2(m-1)(1-a(N+1))}{N}  + \sqrt{\left( 1 - \frac{2(m-1)(1-a(N+1))}{N} \right)^2 + \frac{4(m-1)ma}{N} (1-a(N+1))} \right)$$}

and \eqref{eq517}, \eqref{eq518} and \eqref{eq546bis} hold.
\end{coro}

\dimo  Since ''concentration'' occurs for $v_k$, in the sense that, 

$$
\frac{1}{2\pi} f_{1,k} \rightharpoonup \beta_1 \delta_{z_1} + 4 \sum_{j=1}^{m} \delta_{z_j} \quad\mbox{ weakly in the sense of measure }
$$

then we can simply argue as in the proof of Lemma \ref{lem3024bis} in order to deduce \eqref{eq533*}, we leave the details to the interested reader. The remaining part of the statement has been established above.
\qed

\begin{rem}
Observe that if $u_k$ is \underline{radially symmetric}, then either $S_1 = \emptyset$ or $0 \in S_1$, since the radially symmetric function $v_k$ attains its maximum value at the \underline{origin}. As a consequence, in the radially symmetric case, we get that either $\beta = \frac{4}{a}$ or $\beta = \frac{4}{a} - 4 (N+1)$, consistently with what has been established in \cite{pt2}.  
\end{rem}

Next we wish to investigate further the nature of the conditions \eqref{eq517} and \eqref{eq518}. To this purpose, for $a >0$ we consider the function:

\begin{equation}\label{eq530*}
f(a) := \frac{1}{2a} \left( \sqrt{(1-a(N+1))^2 + \frac{Na}{1-a}}- (1 - a(N+1)) \right).
\end{equation}

 In the {\bf Appendix} we prove the following:\\

\textbf{Calculus Lemma}: Let $N \geq 1$. The function $f(a)$ defined in \eqref{eq530*} is strictly monotone increasing for $a \in \left( 0 , \frac{1}{N+1} \right)$. Furthermore by setting:

\begin{equation}\label{eq532}
a_N = \begin{cases} \frac{4-N}{2(N+1+\sqrt{(N-1)^2+N^2}} \quad &1 \leq N < 4\\
0 &N \geq 4
\end{cases}
\end{equation}

we have:

\begin{enumerate}[(i)]
\item $f(a) \geq 1$ if and only if $N>1$ and $a_N \leq a < \frac{1}{N+1}$
\item If $1 < N \leq 3$ and $a \in \left( 0, \frac{1}{N+1} \right)$ then $f(a) < 2$.
\end{enumerate}

\begin{rem}
From part $(i)$ of the Calculus Lemma it follows tha,t for $1 \leq N < 4$ and $0 < a \leq a_N$ ($a_N$ defined in \eqref{eq532}) then in part $(iii)$ of Theorem \ref{teo5.6} only the \underline{first} alternative can occur. While, if $1 < N \leq 3$, $a_N \leq a < \frac{1}{N+1}$, then in part $(iii)$ of Theorem \ref{teo5.6} the second alternative can occur only with $m=2$, and we can use \eqref{eq533*} in order to locate the blow--up points.
\end{rem}

In conclusion, we can summarize the results established above as follows:

\begin{coro}
If $N \in (0,1)$, $0 < a < \frac{1}{N+1}$ and $S = \emptyset$ then $\beta = \frac{4}{a}$.
\end{coro}

\begin{coro}
If $N \geq 1$, $0 < a < \frac{1}{N+1}$ and $S = \emptyset$, then $\beta$ takes one of the following values:

\begin{enumerate}[i)]
\item $\beta = \frac{4}{a}$,

\item $\beta = \frac{2}{a} (1+\sqrt{1-4 a m (1-a (N+1))})$, with $m \in \mathbb{N}$ and $2 \leq m < N+1$; in particular $N>1$, in this case;

\item $\beta = \frac{2}{a} \left( 1 - \frac{2(m-1)(1-a(N+1))}{N}  + \sqrt{\left( 1 - \frac{2(m-1)(1-a(N+1))}{N} \right)^2 + \frac{4(m-1)ma}{N} (1-a(N+1))} \right)$,

with 

\begin{equation*}
m \in \mathbb{N}, \; 1 \leq m-1 < \frac{1}{2a} \left( \sqrt{(1-a(N+1))^2 + \frac{Na}{1-a}}- (1 - a(N+1)) \right)
\end{equation*}

and \eqref{eq546bis} holds. In particular $N > 1$ and $a_N \leq a < \frac{1}{2(N+1)}$ ($a_N$ defined in \eqref{eq532}), in this case .

\item $\beta = 4 (N+1)$, provided that $N \in \mathbb{N}$ and $\frac{1}{2(N+1)} \leq a < \frac{1}{N+1}$ or
$\beta = \frac{4}{a} - 4 (N+1)$, provided that $0 < a < \frac{1}{2(N+1)}$.

\end{enumerate}
\end{coro}
\qed

\section{The case $\frac{1}{N+1}<a < 1$}
\setcounter{equation}{0}

In view of part $(ii)$ of Corollary \ref{cor14bbbis}, we know that if $0 \notin S$ then $\beta = 4 (N+1)$. Hence, we need to investigate what happens when we assume, 

\begin{equation}\label{eq60}
0 \in S.
\end{equation}

Let,

\begin{equation}\label{eq71}
\beta_0 := \lim_{r \to 0} \lim_{k \to +\infty} \frac{1}{2\pi} \int_{\vert x \vert \leq r} f_{k} (x) \, dx \; 
\end{equation}

so that,

\begin{equation}\label{eq71*}
\beta_0 \in \left[ \frac{4}{a}, 4 (N+1) \right].
\end{equation}

For small $\varepsilon > 0$ (to be specified later) we let $r_k > 0$ be such that:

\begin{equation}\label{eq76}
\frac{1}{2\pi} \int_{\vert x \vert \leq r_k} f_{k} (x) \, dx = \frac{4}{a} - \varepsilon
\end{equation}

Clearly by \eqref{eq71} and \eqref{eq76}, we see that

\begin{equation}\label{eq77}
r_k \to 0, \quad \text{ as } k \to +\infty.
\end{equation}

We set,

\begin{equation}\label{eq78}
v_k (x) = u_k (r_k x) + \frac{2}{a} \log r_k
\end{equation}

which satisfies \eqref{eq53} with, 

\begin{equation}\label{eq712}
\varepsilon_{1,k} = r_k^{\frac{2}{a} (a(N+1) -1)} \to 0, \; \text{ as } k \to +\infty
\end{equation}

We use the notations of the previous section, so for example, we let

\begin{equation}\label{eq713}
\hat{v}_k (x) = u_k \left( \frac{x}{|x|^2} \right) + \beta_k \log \frac{1}{|x|}
\end{equation}

which satisfies \eqref{eq56}. Also we denote with $S_1$ and $\hat{S}_{1}$ the (possibly empty) blow--up set of $v_k$ and $\hat{v}_{k}$ respectively. Clearly Lemma \ref{lem50} applies to $v_k$. 

We start to analyze the case where: 

\begin{equation}\label{eq70}
\frac{1}{N+1}<a \ne 1 \leq \min \left\{ 1, \frac{2}{N+1} \right\} \quad \text{ and } \quad  0 \in S.
\end{equation}

\begin{lemma}\label{lem71}

Assume \eqref{eq70}, then for $\varepsilon > 0$ sufficiently small,

\begin{equation}\label{eq79}
\lim_{R \to +\infty} \lim_{k \to +\infty} \frac{1}{2\pi} \int_{\vert x \vert \geq R} f_{1,k} = 0
\end{equation}

In particular, for $\beta > 2(N+1)$ we have: $0 \notin \hat{S}_1$. 
\end{lemma}

\underline{Proof:} We start to show the following:\\

{\it \textbf{Claim:} If $\frac{1}{N+1}<a < \min \left\{ 1, \frac{2}{N+1} \right\}$ or $N>1$, $a=\frac{2}{N+1}$ and $\beta > 2(N+1)$ then $0 \notin \hat{S}_1$.}\\

We argue by contradiction and assume that $0 \in \hat{S}_1$. In view of our assumption, we can use part $(ii)$ of Lemma \ref{lem50} together with \eqref{eq76} in order to conclude that,

\begin{equation}\label{eq716}
\beta + \varepsilon - \frac{4}{a} \geq \beta_\infty \geq 2 \beta - 4(N+1)
\end{equation}

Therefore, if $\frac{1}{N+1} < a < \frac{2}{N+1}$ then from \eqref{eq716} we deduce that,

\begin{equation}\label{eq716bis}
\frac4a \leq \beta < 4(N+1) - \frac4a + \varepsilon
\end{equation}

which is impossible, provided we choose $\varepsilon \in \left( 0,\frac{8}{a} - 4(N+1) \right)$.

On the other hand, for $N > 1$, $\frac{2}{N+1} = a < 1$ and $\beta > 2(N+1)$ it suffices to choose $\varepsilon \in ( 0,\beta - 2(N+1))$, to obtain from \eqref{eq716bis} a contradiction as well.

So it remains to analyze the following situation,

\begin{equation}\label{eq717}
N > 1, \quad  a=\frac{2}{N+1} < 1 \text{ and } \beta = 2 (N+1) = \frac{4}{a}
\end{equation}

By virtue of \eqref{eq76}, from \eqref{eq717} we find that

$$\frac{1}{2\pi} \int_{|x| \geq 1} f_{1,k} (x) \, dx = \varepsilon$$

Therefore for $\varepsilon > 0$ sufficiently small, we can argue as in Lemma \ref{lem42}, in order to deduce that,

\begin{equation}\label{eq718}
\lim_{R \to +\infty} \lim_{k \to +\infty} \int_{\vert x \vert \geq R} e^{a v_k} = 0.
\end{equation}

On the other hand, since for $\beta = \frac{4}{a}$ we have that,

$$\varepsilon_{1,k} \int_{\mathbb{R}^2} |x|^{2N} e^{v_k} = \int_{\mathbb{R}^2} |x|^{2N} e^{u_k} \to \frac{\beta (a\beta -4)}{4(a(N+1)-1)} = 0,$$ 

and we conclude that  \eqref{eq79} holds in this case as well. \qed

\begin{prop}\label{prop71}
Assume \eqref{eq70} and let $v_k$ be defined in \eqref{eq78} with $S_1$ its (possibly empty) blow--up set.

\begin{enumerate}[i)]

\item $0 \notin S_1$, 

\item if $S_1$ contains at most one point then $\beta = \frac{4}{a}$,

\item if $S_1 = \{ z_1,..., z_m \} \subset \mathbb{R}^2 \setminus \{ 0 \}$ with $m \geq 2$, then $N \geq 1$ and for  $\beta_j = \beta_0 (z_j)$ defined in \eqref{eq515}, $j=1,...,m$, the following holds:

\begin{enumerate}[a)]
\item ''concentration'' occurs (i.e. \eqref{eq522} holds) if and only if

\begin{equation}\label{eq81}
\max_{1 \leq j \leq m} \beta_{j} \geq \frac{2}{a} \; \text{ and } \; 2 \leq m < N+1,
\end{equation}

unless $N \in \mathbb{N}$, $m=N+1$, $\beta = 4 (N+1)$ and in this case $S_1$ is formed (up to rotation) by the $N+1$ roots of the unity.

\item In case ''concentration'' \underline{does not} occur, then $N > 1$, $a \in \left( \frac{1}{N+1}, \frac{2}{N+1} \right] \cap \left( 0, \frac{1}{2} \right)$, $\beta_{j} = 4$, $\; \forall \; j = 1,...,m$; and \eqref{eq516}, \eqref{eq516bis} hold for $\beta$.
\end{enumerate}

\end{enumerate}

\end{prop}

\dimo

Recall that,
 
\begin{equation}\label{eq80}
\frac{1}{2\pi} \int_{\vert x \vert \leq 1} f_{1,k} (x) \, dx = \frac{4}{a} - \varepsilon,
\end{equation}
 so by \eqref{eq5.14.1} we deduce that necessarily $0 \notin S_1$.

Concerning $ii)$, we see that in case $S_1 = \emptyset$, then we argue essentially as in Theorem \ref{teo5.6} and by taking into account \eqref{eq5.36bis} we may conclude that $\beta = \frac{4}{a}$ as claimed. Similarly, in case $S_1 = \{ z_1 \}$, we check first that ''concentration'' must occur. Indeed, if otherwise, then we would get, as in Theorem \ref{teo5.6}, that $\beta=\frac{4}{a}-4$, which is impossible for $a>\frac{1} {N+1}$. Consequently,  \eqref{eq65} holds and yields to the following identity: $\frac{\beta_1 (\beta_1 - 4)}{1-a} =  \frac{\beta (4(N+1) - \beta)}{a(N+1)-1}$, which must hold with $\beta = \beta_1$.  But this is possible only for $\beta = \frac{4}{a}$.

Next, let $S_1 = \{ z_1,..., z_m \} \subset \mathbb{R}^2 \setminus \{ 0 \}$ with $m \geq 2$ and $4 \leq \beta_j \leq \frac{4}{a}$, for all $j=1,...,m$.

We observe that "concentration" occurs (or equivalently \eqref{eq522} holds) if and only if  (by taking into account \eqref{eq5.14.2}) the following holds:

\begin{equation}\label{eq83}
\beta = \sum_{j=1}^{m} \beta_{j}, \quad \text{ and } \quad\quad \frac{1}{1-a} \sum_{j=1}^{m} \beta_{j} (\beta_{j}-4) = \frac{\beta (4(N+1) - \beta)}{a(N+1)-1}.
\end{equation}

In case $\beta_j = 4$, for every $j=1,...,m$, then necessarily $N \in \mathbb{N}$, $m = N+1$ and $z_{j}^{N+1} = \xi$, with suitable $\xi \ne 0$.

But from \eqref{eq80}, we see that necessarily $|\xi| = 1$, and so (up to rotation) $S_1$ is formed by the $(N+1)$-roots of the unity. On the other hand in case,

$$\beta_* = \max \{ \beta_j, j=1,...,m \} > 4$$

then according to \eqref{eq5.18} we see that necessarily $\beta_* \geq \frac{2}{a}$, and $iii)-a)$ is established .

In case ''concentration'' does \underline{not} occur then $\beta_j = 4$, for every $j=1,...,m$,  and as in Theorem \ref{teo5.6} we deduce that,  $\beta<4(N+1)$ and 

\begin{equation}\label{eq84}
\beta = \frac{2}{a} (1+\sqrt{1+4 a m (a (N+1)-1)})
\end{equation}

Since for $a>\frac{1}{N+1}$, the expression in the right hand side of \eqref{eq84} is monotone increasing with respect to $m$, and it attains the value $4(N+1)$ for $m=N+1$, we find that necessarily $2 \leq m < N+1$ as claimed. \qed

\begin{rem}
Notice that for $1 \leq N \leq 3$ and $a \in \left[ \frac{1}{2}, \frac{2}{N+1} \right]$ then only alternative $a)$ can occur in $iii)$. Actually, it is interesting  to investigate whether alternative $b)$ can occur at all, also in view of the geometric problem \eqref{eq566}, which relates to
the existence of a conformal metric in the Riemann sphere with constant curvature equal to one and conical singularities. Indeed, we see that for $a > \frac{1}{N+1}$ the angle $\alpha_\infty$ in \eqref{eq568}, corresponding to the singularity at $\infty$,  is bigger than $2 \pi$, and consequently, the conditions \eqref{eq571} and \eqref{eq571bis} are only sufficient but no longer necessary to guarantee the existence of such metric. 
In such situation, more involved \underline{necessary} conditions about the angles at the singularities have been introduced by \cite{cwx},  \cite{er1,egt1,egt2,egt3} and \cite{uy} in some specific cases. However, it is not clear yet which should be the appropriate set of \underline{necessary} conditions  in general, and in fact, it may be possible to rule out the occurence of \eqref{eq84} in certain cases.
\end{rem}

By combining the results above with those in Lemma \ref{lemma2} we conclude:

\begin{coro}\label{cor72}
Let $\frac{1}{N+1} < a \leq \frac{2}{N+1}$ and suppose that $S\neq\emptyset$ we have:

\begin{enumerate}[i)]
\item if $N \in (0,1)$ then $\beta = \frac{4}{a}$;
\item if $N = 1$ then either $\beta = \frac{4}{a}$ or $\beta = 4(N+1)=8$ and in the latter case $S_1 = \{ z_1,-z_1\}$ with $|z_1|=1$.
\end{enumerate}
\end{coro}
\qed

\begin{rem}
We mention that, by following the arguments in \cite{blt} it is possible to rule out the second alternative in $ii)$, and show that actually the following holds:

\begin{equation}\label{eq85}
\text{ if } \, N \in (0,1] \, \text{ and } \,\frac{1}{N+1} < a \leq \frac{2}{N+1} \, \text{ then } \, \beta = \frac{4}{a}.
\end{equation}

Details will be given in a forthcoming paper.
\end{rem}

\begin{prop}\label{prop73}
Let $N>1$ and $\frac{1}{N+1} < a \leq \frac{2}{N+1}$ then one of the following alternatives holds:

\begin{enumerate}[i)]

\item $\beta = \frac{4}{a}$;

\item $N \in \mathbb{N}$ and $\beta = 4(N+1)$ 

\item $\beta \in \left( \frac{4}{a}, 4(N+1) \right)$ and for some $m \in \mathbb{N}$,

\begin{enumerate}[a)]

\item \underline{either} $a \in \left( \frac{1}{N+1}, \frac{1}{2} \right)$ and 

\begin{equation}\label{eq86}
\beta = \frac{2}{a} (1+\sqrt{1+4 a m (a (N+1)-1)}), \text{ with } 2 \leq m < N+1
\end{equation}

\item \underline{or} $\beta = \sum_{j=1}^{m} \beta_{j}$ with $\beta_j \in \left[ 4, \frac{4}{a}\right]$ and $\max_{j=1,...,m} \beta_j \geq \frac{2}{a}$;
furthermore:

\begin{equation}\label{eq87}
\sum_{j=1}^{m} \beta_{j}^{2} = \frac{\beta (4N - (1-a)\beta)}{a(N+1)-1}, \; \text{ with } 2 \leq m < N+1 - \frac{1}{2a} \max \{ 0, 1-2a \}
\end{equation}

\end{enumerate}
\end{enumerate}

\end{prop}

\dimo By virtue of Proposition \ref{prop71}, it  remains only to check  \eqref{eq87}. To this purpose, it suffices to use \eqref{eq83}, with $\beta = \sum_{j=1}^{m} \beta_j$. For what concerns the claimed upper estimate on $m \in \mathbb{N}$, we recall that $\beta<4(N+1)$ and from \eqref{eq86} and  \eqref{eq87} it is easy to see that necessarily, $2 \leq m < N+1$. On the other hand, in the alternative $iii)-b)$ and for  $a \in \left( \frac{1}{N+1}, \frac{1}{2} \right)$ there holds: 
$\frac{2}{a} + 4 (m-1) \leq \beta < 4(N+1)$, and this yields to the improved estimate for $m$ in \eqref{eq87}.

\qed

\begin{rem}
By using \eqref{eq87} and the fact that $\beta \in \left( \frac{4}{a}, 4(N+1) \right)$ we see that in alternative $iii)-(b)$ the following estimates holds:

\begin{equation}\label{eq88}
\begin{split}
&\text{ if } \, \frac{1}{N+1} < a \leq \frac{2}{N+2} \, \text{ then } \, 16(N+1) < \sum_{j=1}^{m} \beta_{j}^{2} < \left( \frac4a \right)^2\\
&\text{ if } \, \frac{2}{N+2} < a \leq \frac{2}{N+1} \, \text{ then } \, \min \left\{ 16(N+1), \left( \frac4a \right)^2 \right\} < \sum_{j=1}^{m} \beta_{j}^{2} < \frac{(2N)^2}{(1-a)(a(N+1)-1)}.
\end{split}
\end{equation}

\end{rem}

 In the following section, we shall construct an explicit example where the alternative $iii)-b)$ actually occurs with all $\beta_{j}$'s that coincide, namely:  $\beta_1 =...= \beta_m \geq \frac2a$ for every $j=1,...,m$.

For this reason we point out the following general facts concerning such situation.

\begin{coro}\label{cor74}
In altenative $iii)-b)$ of Proposition \ref{prop73} we have that  \underline{all} the $\beta_{j}$'s coincide (i.e. $\beta_1=\beta_2=...=\beta_m$) if and only if, 

\begin{equation}\label{eq89bis}
\beta = \frac{4Nm}{(a(N+1)-1) + m (1-a)}, \: \text{ with } \: m \in \mathbb{N}, \: 2 \leq m < N+1 - \frac{N}{1-a} \max \{ 0, 1-2a \}.
\end{equation}

In particular,

\begin{equation}\label{eq89}
\beta_1 = \beta_2=....=\beta_m=\frac{4N}{(a(N+1)-1) + m (1-a)} \geq \frac{2}{a},
\end{equation}

and the blow--up set $S_1$ of $v_k$ is formed (up to rotation) by the $m$-roots of the unity.
\end{coro}

\underline{Proof}: We observe first that, $\sum_{j=1}^{m} \beta_{j}^{2} \geq \frac{1}{m} \beta^2$ and equality holds if and only if all the $\beta_j$'s coincide.

But from \eqref{eq87} we see that $\sum_{j=1}^{m} \beta_{j}^{2} = \frac{1}{m} \beta^2$ if and only if $\beta$ satisfies \eqref{eq89bis}, and so \eqref{eq89} must hold.

To obtain the given estimate on $m$ in \eqref{eq89bis}, we recall first that $2 \leq m < N+1$, but from \eqref{eq89} we find also that,
\begin{equation}\label{eq89*}
2 \leq m \leq \frac{(N-1)a+1}{1-a}
\end{equation}

and we arrive at \eqref{eq89bis}. Concerning the blow--up set $S_1 = \{ z_1,...,z_m \}\subset \mathbb{R}^2 \setminus \{ 0 \}$ we know that, 

$$\frac{1}{2\pi} f_{1,k} \rightharpoonup \beta_1 \sum_{j=1}^{m}  \delta_{z_j} \quad\mbox{ weakly in the sense of measure, locally in } \RR^2$$

and therefore we can argue exactly as in Lemma \ref{lem3024bis},  to derive that, as $k \to \infty$ 

$$\frac{1}{2\pi} \int_{B_{\delta} (z_j)} \nabla (\log |z|^{2N}) |z|^{2N} e^{v_k} \to \frac{2N \beta_1 (4 - a\beta_1)}{4(1-a)} \frac{z_i}{|z_i|^2} = \beta_{1}^{2} \sum_{i \ne j=1}^{m} \frac{z_i - z_j}{|z_i - z_j|^2}, \; \forall \;  i=1,...,m$$

Thus, by using complex notation and \eqref{eq89}, after explicit calculations, we find that,  

$$(m-1) = 2 \sum_{i \ne j=1}^{m-1} \frac{z_i}{z_i - z_j}, \; \forall \; i=1,2,...,m;$$

and this property just gives the equivalent of \eqref{starpoll2nu} with $m$ in place of $N+1$. Therefore, we can argue exactly as in Theorem \ref{new} in order to conclude that necessarily: $z_{j}^{m} = \xi_0, \; \forall \: j=1,...,m$, and with $\xi_0 = (-1)^{m-1} z_1 \cdot ... \cdot z_m$.

Moreover, in view of \eqref{eq80}, we see that necessarily: $|z_j|= 1$, $\forall \, j=1,...,m$, and the proof is completed. \\
 \qed

Finally, we complete our analysis by considering the case where \eqref{eq70} holds with,

\begin{equation}\label{eq90.1}
N > 1 \; \text{ and } \; \frac{2}{N+1} < a <1
\end{equation}

\begin{prop}
Assume that \eqref{eq60} and \eqref{eq90.1} hold. Then $\beta \geq 2(N+1)$ and either $\beta$ takes one of the values specified in alternative $ii)-iii)$ of Proposition \eqref{prop73} or it satisfies:

\begin{equation}\label{eq91.1}
2(N+1) \leq \beta \leq 4(N+1) - \frac{4}{a}.
\end{equation}

\end{prop}

\dimo
Clearly $\beta \geq 2(N+1)$, and in case $\beta = 2(N+1)$, then obviously it satisfies \eqref{eq91.1}.
Hence, let us assume that $\beta > 2 (N+1)$, and for $\varepsilon > 0$ sufficiently small, we consider the sequence $v_k$ defined by \eqref{eq76} and \eqref{eq78}, and observe that in this case, $0 \notin S_1$.

In case $S_1 = \emptyset$, then necessarily $0 \in \hat{S}_1$ and $\beta = \frac{4}{a} + \beta_\infty$ with $\beta_\infty$ in \eqref{eq5.14.5bis}. Furthermore, by \eqref{eq5.14.6}, we find that: $\beta_\infty \geq 2 \beta - 4 (N+1)$, and we conclude that \eqref{eq91.1} holds in this case. Next, we suppose that $S_1 \ne \emptyset$, and we distinguish between two situations, namely: 

\begin{equation}\label{eq92.1}
\text{ either } \, \exists \; \varepsilon > 0 \text{ as small as needed: } 0 \notin \hat{S}_1
\end{equation}

\begin{equation}\label{eq93.1}
\text{ or } \; \forall \, \varepsilon > 0 \text{ sufficiently small: } 0 \in \hat{S}_1 
\end{equation}

If \eqref{eq92.1} holds then $S_1 = \{ z_1,...,z_m \} \subset \mathbb{R}^2 \setminus \{ 0 \}$ with $m \geq 2$, and as above, we may conclude that in this case, $\beta$ must satisfy either one of the alternatives $ii)$ and $iii)$ in Proposition \eqref{prop73}.

In case \eqref{eq93.1} holds, then we can apply part $iii)$ of Lemma \ref{lem50} to $v_k$, and obtain that for $\beta_\infty$ in \eqref{eq5.14.5bis} we have:

$$\beta_\infty \leq \beta - \frac{4}{a} + \varepsilon \; \text{ and } \; \beta_\infty \geq 2 \beta - 4(N+1),\; \forall \varepsilon > 0;$$

and from such estimates we derive the upper bound on $\beta$ as given in \eqref{eq91.1}.

\qed

\subsection{An example}\label{example}

 We recall that for $a \in (0,1)$ the problem:

\begin{equation}\label{eq90}
\begin{cases}
-\Delta U = e^{aU} + e^U \quad\mbox{in }\mathbb{R}^2,\\
\displaystyle{\frac{1}{2\pi} \int_{\mathbb{R}^2}( e^{aU} + e^U ) \, dx} = \beta_0,
\end{cases}
\end{equation}

admits a solution if and only if: $\max \left\{ 4, \frac{4(1-a)}{a} \right\} < \beta_0 < \frac{4}{a}$. Furthermore, up to translation, the solution is radially symmetric about the origin and it is uniquely identified by $\beta_0$. We refer to \cite{pot,pt2} for details.

Hence, we fix $\beta_0$ such that:

\begin{equation}\label{eq91}
\max \left\{ 4, \frac{4(1-a)}{a} \right\} < \beta_0 < \frac{4}{a},
\end{equation}

and let $u_0 (z) = u_0 (r)$ be the corresponding unique radial solution for \eqref{eq90}. As usual we are using complex notations, and so for given $ m_1 \in \mathbb{N}$  we define:

\begin{equation}\label{eq92}
U_k (z) = u_0 (z^{m_1 + 1} - \xi_k) : \xi_k \in \mathbb{C} \setminus \{ 0 \} 
\end{equation}

Clearly, $U_k$ is {\bf {not}} radially symmetric about any point and it satisfies:
\begin{equation}
\begin{cases}
-\Delta U_k=(m_1 +1)^{2}\left|z\right|^{2 m_1}\left(e^{a U_k}+e^{U_k}\right),\\
\frac{1}{2\pi}\displaystyle{\int_{\RR^{2}}\left(m_1 +1\right)^{2}\left|z\right|^{2m_1}(e^{aU_k}+e^{U_k})}=(m_1 +1)\beta_0.
\end{cases}
\end{equation}

In turn, if we let

\begin{equation}\label{eq93}
u_k (z) := U_k \left(\frac{z}{\left|z\right|^{2}}\right)+(m_1 +1)\beta_0\log\left(\frac{1}{\left|z\right|}\right) = u_0 \left(\frac{1}{\bar{z}^{m_1+1}} - \xi_k \right) - \beta_0\log|\bar{z}^{m_1+1}|
\end{equation}

then $u_k$ can be extended smoothly at the origin, and it satisfies:

\begin{equation*}
\begin{cases}
-\Delta u_k =(m_1 +1)^{2}\left(\left|z\right|^{a(m_1 +1)\beta_0-2(m_1 +2)}e^{a u_k }+\left|z\right|^{(m_1 +1)\beta_0-2(m_1 +2)}e^{u_k }\right),\\
\frac{1}{2\pi}(m_1 +1)^{2}\displaystyle{\int_{\RR^{2}}\left(\left|z\right|^{a(m_1 +1)\beta_0-2(m_1 +2)}e^{au_k}+\left|z\right|^{(m+1)\beta_0-2(m+2)}e^{u_k}\right)=(m_1 +1)\beta_0.}
\end{cases}
\end{equation*}

Therefore, whenever it is possible to choose $\beta_0$ such that,

\begin{equation}\label{eq96}
\beta_0=\frac{2(m_1 +2)}{(m_1 +1)a}\in\left(\max\left\{4, \frac{4(1-a)}{a}\right\},\frac{4}{a}\right)
\end{equation}


with $0<a<1$ and $m_1 \in \NN$, and if we replace $u_k (z)$ with $u_k \left( \frac{z}{m+1} \right)$, then we are able to obtain a sequence of non--radial solutions for  the problem:

\begin{equation}\label{eq94}
\begin{cases}
-\Delta u_k = \left(e^{u_k}+\left|z\right|^{2N} e^{u_k} \right),\\
\frac{1}{2\pi} \int_{\RR^{2}}\left(e^{a u_k}+\left|z\right|^{2N}e^{u_k}\right)=\frac{2N}{1-a}
\end{cases}
\end{equation}

with

\begin{equation}\label{eq95}
N = N(m_1,a)=\frac{(m_1 +2)(1-a)}{a}.
\end{equation}
To this purpose, for $a \in \left( \frac{1}{4}, \frac{3}{4} \right)$ we set

\begin{equation}\label{eq97}
m_a = \begin{cases}
\frac{2a}{1-2a} \quad &\text{ if } \frac{1}{4} < a < \frac{1}{2}\\
+\infty \quad &\text{ if } a = \frac{1}{2}\\
\frac{2(1-a)}{2a-1} \quad &\text{ if } \frac{1}{2} < a < \frac{3}{4},
\end{cases}
\end{equation}

then $m_a = m_{1-a}$, and by direct calculations, we check that, \eqref{eq96} holds if and only if 

\begin{equation}\label{eq98}
a \in \left( \frac{1}{4}, \frac{3}{4} \right) \text{ and }m_1 \in \mathbb{N} : 1 \leq m_1 < m_a.
\end{equation}

At this point, by assuming \eqref{eq98}, we take, 

\begin{equation}\label{eq99}
\xi_k \in \mathbb{C}: | \xi_k | \to + \infty, \text{ as } k \to \infty
\end{equation} 

and we show that,

\begin{equation}\label{eq100}
u_k (z) := u_0 \left(\left(\frac{m_1 + 1}{\bar{z}}\right)^{m_1 + 1} - \xi_k \right) - \beta_0 \log \left| \left( \frac{\bar{z}}{m_1 + 1}\right)^{m_1 + 1} \right|
\end{equation}

defines a solution sequence of  \eqref{eq94}-\eqref{eq95}, 
which admits a blow--up point exactly at the origin, namely $S = \{ 0 \}$.

Furthermore, by setting: 

\begin{equation}\label{eq102}
r_k = (m_1 + 1) |\xi_k|^{-\frac{1}{m_1 +1}}, \text{ as } k \to +\infty
\end{equation}

and 

\begin{equation}\label{eq103}
v_k (z) = u_k (r_k z) + \frac{2}{a} \log (r_k)
\end{equation}

we shall show that $v_k$ blows up at exactly $(m_1+1)$-points, say $z_1,...,z_{m_1+1}$, which (after a rotation) correspond to the $(m_1+1)$-roots of the unity. Moreover, for the relative  $\beta_j$ defined in \eqref{eq515}, we have:  $\beta_j = \beta_0$, $\, \forall \, j=1,...,m_1 + 1$, with $\beta_0$ in \eqref{eq96}.

In other words, $v_k$ admits precisely the blow--up behavior described by Corollary \ref{cor74}, and accordingly the value of $\beta_0$ in \eqref{eq96} matches with that in \eqref{eq89} when $m = m_1 +1$ and $N$ in given by \eqref{eq94}.

To check that $u_k$ blows--up at the origin, it suffices to take $z_k \in \mathbb{C}$, such that: $\left( \frac{\bar{z}_k}{m_1 + 1}\right)^{m_1 + 1} = \frac{\xi_k}{|\xi_k|^2}$, so that $z_k \to 0$ and $u_k (z_k) = u_0 (0) + \beta_0 \log |\xi_k| \to +\infty$, as $k \to +\infty$.

Furthermore, by recalling that, for every $R>0$ there exists a suitable constant $C = C(R) > 0$:

\begin{equation}\label{eq104}
u_0 (\eta) + \beta_0 \log |\eta| \leq C, \quad \forall \; |\eta| \geq R
\end{equation}

 we see that, for every compact set $K \subset\subset \mathbb{C} \setminus \{ 0 \}$ there holds:
$$\sup_{K} u_k \to - \infty \quad \text{ as } k \to +\infty,$$
and we conclude that $S = \{ 0 \}$.

Concerning $v_k = v_k (z)$, in \eqref{eq102}-\eqref{eq103}, we find:
 
\begin{equation}\label{eq105}
\begin{split}
v_k (z) &= u_0 \left( |\xi_k| \left(\frac{m_1 + 1}{\bar{z}}\right)^{m_1 + 1} - \frac{\xi_k}{|\xi_k|} \right)\\
&+ \beta_0 \log \left( \frac{|\xi_k|}{|z|^{m_1 + 1}}\right) - \frac{2}{a(m_1 + 1)} \log |\xi_k| + \frac{2}{a} \log (m_1 + 1) = \\
&=u_0 \left( |\xi_k| \left(\frac{m_1 + 1}{\bar{z}}\right)^{m_1 + 1} - \frac{\xi_k}{|\xi_k|} \right) - \beta_0 (m_1 + 1) \log |z| + \frac{2}{a} \log (|\xi_k|(m_1 + 1))
\end{split}
\end{equation}

Therefore, by setting: $\frac{\xi_k}{|\xi_k|} = e^{i \theta_k}$, $\theta_k \in [0,2\pi)$, and assuming that: $\theta_k \to \theta_0 \in [0,2\pi)$, then we see that,  

$$z_{j,k} := e^{i \frac{\theta_k + 2 \pi (j-1)}{m_1 + 1}}\to e^{i \frac{\theta_0 + 2 \pi (j-1)}{m_1 + 1}}:= z_j, \quad j=1,...,m_1+1,$$ 

and,


$$v_k (z_{j,k}) = u_0 (0) + \frac{2}{a} \log (|\xi_k|(m_1 + 1)) \to + \infty,$$
as $ k \to +\infty$.

Moreover, we can use once more \eqref{eq104}, in order to check that, for every compact set $K \subset \mathbb{C} \setminus \{ z_1,...,z_{m_1 + 1} \}$ we have:
$$\sup_K v_k \to - \infty$$
as $k \to +\infty$.

Hence the blow--up set $S_1$ of $v_k$ is formed (up to rotation) by the $(m_1+1)$-roots of the unity and ''concentration'' occurs. Similarly, we may also check that $0 \notin \hat{S}_1$.

Finally, by using \eqref{eq105} and the natural change of variable: $z \to \left( z_{j,k} + \frac{\eta}{(m_1 + 1) |\xi_k|}\right)$ we find:

$$\frac{1}{2\pi} \int_{B_\delta (z_{j,k})} (e^{a v_k} + r_{k}^{\frac{2}{a} (a (N+1) -1)} |z|^{2N} e^{v_k} ) \, |dz| =$$
$$\frac{1}{2\pi} \int_{|\eta| \leq \delta |\xi_k|(m_1 +1)} (e^{a u_0} + e^{u_0} ) \, |d \eta| + o(1) \to \frac{1}{2\pi} \int_{\mathbb{R}^2} (e^{a u_0} + e^{u_0} )\, |d \eta| = \beta_0, \text{ as } k \to \infty,$$
and so we obtain that,

$$\beta_j = \lim_{\delta \to 0} \lim_{k \to \infty} \frac{1}{2\pi} \int_{B_\delta (z_{j})} (e^{a v_k} + r_{k}^{\frac{2}{a} (a (N+1) -1)} |z|^{2N} e^{v_k} ) = \beta_0, \; \forall \, j=1,...,m_1+1$$

as claimed.

\section{The case $a > 1$}\label{amagg1}

In view of part $i)$ of Corollary \ref{cor72} and Corollary \ref{coro415}, it remains to analyze the case where:

\begin{equation}\label{eq10.1}
a > 1 \; \text{ and } \; S = \emptyset
\end{equation}

To this purpose, for $\varepsilon > 0$ sufficiently small, we let $r_k > 0$ be defined in \eqref{eq76}. In view of \eqref{eq10.1}  we have:

\begin{equation}\label{eq10.2}
r_k \to +\infty.
\end{equation}

Therefore, in this case we consider the following new scaling for $u_k$:

\begin{equation}\label{eq10.3}
v_k (x) = u_k (r_k x) + 2(N+1) \log r_k
\end{equation}

which satisfies:

\begin{equation}\label{eq10.4}
\begin{cases}
-\Delta v_k=\varepsilon_{2,k}e^{av_k}+\left|x\right|^{2N}e^{v_k}:=f_{2,k}(x) \mbox{ in } \mathbb{R}^2\\
\beta_k:=\dfrac{1}{2\pi}\displaystyle{\int_{\RR^{2}}f_{2,k}(x)dx}.
\end{cases}
\end{equation}

with $\varepsilon_{2,k}= r_{k}^{2(1-a(N+1))},$ and so,

\begin{equation}\label{eq10.5}
\varepsilon_{2,k} \to 0, \; \text{ as } \; k \to \infty.
\end{equation}

By using the notation above, we let $\hat{v}_k$ as defined in \eqref{eq713} and denote by $S_1$ and  $ \hat{S}_1$  the (possibly empty) blow-up set of $v_k$ and $\hat{v}_k$ respectively.

Clearly problem \eqref{eq10.4}--\eqref{eq10.5} satisfies the assumption \eqref{conclu1}, and therefore  we may apply Proposition \ref{prop29} and its consequences to $v_k$ or  $\hat{v}_k$.
Indeed, exactly as for Lemma \ref{lem50}, we can use Proposition \ref{prop33}, Proposition \ref{prop33bis} together with Corollary \ref{coro31}, Corollary \ref{coro34bis} and Proposition \ref{prop35} in order to obtain that the analogous of the properties in \eqref{eq5.14.1}, \eqref{eq5.14.2}, \eqref{eq5.14.6} and \eqref{eq5.14.7} hold for $v_k$ and $\hat{v}_k$ respectively. In particular, we can check that the following holds:

\begin{lemma}\label{lem71bis}
Let $0<a\ne\frac{1}{N+1}$. If $x_0 \in S_1$ and 

\begin{equation*}
\begin{split}
&\beta(x_0) = \lim_{r \to 0} \lim_{k \to +\infty} \frac{1}{2\pi} \int_{B_r (x_0)} f_{2,k} (x) \, dx \\
&=\lim_{r\to0}\lim_{k\to+\infty} \left( \frac{1}{2\pi}\int_{B_{r} (x_0)} \varepsilon_{2,k} e^{a v_k(x)} \, dx + \frac{1}{2\pi}\int_{B_{r} (x_0)} |x|^{2N}  e^{a v_k(x)} \, dx \right)  := \beta_{0,1} (x_0) + \beta_{0,2} (x_0) 
\end{split}
\end{equation*}

then, by setting 

$$N(x_0) = \begin{cases} N \quad \text{ if } x_0 = 0\\
0 \quad \text{ if } x_0 \ne 0
\end{cases}$$

we have that $\beta_{0,1} (x_0)$ and $\beta_{0,2} (x_0)$ satisfy \eqref{eq5.14.1} with $N = N(x_0)$, and in particular, 

\begin{equation}\label{eq10.6}
\frac{4}{a} \leq \beta (x_0) \leq 4(N(x_0)+1).
\end{equation}

Furthermore,

\begin{itemize}
\item \begin{equation}\label{eq10.7}
\text{ if } a \leq 2 \text{ or } \beta (x_0) \geq 2  \text{ then } f_{2,k} \rightharpoonup \beta_0 \delta_{x_0} \text {weakly in the sense of measures in } B_\delta (x_{0})
\end{equation}

\item \begin{equation}\label{eq10.8}
\text{ if } \beta (x_0) < 2  \text{ then }\beta (x_0)=\frac{4}{a} \text{ and } a > 2 \text{ in this case.}
\end{equation} 

\end{itemize}

If $0 \in \hat{S}_1$ and $\beta > \max\{ 2(N+1), \frac{2}{a} \}$ then

\item \begin{equation}\label{eq10.9}
\beta_\infty := \lim_{R \to +\infty} \lim_{k \to +\infty} \frac{1}{2\pi} \int_{|x| \geq R} f_{2,k} (x) \, dx \geq 2\beta -  \max\{ 4(N+1), \frac{4}{a} \}
\end{equation} 

\end{lemma}
\qed

By applying the results above to the sequence $v_k$ in \eqref{eq10.3} we find:

\begin{prop}\label{prop7.2}
Assume \eqref{eq10.1}, we have:
\begin{itemize}

\item If $S_1 = \emptyset$ then $\beta = 4 (N+1)$;

\item If $S_1 \ne \emptyset$ then $S_1 = \{ z_1,..., z_m\}$ with $|z_j|\geq 1$ and $\beta_j := \beta (z_j) \in \left[ \frac{4}{a}, 4 \right] \; \forall \, j=1,...,m$.

\end{itemize}

Moreover,

\begin{enumerate}[i)]

\item if $N \in (0,1)$ and $0 < a \leq \frac{2}{N+1}$ then $\beta = \frac{4}{a}$ and $0 \notin \hat{S}_1$ for $a \ne \frac{2}{N+1}$ ;

\item if $a> \max\{ 1, \frac{2}{N+1} \}$ then \underline{either} $\beta$ satisfies \eqref{eq91.1},  \underline{or} $0 \notin \hat{S}_1$ and $\beta$ satisfies one of the following:\\
\begin{equation}\label{eq10.13}
 \beta = 2 (N+1+\sqrt{(N+1)^2-\frac{4m}{a^2}(a(N+1)-1)}),
\end{equation}

with $m \in \mathbb{N}$ and $1 \leq m \leq \left( \frac{(N+1)a}{2} \right)^2 \frac{1}{a(N+1)-1}$ and $a>2$;\\

\begin{equation}\label{eq10.14}
 \beta = \sum_{j=1}^{m} \beta_{j}, \; \max_{j=1,...,m} \beta_j \geq 2 \text{ and } \sum_{j=1}^{m} \beta_{j}^{2} = \frac{\beta (4N - (1-a)\beta)}{a(N+1)-1},
\end{equation}

with $m \in \mathbb{N}$ and $2 \leq m < \ a(N+1)- max\left\{ 0, \frac{a}{2}-1 \right\} $.\\

\end{enumerate}

\end{prop}

\dimo If $S_1 = \emptyset$ then (along a subsequence) we have:
$$
v_k\to v \mbox{ in } C^{2}_{loc}(\mathbb{R}^2),
$$
with $v$ satisfying: $-\Delta v=|x|^{2N} e^{v}$ in $\mathbb{R}^2$ and such that,  $|x|^{2N} e^{v} \in L^1 (\mathbb{R}^2)$. Consequently, $\beta \geq \frac{1}{2\pi} \int_{\mathbb{R}^2} |x|^{2N} e^{v}=4(N+1)$, and since $\beta \leq 4(N+1)$, we may conclude that, $\beta=4(N+1)$ as claimed.

If $S_1 \ne \emptyset$ then, by recalling that,

\begin{equation}\label{eq10.16}
\frac{1}{2\pi} \int_{|x|\leq1} |x|^{2N} e^{v} = \frac{4}{a} - \varepsilon
\end{equation}

we may use \eqref{eq10.6}  to conclude that if $x_0 \in S_1$ then $|x_0| \geq 1$ and $ \beta (x_0) \in \left[ \frac{4}{a}, 4 \right], \; \forall  j=1,...,m$. So $S_1 = \{ z_1,..., z_m\}$,  and if we assume that,  
\begin{equation}\label{eq10.15bis}
N \in (0,1)   \text { and }   0 < a \leq \frac{2}{N+1}, 
\end{equation}

then we can use  \eqref{eq10.7} in order to show that, 

\begin{equation}\label{eq10.16bis}
\frac{1}{2\pi} f_{2,k} \rightharpoonup \sum_{j=1}^{m} \beta_j \delta_{z_j}\text{ weakly in the sense of measure },
\end{equation}
locally in $\mathbb{R}^2$. We claim that,

\begin{equation}\label{eq10.19}
\lim_{R\to+\infty} \left(\lim_{k\to+\infty} \frac{1}{2\pi} \int_{|x|\geq R} f_{2,k} (x) \, dx \right) = 0.
\end{equation}
To check \eqref{eq10.19}, we observe first that, $0 \notin \hat{S}_1$ when $1<a<\frac{2}{N+1}$. Indeed in this case, $\beta \geq \frac{4}{a} > 2(N+1)$, and if otherwise then from \eqref{eq10.9} we would find,

\begin{equation}\label{eq10.18}
\beta_\infty \leq \beta - \frac{4}{a} \text{ and } \beta_\infty \geq 2 \beta - 4(N+1),
\end{equation}

that is: $\frac{4}{a} < 2(N+1) < \beta \leq 4(N+1) - \frac{4}{a}$, in contradiction with \eqref{eq10.15bis}. Similarly, we check that $0 \notin \hat{S}_1$  also in case, $a = \frac{2}{N+1}$ but $\beta > 2(N+1)$. Therefore,  \eqref{eq10.19} holds in such situations.

Finally, if $a = \frac{2}{N+1}$ and $\beta = 2(N+1) = \frac{4}{a}$,  then $m=1$ and $\beta = \beta_1$, and so  \eqref{eq10.19} is verified in this case as well.
As a consequence of  \eqref{eq10.19} and Lemma \ref{lem71bis}  we obtain that \eqref{eq83} holds together with the following:

\begin{equation}\label{eq10.17}
\beta = \sum_{j=1}^{m} \beta_{j}, \; \beta_j \in \left[ \frac{4}{a}, 4 \right], \; \forall \, j=1,...,m.
\end{equation}
On the other hand, from \eqref{eq10.17}  we obtain: $\frac{4}{a} m \leq \beta \leq 4(N+1)$, and by virtue of \eqref{eq10.15bis}, we find that necessarily,  $1 \leq m \leq 2$. We can rule out the possibility that $m = 2$, since it can be attained only for  $ a=\frac{2}{N+1}$,  $\beta_1 = \beta_2 = \frac{4}{a}$ and $\beta=4(N+1)$, but these values would make \eqref{eq83} fail. Thus,  $m = 1$ and by \eqref{eq83} we find  $\beta = \frac{4}{a}$ as claimed. \\
To establish $ii)$, we simply need to treat the case where: $\beta > 2 (N+1)$ and observe that, if in addition, we have that $0 \in \hat{S}_1$ then \eqref{eq10.18} holds, and it implies \eqref{eq91.1}.
Thus to complete the proof we need to show that, if  $0 \notin \hat{S}_1$ and $\beta > 2(N+1)$, then either one of  \eqref{eq10.13} and  \eqref{eq10.14} must hold. To this purpose, we observe first that in this case \eqref{eq10.19} holds.\\
We start to analyze what happens if $\max_{j=1,...,m} \beta_j \geq 2$.\\
As in the proof of Proposition \ref{prop71}, we see that in this case "concentration" occurs and  \eqref{eq83} holds, and it implies \eqref{eq10.14}. 
Since if we take $m=1$ in  \eqref{eq10.14} we obtain the value  $\beta = \frac{4}{a}$, while under the given assumption we have : $\beta \geq 2(N+1) > \frac{4}{a}$, we deduce that necessarily: $m \geq 2$. Also notice that, if  $\beta_j < 2$ for some $j \in \{ 1,...,m \}$ then necessarily $\beta_j = \frac{4}{a}$ and $a > 2$.Thus, when $\max \left\{1, \frac{2}{N+1} \right\} < a \leq 2$, then for \eqref{eq10.17} we derive the estimate: $\frac{4}{a}m \leq \beta \leq 4 (N+1)$, which implies that,  $2 \leq m \leq a(N+1)$.

While, for $a>2$, then from  \eqref{eq10.17} we obtain the estimate: $2 + \frac{4}{a} (m-1) \leq \beta \leq 4 (N+1)$, which implies a better estimate on $m$, namely:  $1 \leq (m-1) \leq \frac{a}{2} (2N+1)$. Thus, by combining such information, we check that  \eqref{eq10.14} holds in this case.

Finally, in case:  $\max_{j=1,...,m} \beta_j < 2$ then $a>2$ and $\beta_j = \frac{4}{a}$, for every $j=1,...,m$,  and in view of part $ii)$ of Lemma \ref{lem3024bis}, we see that ''concentration'' cannot occur.

Therefore,

$$\frac{4}{a} m = \lim_{k\to+\infty} \left(\varepsilon_{1,k} \frac{1}{2\pi} \int_{\mathbb{R}^2} e^{a v_k}\right) = \lim_{k\to+\infty} \left(\frac{1}{2\pi} \int_{\mathbb{R}^2} e^{a u_k}\right) = \frac{a \beta (4(N+1) - \beta)}{4(a(N+1)-1)}$$

and from this identity easily we derive \eqref{eq10.13}.

\qed

\begin{rem}\label{rem7.1}
As before, the possibility for $\beta$ to attain the value in \eqref{eq10.13} is linked to the solvability of the following singular Liouville equation:

\begin{equation}\label{eq10.20}
\begin{cases}
-\Delta v = |x|^{2N} e^{v} + \frac{8 \pi}{a} \sum_{j=1}^{m} \delta_{z_j} \mbox{ in } \mathbb{R}^2\\
\dfrac{1}{2\pi}\displaystyle{\int_{\RR^{2}} |x|^{2N} e^{v} \, dx} = \beta - \frac{4}{a} m.
\end{cases}
\end{equation}

with $\beta$ specified in \eqref{eq10.13}. 
\end{rem}

We know that, a solution of problem \eqref{eq10.20} can be interpreted as the  conformal factor for a metric on the Riemann sphere  with constant curvature equal to one and  assigned conical singularities. In spite of the results in \cite{er1}, \cite {egt1}, \cite {egt2}, \cite {egt3}, \cite{cwx}, \cite{troy1}, \cite {troy2} and \cite {lt}, which characterize the existence of such metric in some cases, still it remains a challenging open problem to determine the appropriate set of necessary and sufficient conditions about the angle at the singularities  so that  such  metric exists in general. \\ 

Notice also that the value of $\beta$ specified in \eqref{eq10.13} belongs in the interval $\left( 2(N+1), 4(N+1) - \frac{4}{a} \right]$, and actually it coincides with the value: $4(N+1) - \frac{4}{a}$, exactly when $m=1$. Therefore, it is reasonable to expect that \eqref{eq10.13} may occur only with $m=1$.

\section{Appendix}\label{appendix}

In this Appendix we provide the proof of some simple but technical Calculus lemmata, that we have used above.

\textbf{Calculus Lemma}: Let $N \geq 1$. The function $f(a) := \frac{1}{2a} \left( \sqrt{(1-a(N+1))^2 + \frac{Na}{1-a}}- (1 - a(N+1)) \right)$ satisfies the following:

\begin{enumerate}[(i)]
\item if $a \in \left( 0 , \frac{1}{N+1} \right)$ then $f'(a) > 0$ and $\frac{N}{4} < f(a) < \frac{N+1}{2}$;
\item $f(a) \geq 1$ if and only if $N>1$ and $a_N \leq a < \frac{1}{N+1}$, with $a_N$ defined in \eqref{eq532}.
\item If $1 < N \leq 3$ and $a \in \left( 0, \frac{1}{N+1} \right)$ then $f(a) < 2$.
\end{enumerate}

\underline{Proof}: It is straightforward to compute: 

$$f'(a) = \frac{1}{2a^2} \left(1-  \frac{(1-a(N+1)) + \frac{Na}{1-a} \left( 1 - \frac{1}{2(1-a)} \right)}{\sqrt{(1-a(N+1))^2 + \frac{Na}{1-a}}} \right)$$

Therefore $f'>0$ if and only if

$$(1-a(N+1)) + \frac{Na}{1-a} \left( 1 - \frac{1}{2(1-a)} \right) < \sqrt{(1-a(N+1))^2 + \frac{Na}{1-a}}$$

or equivalently,

\begin{equation}\label{eq533}
2 (1-a(N+1)) \left( 1 - \frac{1}{2(1-a)} \right) + \frac{Na}{1-a} \left( 1 - \frac{1}{2(1-a)} \right)^2 < 1
\end{equation}

To check \eqref{eq533} we use \eqref{eq529},  so that, for $0 < a < \frac{1}{N+1}$, we know that,

\begin{equation}\label{eq534}
\frac{Na}{1-a} < 1 - (1-a(N+1))^2
\end{equation}

Thus, we can insert \eqref{eq534} into \eqref{eq533} and derive the following estimate:

\begin{equation}\label{eq534bis}
\begin{split}
&\frac{Na}{1-a} \left( 1 - \frac{1}{2(1-a)} \right)^2 + 2 (1-a(N+1)) \left( 1 - \frac{1}{2(1-a)} \right) <\\
&\left( 1 - \frac{1}{2(1-a)} \right)^2 - (1-a(N+1))^2 \left( 1 - \frac{1}{2(1-a)} \right)^2 + 2 (1-a(N+1)) \left( 1 - \frac{1}{2(1-a)} \right) =\\
&= \left( 1 - \frac{1}{2(1-a)} \right)^2 - \left( (1-a(N+1)) \left( 1 - \frac{1}{2(1-a)} \right) - 1 \right)^2 + 1 =\\
&= 1 - \left(1+ \left( 1 - \frac{1}{2(1-a)}\right) a(N+1) \right) \left(\frac{a}{1-a} +  a(N+1) \left( 1 - \frac{1}{2(1-a)}\right) \right).
\end{split}
\end{equation}

Since, for $N \geq 1$ and $0 < a < \frac{1}{N+1}$ we have: $\frac{1}{2(1-a)} < 1$, we see that \eqref{eq533} follows from \eqref{eq534bis}.

Next, we observe that $f(0):=\lim_{a \to 0^+} f(a) = \frac{N}{4}$ and $f \left( \frac{1}{N+1} \right) = \frac{N+1}{2}$, and therefore for $a \in \left( 0, \frac{1}{N+1} \right)$ we have: $\frac{N}{4} < f(a) < \frac{N+1}{2}$.

Consequently, for $N \in (0,1]$ we see that necessarily: $f(a) < 1$, for all $a \in \left( 0, \frac{1}{N+1} \right)$. 

While, if $N \geq 4$ and $a \in \left( 0, \frac{1}{N+1} \right)$then $f(a) > \frac{N}{4} \geq 1$.

On the other hand, for $1 < N < 4$ we find that, $f(a) = 1$ if and only if $a = a_N = \frac{4-N}{2(N+1+\sqrt{(N-1)^2+N^2)}} \in \left( 0, \frac{1}{N+1} \right)$ and so in this case: $f(a) > 1$ if and only if $a_N < a < \frac{1}{N+1}$, as claimed in $(i)$.

Finally, for $1 \leq N \leq 3$ and $0 < a < \frac{1}{N+1}$ we have: $f(a) < \frac{N+1}{2} \leq 2$. \\
\qed

Next, we point out some monotonicity property for the functions related to the expression \eqref{eq527} and \eqref{eq527bis}. To this purpose, we consider

\begin{equation}\label{A.1}
\phi (t) = 1-t + \sqrt{(1-t)^2 + 2 a (N+1) t - \frac{Na}{1-a} t^2}
\end{equation}

which, for $t = \frac{2(m-1)(1-a)}{N}$, coincides with $\frac{a}{2} \beta_{j_0}$ and $\beta_{j_0}$ defined in \eqref{eq527}.

Similarly we consider

\begin{equation}\label{A.2}
\psi (t) = 1-s + \sqrt{(1-s)^2 + 2 a s + \frac{Na}{1-a(N+1)} s^2}
\end{equation}

and notice that, $\psi(t)$ for $s = \frac{2(m-1)(1-a(N+1))}{N}$ coincides with the expression of $\frac{a}{2} \beta$, with $\beta$ given by \eqref{eq527bis}.

Clearly, for $0 < a < \frac{1}{N+1}$, both functions $\phi$ and $\psi$ are always well defined. \\
We prove:

\textbf{Calculus Lemma 1}: Let $a \in \left( 0 , \frac{1}{N+1} \right)$ then the following holds:

\begin{enumerate}[(i)]

\item $\phi'(t) < 0$ and we have

\begin{equation}\label{A.3}
\phi(t) \geq 1 \Leftrightarrow t \leq  \frac{1-a}{aN} \left( \sqrt{(1-a(N+1))^2 +\frac{aN}{1-a}} - (1 - a (N+1))\right)=:t_a;
\end{equation}

\item if $0 < s \leq \frac{1-a(N+1)}{1-a} t_a=: s_a$ then $\psi'(s) < 0$, and in particular

\begin{equation}\label{A.4}
\psi(s) \geq \psi(s_a) = a(N+1) + \left( \sqrt{(1-a(N+1))^2 +\frac{Na}{1-a}}\right).
\end{equation}

\end{enumerate}

\underline{Proof}: By straightforward calculations, we obtain:

$$\phi'(t) = \frac{\left( \frac{t}{1-a} - 1 \right)(1-a(N+1)) - \sqrt{(1-t)^2 + 2 a (N+1) t - \frac{Na}{1-a} t^2}}{\sqrt{(1-t)^2 + 2 a (N+1) t - \frac{Na}{1-a} t^2}}$$

and readily we see that, $\phi'(t)<0$ for $t \leq 1-a$. On the other hand, if $t > 1-a$ then

{\footnotesize
\begin{equation}\label{A.5}
\phi'(t) < 0 \Leftrightarrow \left( \frac{t}{1-a} - 1 \right)^2 (1-a(N+1))^2 < \left( \frac{t}{1-a} - 1 \right)(1-a(N+1))^2 (1-a)+ a (N+2-a(N+1))
\end{equation}}

and so for $0<a<\frac{1}{N+1}$  the inequality in \eqref{A.5} certainly holds . 

Furthermore, by simple calculations we also find that $\phi(t) \geq 1$ if and only if $t \leq t_a$, with $t_a$ as given in \eqref{A.3}.

Concerning $(ii)$, we compute:

$$\psi'(s) = - 1 + \frac{(1-a) \left( \frac{s}{1-a(N+1)} - 1 \right)}{\sqrt{(1-s)^2 + 2 a s + \frac{Na}{1-a(N+1)} s^2}} = - 1 + \frac{(1-a) \left( \frac{s}{1-a(N+1)} - 1 \right)}{\sqrt{1-2(1-a)s + \frac{(1-a)s^2}{1-a(N+1)}}}$$

and after straightforward calculations we find that,  $\psi'(s) < 0$ if and only if :

\begin{equation}\label{A.8}
 \frac{s}{1-a(N+1)} < 1+\sqrt{1 + \frac{2-a}{N(1-a)}}.
\end{equation}
Therefore, to establish \eqref{A.4} we only need to check that,

\begin{equation}\label{A.8bis}
s_a \leq (1-a(N+1)) \left(1+\sqrt{1 + \frac{2-a}{N(1-a)}}\right).
\end{equation}

In other words, we must show that the following inequality:

\begin{equation}\label{A.9}
\frac{1}{Na} \left( \sqrt{(1-a(N+1))^2 +\frac{Na}{1-a}} - (1-a(N+1)) \right) < 1+\sqrt{1 + \frac{2-a}{N(1-a)}}
\end{equation}
 
is satisfied for any $a \in \left( 0, \frac{1}{N+1} \right)$.

But this is established easily once we observe that \eqref{A.9} can be written equivalently as follows:

\begin{equation}\label{A.10}
\frac{1}{Na} \left( \sqrt{(1-a(N+1))^2 +\frac{Na}{1-a}} - 1 \right) < \sqrt{1 + \frac{2-a}{N(1-a)}} - \frac{1}{N}
\end{equation}

and for $a \in \left( 0, \frac{1}{N+1} \right)$ and $N \geq 1$, we see that  the left hand side of \eqref{A.10} is negative while the right hand side is positive.
\qed

\begin{rem}\label{remA}
On the basis of Calculus Lemma 1, we check that for $\beta_{j_0} \geq \frac{2}{a}$ is satisfied if and only if \eqref{eq517} holds, and it implies that, $\beta \geq 2 (N+1) + \frac{2}{a} \sqrt{(1-a(N+1))^2 + \frac{Na}{1-a}}$ as claimed in Theorem \eqref{teo5.6}.
\end{rem}


\begin{thebibliography}{99}



\bibitem{ao1} J. Ambjorn, P. Olesen,
{\it Anti--screening of large magnetic fields by vector bosons}, Phys. Lett., {\bf{B 214}}, 565--569. (1988).

\bibitem{ban} C. Bandle, 
{\it Isoperimetric Inequalities and Applications}, Pitman, Boston, (1980).

\bibitem{BaCa}  D. Bartolucci, D. Castorina,  {\it Self gravitating cosmic strings and the Alexandrov's inequality for Liouville-type equations}, preprint (2014), arXiv:1409.3135.

\bibitem{barm} D. Bartolucci, E. Montefusco, 
{\it Blow--up analysis, existence and qualitative properties of solutions for the two dimensional Emden--Fowler equation with singular potential}, Math. Meth. in Appl. Sci., {\bf{30}}(18), (2007), 706--729.

\bibitem{blt} D. Bartolucci, C. S. Lin, G. Tarantello,
{\it Uniqueness and symmetry results for solutions of a mean field equation on $S^{2}$ via a new bubbling phenomenon}, Comm. Pure Appl. Math., {\bf{64}} (12), (2011), 1677--1730.

\bibitem{bt} D. Bartolucci, G. Tarantello, {\sl Liouville type equations with singular data and their applications to periodic multivortices for the
electroweak theory}, Comm. Math. Phys., {\bf 229} (2002), 3--47.

\bibitem{bm} H. Brezis, F. Merle,
{\it Uniform estimates and blow--up behavior for solutions of $-\Delta u =V(x)e^{u}$ in two dimensions}, Comm. Partial Differential Equations {\bf{16}}, (1991), 1223--1253.

\bibitem{bp} S. Baraket, F. Pacard,
{\it Construction of a singular limit for a semilinear elliptic equation in dimension $2$}, Calc. Var. P.D.E., {\bf{6}}, (1998), 1--38.

\bibitem{c} D. Chae, 
{\it Existence of a semilinear elliptic system with exponential nonlinearities}, Discrete Contin. Dyn. Syst., {\bf{18}}, (2007), 709--718.

\bibitem{c1} D. Chae, 
{\it Existence of multistring solutions of a self--gravitating massive W--boson}, Lett. Math. Phys., {\bf{73}}, (2005), 123--134.

\bibitem{cgs} R. M. Chen, Y. Guo, and D. Spirn,
{\it Asymptotic behaviour and symmetry of condensate solutions in electroweak theory}, J. Anal. Math., {\bf{117}}, 47--85.

\bibitem{ck} S. Chanillo, M. Kiessling, {\it Conformally invariant systems of nonlinear PDE of Liouville type}, Geom. Funct. Analysis, {\bf{5}} (1995), 924--947.
 
\bibitem{ck1} S. Chanillo, M. Kiessling, {\it Surfaces with prescribed scalar curvature}, Duke Math. J., {\bf{105}} (2002), 309--353. 

\bibitem{ck2} S. Chanillo, M. Kiessling, {\it Rotational symmetry of solutions of some nonlinear problems in statistical mechanics and in geometry}, Comm. Math. Phys., {\bf{160}} (1994), 217--238. 

\bibitem{ct} D. Chae, G. Tarantello, 
{\it Selfgravitating electroweak strings}, J. Differ. Equations, {\bf{213}}, (2005), 146--170.

\bibitem{ct1} D. Chae, G. Tarantello, 
{\it On planar selfdual electroweak vortices}, Ann. Inst. H. Poincare Anal. Non Lineaire, {\bf{21}}, (2004), 187--207.

\bibitem{chenlin1} C.C. Chen, C.S. Lin, 
{\it Sharp estimates for solutions of multi--bubbles in compact Riemann surfaces}, Comm. Pure Appl. Math., {\bf{55}}, (2002), 728--771.

\bibitem{chenlin2} C.C. Chen, C.S. Lin,
{\it Topological degree for a mean field equation on Riemann surfaces}, Comm. Pure Appl. Math., {\bf{56}}, (2003), no.12, 1667--1727.

\bibitem{chenlin3} C.C. Chen, C.S. Lin,
{\it Mean field equations of Liouville type with singular data: sharper estimates}, Discr. Cont. Dyn. Sist., {\bf{28}}, (2010), no.3, 1237--1272.

\bibitem{cwx} Q. Chen, Y. Wu and B. Xu, 
{\it Conformal metrics with constant curvature one and finite conic singularities on compact Riemann surfaces}, Pacific J. Math., {\bf 273} (2015), n.1,  75--100.

\bibitem{chlin} K.S. Cheng, C.S. Lin,
{\it On the asymptotic behavior of the conformal Gaussian curvature equations in $\in\RR^{2}$}, Math. Ann., {\bf{308}}(1), (1997), 119--139.

\bibitem{cl} W. Chen, C. Li,
{\it Qualitative properties of solutions to some nonlinear elliptic equations in $\RR^{2}$}, Duke Math. J., {\bf{71}}, (1993), no.2, 427--439.

\bibitem{cl1} W. Chen, C. Li,
{\it Classification of solutions of some nonlinear elliptic equations}, Duke Math. J., {\bf{63}}, (1991), no.2, 615--623.

\bibitem{csw} M. Chipot, I. Shafrir and G. Wolansky,
{\it On the solutions of the Liouville systems}, J. Differ. Equations, {\bf{140}}, (1997),  no.1, 59--105.

\bibitem {sg} H. P. de Saint-Gervais, 
{\it Uniformisation des surfaces de Riemann}, ENS Editions, 2010.

\bibitem{pkm} M. del Pino, M. Kowalczyk, M. Musso,
{\it Concentrating solutions in a two--dimensional elliptic problem with exponential Neumann data}, J. Funct. Anal., {\bf{227}}, (2005), 430--490.


\bibitem{egp} P. Esposito, M. Grossi, A. Pistoia, 
{\it On the existence of blowing--up solutions for a mean field equation}, Ann. I.H.P. Anal. Non Lineare, {\bf{22}}, n.2, (2005), 227--257.

\bibitem{er1} A. Eremenko,
{\it Metrics of positive curvature with conical singularities on the sphere}, Proc. AMS, {\bf {132}} (2004), 3349--3355.

\bibitem{egt1} A. Eremenko, A. Gabrielov and V. Tarasov,
{\it Metrics with conic singularities and spherical polygons}, arxiv:1405.1738

\bibitem {egt2} A. Eremenko, A. Gabrielov and V. Tarasov, 
{\it Metrics with four conic singularities and spherical quadrilaterals}, arxiv: arXiv:1409.1529.

\bibitem {egt3} A. Eremenko, A. Gabrielov and V. Tarasov, 
{\it Spherical quadrilaterals with three non-integer angles}, preprint 2015.

\bibitem{f} R. Fortini, 
{\it The role of Liouville type systems in the analysis of selfdual
gauge field theories}, PhD thesis Dipartimento di Matematica Universit\'a di Roma Tor Vergata (2013).

\bibitem{ft} R. Fortini, G. Tarantello,
{\it The role of Liouville systems in the study of non--abelian Chern--Simons vortices}, Proceeding ICMP12, XVIIth International Congress on Mathematical Physics, Aalborg (2012), 383--390.

\bibitem {fuj} S. Fujimori, Y. Kawakami, M. Kokubu, W. Rossman, M. Umehara and K. Yamada, 
{\it CMC-1 trinoids in hyperbolic 3-space and metrics of constant curvature one with conic singularities on the 2-sphere,} Proc. Japan
Acad., {\bf{87}} (2011), 144--149.

\bibitem{jost} J. Jost, G. Wang,
{\it Analytic aspects of the Toda system: I. A Moser--Trudinger inequality, } Comm. Pure Appl. Math., {\bf{54}}, (2001), 1289--1319.

\bibitem{ls} Y.Y. Li \& I.Shafrir, {\sl Blow-up analysis for Solutions of $-\Delta u = V(x)e^{u}$
in dimension two}, {Ind. Univ. Math. J.},  {\bf {43}}(4) (1994), 1255--1270.

\bibitem{lin} C.S. Lin,
{\it Uniqueness of solutionss to the mean field equation for the spherical Onsager vortex}, Arch. Ration. Mech. Anal., {\bf{153}}, (2000), 153--176.

\bibitem{lwy} C. S. Lin, J. Wei and D. Ye,
{\it Classification and nondegeneracy of $SU(n+1)$ Toda system with singular sources}, Invent. Math., {\bf{190}}, issue~1, (2012), pp. 169--207. 

\bibitem{lz} C. S. Lin, L. Zhang,
{\it Profile of bubbling solutions to a Liouville system}, Ann. Inst. Henri Poincar'e, Anal. Nonlin'eaire {\bf{27}}, (2010), no.1, 117--143.

\bibitem {lt} F. Luo and G. Tian, 
{\it Liouville equation and spherical convex polytopes}, Proc. Amer. Math. Soc. {\bf{116}} (1992), no. 4, 1119--1129.

\bibitem {mco} R. McOwen, 
{\it Point singularities and conformal metrics on Riemann surfaces}, Proc. Amer. Math. Soc. {\bf{103}} (1988), 222--224.


\bibitem{pt} J. Prajapat, G. Tarantello,
{\it On a class of elliptic problems in $\RR^{2}$: Symmetry and Uniqueness results}, Proc. Roy. Soc. Edinburgh, {\bf{131A}}, (2001), 967--985.

\bibitem{pot} A. Poliakovsky, G. Tarantello,
{\it On a planar Liouville--type problem in the study of selfgravitating strings}, J. Diff. Equations, {\bf{252}}, (2012),  no.5, 3668--3693.

\bibitem{pt2} A. Poliakovsky, G. Tarantello,
{\it On singular Liouville systems}, Analysis and Topology in Nonlinear Differential Equations, PNDLE 85 Birkhauser Basel, (2014).

\bibitem{sw1} I. Shafrir, G. Wolansky,
{\it Moser--Trudinger type inequalities for systems in two dimensions}, C.R. Math Acad. Sci. Paris, {\bf{333}}, (2001), 439--443.

\bibitem{sw2} I. Shafrir, G. Wolansky,
{\it Moser--Trudinger and logarithmic HLS inequalities for systems},J. Eur. Math. Soc., {\bf{7}}, (2005), 413--448.

\bibitem{s} T. Suzuki,
{\it Global analysis for a two--dimensional elliptic eigenvalue problem with the exponential nonlinearity}, Ann. Inst. Henri Poincar'e, Anal. Nonlin'eaire, {\bf{9}} (1992), 367--398.

\bibitem{tar} G. Tarantello, {\it Analytical aspects of Liouville-type equations with singular sources},  Stationary
Partial Differential Equations and their Applications  72. (2007) Boston: Birkh\"auser.

\bibitem{tar2} G. Tarantello, {\it Self--dual Gauge field vortices: an analytical approach}, PNDLE, 72 Birkhauser Boston, Inc. Boston MA, (2008).

\bibitem{tar?} G. Tarantello, {\it Radial symmetry for solutions of a cosmic strings equation}, work in progress.

\bibitem {troy1} M. Troyanov, 
{\it Metrics of constant curvature on a sphere with two conial singularities}, LNM 1410, Springer, NY, 1989, 296-308.

\bibitem {troy2} M. Troyanov, 
{\it Prescribing curvature on compact surfaces with conial singularities}, Trans. AMS, {\bf{324}} (1991),
793-821.

\bibitem {uy} M. Umehara and K. Yamada, 
{\it Metrics of constant curvature 1 with three conial singularities on the 2-sphere}, Illinois J. Math., {\bf{44}} (2000), 72-94.

\bibitem{y} Y. Yang, 
{\it Solitons in field theory and nonlinear analysis},  
Springer Monographs in Mathematics, Spinger--Verlag, New York (2001).

\end{thebibliography}
\end{document}